\documentclass[11pt]{article}

\usepackage{amsmath,amsthm,amsfonts,amssymb,mathtools,bm}
\usepackage{mathrsfs,bbm,nicefrac,stmaryrd}

\usepackage{xcolor,verbatim,xparse,etoolbox,cite,xurl}
\usepackage[utf8]{inputenc}

\usepackage[margin=1in]{geometry}
\usepackage{titlesec,setspace}
\usepackage[shortlabels]{enumitem}

\titleformat{\section}{\bfseries\large}{\thesection.}{1em}{}
\titleformat{\subsection}{\bfseries\normalsize}{\thesubsection.}{1em}{}

\usepackage[colorlinks,
    linkcolor={red!80!black},
    citecolor={green},
    urlcolor={blue!80!black},
    hypertexnames=false
]{hyperref}

\usepackage[nameinlink,capitalise,sort]{cleveref}

\crefname{enumi}{item}{items}
\crefname{equation}{}{}
\crefname{subsection}{Subsection}{Subsections}

\usepackage{thmtools}

\declaretheorem[style=plain,numberwithin=section,name=Theorem]{theorem}

\declaretheorem[style=plain,sibling=theorem,name=Proposition]{prop}

\declaretheorem[style=plain,sibling=theorem,name=Setting]{setting}

\declaretheorem[style=definition,sibling=theorem,name=Definition]{definition}

\crefname{theorem}{Theorem}{Theorems}
\Crefname{theorem}{Theorem}{Theorems}

\crefname{lemma}{Lemma}{Lemmas}
\Crefname{lemma}{Lemma}{Lemmas}

\crefname{prop}{Proposition}{Propositions}
\Crefname{prop}{Proposition}{Propositions}

\crefname{cor}{Corollary}{Corollaries}
\Crefname{cor}{Corollary}{Corollaries}

\crefname{definition}{Definition}{Definitions}
\Crefname{definition}{Definition}{Definitions}

\crefname{remark}{Remark}{Remarks}
\Crefname{remark}{Remark}{Remarks}

\crefname{setting}{Setting}{Settings}
\Crefname{setting}{Setting}{Settings}

\crefname{conjecture}{Conjecture}{Conjectures}
\Crefname{conjecture}{Conjecture}{Conjectures}

\crefname{figure}{Figure}{Figures}
\Crefname{figure}{Figure}{Figures}

\DeclareMathAlphabet{\mathpzc}{OT1}{pzc}{m}{it}

\DeclareFontEncoding{LS1}{}{}
\DeclareFontSubstitution{LS1}{stix}{m}{n}
\DeclareMathAlphabet{\mathscr}{LS1}{stixscr}{m}{n}
\newcommand{\imag}{\mathbf{i}}

\newcommand{\E}{\mathbb{E}}
\renewcommand{\P}{\mathbb{P}}

\newcommand{\R}{\mathbb{R}}
\newcommand{\N}{\mathbb{N}}
\newcommand{\Z}{\mathbb{Z}}

\newcommand{\bbO}{\mathbb{O}}

\newcommand{\Reli}[4]{\mathbf{N}^{#3}_{#2}}
\newcommand{\Relii}[5]{\mathbf{N}^{#3}_{#2,#5}}

\renewcommand{\d}{ \mathrm{d}}

\renewcommand{\c}[1]{\mathfrak{c}^{#1}}

\newcommand{\g}{\mathscr{g}}
\newcommand{\vertiii}[1]{{\left\vert\kern-0.25ex\left\vert\kern-0.25ex\left\vert #1 
    \right\vert\kern-0.25ex\right\vert\kern-0.25ex\right\vert}}

\newcommand{\scrl}{\mathscr{l}}

\newcommand{\cF}{\mathcal{F}}
\newcommand{\cG}{\mathcal{G}}

\newcommand{\cS}{\mathcal{S}}

\newcommand{\bbF}{\mathbb{F}}
\newcommand{\polycons}{q}

\newcommand{\bfa}{\mathbf{a}}

\newcommand{\bfd}{\mathbf{d}}
\newcommand{\bfe}{\mathbf{e}}

\newcommand{\bfl}{\mathbf{l}}

\newcommand{\bbA}{\mathbb{A}}

\newcommand{\bfF}{\mathbf{F}}

\newcommand{\bfS}{\mathbf{S}}

\newcommand{\bfW}{\mathbf{W}}
\newcommand{\bfX}{\mathbf{X}}

\newcommand{\scrA}{\mathscr{A}}

\newcommand{\fC}{\mathfrak{C}}

\newcommand{\fF}{\mathfrak{F}}

\newcommand{\fJ}{\mathfrak{J}}

\newcommand{\fL}{\mathfrak{L}}

\newcommand{\fd}{\mathfrak{d}}

\newcommand{\fl}{\mathfrak{l}}

\newcommand{\scrc}{\mathscr{c}}

\renewcommand{\emptyset}{\varnothing}

\newcommand{\Act}{\mathbb{A}}

\newcommand{\bbS}{\mathbb{S}}
\newcommand{\bbT}{\mathbb{T}}
\newcommand{\bbX}{\mathbb{X}}

\newcommand{\ceil}[1]{ \left\lceil #1 \right\rceil}

\newcommand{\qqandqq}{\qquad\text{and}\qquad}

\newcommand{\indicator}[1]{\mathbbm{1}_{\smash{#1}}}

\usepackage{cleveref}

\ExplSyntaxOn

\seq_new:N \g_cflist_loaded
\seq_new:N \g_cflist_pending

\NewDocumentCommand{\cfadd} { m } {
  \seq_if_in:NnF \g_cflist_loaded { #1 } {
    \seq_if_in:NnF \g_cflist_pending { #1 } {
      \seq_gput_right:Nn \g_cflist_pending { #1 }
    }
  }
}

\NewDocumentCommand{\cfconsiderloaded} { m } {
  \seq_gput_right:Nn \g_cflist_loaded {#1}
}

\NewDocumentCommand{\cfremove} { m } {
  \seq_gremove_all:Nn \g_cflist_pending { #1 }
}

\NewDocumentCommand{\cfload} { o } {
  \seq_if_empty:NTF \g_cflist_pending {
    \IfValueTF{#1}{\ignorespaces}{\unskip}
  } {
    (cf.\ \cref{\seq_use:Nn \g_cflist_pending {,}})\IfValueTF{#1}{#1~}{\unskip}
    \seq_gconcat:NNN \g_cflist_loaded \g_cflist_loaded \g_cflist_pending
    \seq_gclear:N \g_cflist_pending
    \IfValueT{#1}{\ignorespaces}
  }
}

\NewDocumentCommand{\cfclear} {} {
  \seq_gclear:N \g_cflist_loaded
  \seq_gclear:N \g_cflist_pending
}

\NewDocumentCommand{\cfout} { o } {
  \seq_if_empty:NTF \g_cflist_pending {\unskip\IfValueT{#1}{\ignorespaces}} {
    (cf.\ \cref{\seq_use:Nn \g_cflist_pending {,}})\IfValueTF{#1}{#1~}{\unskip}
    \seq_gclear:N \g_cflist_pending
    \IfValueT{#1}{\ignorespaces}
  }
}

\NewDocumentCommand{\ifnocf} { m } {
  \seq_if_empty:NT \g_cflist_pending { #1 }
}

\ExplSyntaxOff

\ExplSyntaxOn

\bool_new:N \g_noteobserve

\NewDocumentCommand{\setnote}{}{
  \bool_gset_true:N \g_noteobserve
}

\NewDocumentCommand{\setobserve}{}{
  \bool_gset_false:N \g_noteobserve
}

\NewDocumentCommand{\nobs}{ o }{
  \IfValueT{#1}{
    \str_if_eq:noTF {note} {#1} {
      \bool_gset_true:N \g_noteobserve
    } {
      \str_if_eq:noTF {Note} {#1} {
        \bool_gset_true:N \g_noteobserve
      } {
        \bool_gset_false:N \g_noteobserve
      }
    }
  }
  \bool_if:nTF { \g_noteobserve } {
    \bool_gset_false:N \g_noteobserve
    note
  } {
    \bool_gset_true:N \g_noteobserve
    observe
  }
  \IfValueF{#1}{~}
}

\NewDocumentCommand{\Nobs}{ o }{
  \IfValueT{#1}{
    \str_if_eq:noTF {note} {#1} {
      \bool_gset_true:N \g_noteobserve
    } {
      \str_if_eq:noTF {Note} {#1} {
        \bool_gset_true:N \g_noteobserve
      } {
        \bool_gset_false:N \g_noteobserve
      }
    }
  }
  \bool_if:nTF { \g_noteobserve } {
    \bool_gset_false:N \g_noteobserve
    Note
  } {
    \bool_gset_true:N \g_noteobserve
    Observe
  }
  \IfValueF{#1}{~}
}

\ExplSyntaxOff

\ExplSyntaxOn

\bool_new:N \g_hencetherefore

\NewDocumentCommand{\hence}{ o }{
  \IfValueT{#1}{
    \str_if_eq:noTF {hence} {#1} {
      \bool_gset_true:N \g_hencetherefore
    } {
      \str_if_eq:noTF {Hence} {#1} {
        \bool_gset_true:N \g_hencetherefore
      } {
        \bool_gset_false:N \g_hencetherefore
      }
    }
  }
  \bool_if:nTF { \g_hencetherefore } {
    \bool_gset_false:N \g_hencetherefore
    hence
  } {
    \bool_gset_true:N \g_hencetherefore
    therefore
  }
  \IfValueF{#1}{~}
}

\NewDocumentCommand{\Hence}{ o }{
  \IfValueT{#1}{
    \str_if_eq:noTF {hence} {#1} {
      \bool_gset_true:N \g_hencetherefore
    } {
      \str_if_eq:noTF {Hence} {#1} {
        \bool_gset_true:N \g_hencetherefore
      } {
        \bool_gset_false:N \g_hencetherefore
      }
    }
  }
  \bool_if:nTF { \g_hencetherefore } {
    \bool_gset_false:N \g_hencetherefore
    Hence,~we~obtain
  } {
    \bool_gset_true:N \g_hencetherefore
    Therefore,~we~obtain
  }
  \IfValueF{#1}{~}
}

\ExplSyntaxOff

\ExplSyntaxOn

\seq_const_from_clist:Nn \g_prove_mru {
  establish,
  demonstrate,
  prove,
  show,
  imply,
  ensure
}

\prop_new:N \l__verbs
\prop_put:Nnn \l__verbs {show} {shows}
\prop_put:Nnn \l__verbs {imply} {implies}
\prop_put:Nnn \l__verbs {demonstrate} {demonstrates}
\prop_put:Nnn \l__verbs {prove} {proves}
\prop_put:Nnn \l__verbs {establish} {establishes}
\prop_put:Nnn \l__verbs {ensure} {ensures}
\prop_put:Nnn \l__verbs {assure} {assures}

\tl_new:N \g_wordtmp
\seq_new:N \l_mytmps

\cs_generate_variant:Nn \str_if_in:nnTF { nVTF }
\cs_generate_variant:Nn \str_if_in:nnTF { xVTF }

\NewDocumentCommand{\prove}{ o }{
  \IfValueTF{#1}{
    \seq_clear:N \l_mytmps
    \seq_map_inline:Nn \g_prove_mru {
      \str_if_eq:nnTF {##1} {ensure} {
        \str_set:Nn \l_temps {n}
      } {
        \str_set:Nx \l_temps {\str_head_ignore_spaces:n {##1}}
      }
      \str_if_in:xVTF {#1} \l_temps {
        \seq_put_right:Nn \l_mytmps {##1}
      } { }
    }
    \seq_get_right:NN \l_mytmps \g_wordtmp
  } {
    \seq_get_right:NN \g_prove_mru \g_wordtmp
  }
  \tl_use:N \g_wordtmp
  \IfValueTF{#1}{}{~}
  \seq_gput_left:NV \g_prove_mru \g_wordtmp
  \seq_gremove_duplicates:N \g_prove_mru
}

\NewDocumentCommand{\proves}{ o }{
  \IfValueTF{#1}{
    \seq_clear:N \l_mytmps
    \seq_map_inline:Nn \g_prove_mru {
      \str_if_eq:nnTF {##1} {ensure} {
        \str_set:Nn \l_temps {n}
      } {
        \str_set:Nx \l_temps {\str_head_ignore_spaces:n {##1}}
      }
      \str_if_in:xVTF {#1} \l_temps {
        \seq_put_right:Nn \l_mytmps {##1}
      } { }
    }
    \seq_get_right:NN \l_mytmps \g_wordtmp
  } {
    \seq_get_right:NN \g_prove_mru \g_wordtmp
  }
  \str_set:NV \l_tmpa_str \g_wordtmp
  \prop_get:NVN \l__verbs \l_tmpa_str \l_tmpa_tl
  \tl_use:N \l_tmpa_tl
  \IfValueTF{#1}{}{~}
  \seq_gput_left:NV \g_prove_mru \g_wordtmp
  \seq_gremove_duplicates:N \g_prove_mru
}

\newcommand{\llabel}[1]{\savelabel{#1}\label{\loc.#1}\ignorespaces}

\clist_new:N \l_localreflist
\clist_new:N \l_reflist

\NewDocumentCommand{\lref} { m } {
  \clist_set:No \l_localreflist {#1}
  \clist_clear:N \l_reflist
  \clist_map_inline:Nn \l_localreflist { \clist_put_right:Nn \l_reflist {\loc.##1} }
  \cref{\l_reflist}
}

\NewDocumentCommand{\Lref} { m } {
  \clist_set:No \l_localreflist {#1}
  \clist_clear:N \l_reflist
  \clist_map_inline:Nn \l_localreflist { \clist_put_right:Nn \l_reflist {\loc.##1} }
  \Cref{\l_reflist}
}

\NewDocumentCommand{\itref}{ m m }{
  \clist_set:No \l_localreflist {#2}
  \clist_clear:N \l_reflist
  \clist_map_inline:Nn \l_localreflist { \clist_put_right:Nn \l_reflist {#1.##1} }
  \cref{\l_reflist}~in~\cref{#1}
}

\seq_new:N \l_enum_seq
\int_new:N \l_num_items

\bool_new:N \g_commaused_bool

\providecommand{\comma}{}

\cs_new:Nn \enum_it:nn {
  \int_case:nnF {\l_num_items - #1} {
    {0} {
      \renewcommand{\comma}{}
      #2\space
    }
    {1} {
      \bool_gset_false:N \g_commaused_bool
      \renewcommand{\comma}{,~\bool_gset_true:N \g_commaused_bool}
      #2
      \bool_if:NTF \g_commaused_bool {} {,~}
      and~
    }
  } {
    \bool_gset_false:N \g_commaused_bool
    \renewcommand{\comma}{,~\bool_gset_true:N \g_commaused_bool}
    #2
    \bool_if:NTF \g_commaused_bool {} {,~}
  }
}

\cs_new:Nn \enum_it_U:nn {
  \int_case:nnF {\l_num_items - #1} {
    {0} {
      \renewcommand{\comma}{}
      #2
      \space
    }
    {1} {
      \bool_gset_false:N \g_commaused_bool
      \renewcommand{\comma}{,~\bool_gset_true:N \g_commaused_bool}
      #2
      \bool_if:NTF \g_commaused_bool {} {,~}
      and~
    }
  } {
    \bool_gset_false:N \g_commaused_bool
    \renewcommand{\comma}{,~\bool_gset_true:N \g_commaused_bool}
    \int_compare:nTF {#1=1} {
      \text_titlecase_first:n {#2}
    } {
      #2
    }
    \bool_if:NTF \g_commaused_bool {} {,~}
  }
}

\cs_generate_variant:Nn \tl_if_eq:nnTF {onTF}

\cs_new:Nn \enum:nnnn {
  \seq_set_split:Nnn \l_enum_seq ; {#1}
  \seq_remove_all:Nn \l_enum_seq { }
  \seq_remove_all:Nn \l_enum_seq {#2}
  \seq_log:N \l_enum_seq
  \int_set:Nn \l_num_items {\seq_count:N \l_enum_seq}
  \int_log:N \l_num_items
  \int_case:nnF {\l_num_items} {
    { 0 } { 0 }
    { 1 } {
      \IfBooleanTF{#4} {
        \tl_set:Nn \l_text_case_exclude_arg_tl {\cref}
        \bool_set_false:N \l_text_titlecase_check_letter_bool
        \text_titlecase_first:n {\seq_use:Nn \l_enum_seq {}}
      } {
        \seq_use:Nn \l_enum_seq {}
      }
      \space
      \tl_if_eq:onTF{#3}{-}{}{
        \bool_if:NTF \l_plural_bool {
          \prove[#3]~
        } {
          \proves[#3]~
        }
      }
    }
    { 2 } {
      \IfBooleanTF{#4} {
        \tl_set:Nn \l_text_case_exclude_arg_tl {\cref}
        \bool_set_false:N \l_text_titlecase_check_letter_bool
        \text_titlecase_first:n {\seq_item:Nn \l_enum_seq {1}}
        {} ~and~
        \seq_item:Nn \l_enum_seq {2}
      } {
        \seq_use:Nn \l_enum_seq {~and~}
      }
      \space
      \tl_if_eq:onTF{#3}{-}{}{
        \prove[#3]~
      }
    }
  } {
    \IfBooleanTF{#4} {
      \tl_set:Nn \l_text_case_exclude_arg_tl {\cref}
      \seq_indexed_map_function:NN \l_enum_seq \enum_it_U:nn
    } {
      \seq_indexed_map_function:NN \l_enum_seq \enum_it:nn
    }
    \tl_if_eq:onTF{#3}{-}{}{
      \prove[#3]~
    }
  }
}

\cs_generate_variant:Nn \enum:nnnn {nxnn}
\cs_generate_variant:Nn \enum:nnnn {nxxn}

\NewDocumentCommand{\enum}{O{} m O{-} s}{
  \IfBooleanTF{#4}{
    \enum:nxnn {#2} {#1} {sindep} \BooleanFalse
  } {
    \enum:nxxn {#2} {#1} {#3} \BooleanFalse
  }
}

\NewDocumentCommand{\dott}{}{\ifnocf{.}\space}

\bool_new:N \g_arg_start_bool
\bool_gset_true:N \g_arg_start_bool

\NewDocumentCommand{\startnewargseq}{}{\bool_gset_true:N \g_arg_start_bool \tl_set:Nn \g_label_tl {}}

\cs_generate_variant:Nn \seq_if_in:NnTF {NxTF}
\cs_generate_variant:Nn \seq_remove_all:Nn {Nx}

\int_new:N \l_random_int

\cs_generate_variant:Nn \tl_if_head_eq_catcode:nNTF {oNTF}
\cs_generate_variant:Nn \tl_if_head_eq_catcode:nNTF {VNTF}

\cs_generate_variant:Nn \tl_if_head_eq_catcode:nNTF {eNTF}
\cs_generate_variant:Nn \tl_log:n {o}
\cs_generate_variant:Nn \tl_log:n {f}
\cs_generate_variant:Nn \tl_log:n {x}
\cs_generate_variant:Nn \tl_log:n {e}

\cs_generate_variant:Nn \tl_if_in:nnTF {onTF}
\cs_generate_variant:Nn \tl_if_in:NnTF {NeTF}

\cs_generate_variant:Nn \tl_if_head_eq_meaning:nNTF {VNTF}

\seq_const_from_clist:Nn \g_arg_mru_this {
  Ahpr,
  Tapr,
  Ctapr,
  H
}

\seq_const_from_clist:Nn \g_arg_mru_nothis {
  Ia,
  Nwc,
  N,
  Itns,
  Fm,
  Itnswc,
  Mo
}

\seq_new:N \l_arg_seq
\tl_new:N \l_cons_tl
\tl_new:N \l_dummy_tl

\bool_new:N \g_debug_bool
\bool_gset_false:N \g_debug_bool

\bool_new:N \l_insidearg_bool

\bool_new:N \g_firstargletter_bool

\sys_gset_rand_seed:n {0903}

\bool_new:N \l_plural_bool
\tl_new:N \l_arg_verbs_tl

\NewDocumentCommand{\argument}{mom}{
\color{black}
  \bool_set_false:N \l_plural_bool
  \tl_set:Nn \l_arg_verbs_tl {sindep}
  \keys_define:nn { benno/argument } {
    plural .value_forbidden:n = true,
    plural .code:n = {\bool_set_true:N \l_plural_bool},
    verbs .value_required:n = false,
    verbs .tl_set:N = \l_arg_verbs_tl,
  }
  \IfValueT{#2}{
    \keys_set:nn { benno/argument } {#2}
  }
  \bool_log:N \l_plural_bool
  \bool_gset_true:N \l_insidearg_bool
  \seq_set_split:Nnn \l_arg_seq ; {#1}
  \seq_remove_all:Nn \l_arg_seq { }
  \seq_log:N \l_arg_seq
  \tl_set:Nn \l_cons_tl {#3}
  \tl_trim_spaces:N \l_cons_tl
  \seq_if_in:NxTF \l_arg_seq {\lref{\g_label_tl}} {
    \seq_remove_all:Nx \l_arg_seq {\lref{\g_label_tl}}
    \seq_get_left:NNTF \l_arg_seq \l_dummy_tl {
      \tl_trim_spaces:N \l_dummy_tl
      \bool_gset_false:N \g_firstargletter_bool
      \tl_if_head_eq_catcode:VNTF \l_dummy_tl a {
        \bool_gset_true:N \g_firstargletter_bool
      } {
        \tl_if_head_eq_meaning:VNTF \l_dummy_tl {\cref} {
          \tl_set:Nx \l_tmpa_tl {\tl_tail:N \l_dummy_tl}
          \tl_set:Nx \l_tmpb_tl {\tl_head:N \l_tmpa_tl}
          \bool_gset_true:N \g_firstargletter_bool
          \tl_if_in:NeTF \l_tmpb_tl {lem\c_colon_str} {} {
            \tl_if_in:NeTF \l_tmpb_tl {thm\c_colon_str} {} {
              \tl_if_in:NeTF \l_tmpb_tl {prop\c_colon_str} {} {
                \tl_if_in:NeTF \l_tmpb_tl {cor\c_colon_str} {} {
                  \bool_gset_false:N \g_firstargletter_bool
                }
              }
            }
          }
        } {
        }
      }
      \bool_if:NTF \g_firstargletter_bool {
        \seq_set_eq:NN \l_tmpa_seq \g_arg_mru_this
        \seq_remove_all:Nn \l_tmpa_seq {H}
        \seq_get_right:NN \l_tmpa_seq \l_tmpa_tl
        \int_case:nnF {\seq_count:N \l_arg_seq} {
          {1} {
            \str_case:VnF {\l_tmpa_tl} {
              {Ahpr} {
                \bool_if:NT \g_debug_bool {C1.1}
                \seq_gput_left:Nn \g_arg_mru_this {Ahpr}
                \seq_gremove_duplicates:N \g_arg_mru_this
                \enum:nxnn {#1} {\lref{\g_label_tl}} {-} {\BooleanTrue}
                \hence~
                \bool_if:NTF \l_plural_bool {
                  \prove[\l_arg_verbs_tl]~\ignorespaces #3
                } {
                  \proves[\l_arg_verbs_tl]~\ignorespaces #3
                }
              }
              {Tapr} {
                \bool_if:NT \g_debug_bool {C1.2}
                \seq_gput_left:Nn \g_arg_mru_this {Tapr}
                \seq_gremove_duplicates:N \g_arg_mru_this
                \enum[\lref{\g_label_tl}]{
                  This;
                  #1
                }[\l_arg_verbs_tl]\ignorespaces #3
              }
              {Ctapr} {
                \bool_if:NT \g_debug_bool {C1.3}
                \seq_gput_left:Nn \g_arg_mru_this {Ctapr}
                \seq_gremove_duplicates:N \g_arg_mru_this
                Combining~
                \enum[\lref{\g_label_tl}]{
                  this;
                  #1
                } \proves[\l_arg_verbs_tl]~\ignorespaces #3
              }
            } {}
          }
        } {
          \str_case:VnF {\l_tmpa_tl} {
             {Ahpr} {
              \bool_if:NT \g_debug_bool {C2.1}
              \seq_gput_left:Nn \g_arg_mru_this {Ahpr}
              \seq_gremove_duplicates:N \g_arg_mru_this
              \enum:nxnn {#1} {\lref{\g_label_tl}} {-} {\BooleanTrue}
              \hence~
              \prove[\l_arg_verbs_tl]~\ignorespaces #3
            }
            {Tapr} {
              \bool_if:NT \g_debug_bool {C2.2}
              \seq_gput_left:Nn \g_arg_mru_this {Tapr}
              \seq_gremove_duplicates:N \g_arg_mru_this
              \enum[\lref{\g_label_tl}]{
                This;
                #1
              }[\l_arg_verbs_tl]\ignorespaces #3
            }
            {Ctapr} {
              \int_case:nn {\int_rand:nn {0} {1}} {
                {0} {
                  \bool_if:NT \g_debug_bool {C2.3}
                  \seq_gput_left:Nn \g_arg_mru_this {Ctapr}
                  \seq_gremove_duplicates:N \g_arg_mru_this
                  Combining~
                  \enum[\lref{\g_label_tl}]{
                    this;
                    #1
                  } \proves[\l_arg_verbs_tl]~\ignorespaces #3
                }
                {1} {
                  \bool_if:NT \g_debug_bool {C2.4}
                  \seq_gput_left:Nn \g_arg_mru_this {Ctapr}
                  \seq_gremove_duplicates:N \g_arg_mru_this
                  Combining~
                  \enum:nxnn {#1} {\lref{\g_label_tl}} {-} {\BooleanFalse}
                  \hence~
                  \proves[\l_arg_verbs_tl]~\ignorespaces #3
                }
              }
            }
          } {}
        }
      } {
        \seq_set_eq:NN \l_tmpa_seq \g_arg_mru_this
        \seq_remove_all:Nn \l_tmpa_seq {H}
        \seq_remove_all:Nn \l_tmpa_seq {Ahpr}
        \seq_get_right:NN \l_tmpa_seq \l_tmpa_tl
        \int_case:nnF {\seq_count:N \l_arg_seq} {
          {1} {
            \str_case:VnF {\l_tmpa_tl} {
              {Tapr} {
                \bool_if:NT \g_debug_bool {C3.1}
                \seq_gput_left:Nn \g_arg_mru_this {Tapr}
                \seq_gremove_duplicates:N \g_arg_mru_this
                \enum[\lref{\g_label_tl}]{
                  This;
                  #1
                }[\l_arg_verbs_tl]\ignorespaces #3
              }
              {Ctapr} {
                \bool_if:NT \g_debug_bool {C3.2}
                \seq_gput_left:Nn \g_arg_mru_this {Ctapr}
                \seq_gremove_duplicates:N \g_arg_mru_this
                Combining~
                \enum[\lref{\g_label_tl}]{
                  this;
                  #1
                } \proves[\l_arg_verbs_tl]~\ignorespaces #3
              }
            } {}
          }
        } {
          \str_case:VnF {\l_tmpa_tl} {
            {Tapr} {
              \bool_if:NT \g_debug_bool {C4.1}
              \seq_gput_left:Nn \g_arg_mru_this {Tapr}
              \seq_gremove_duplicates:N \g_arg_mru_this
              \enum[\lref{\g_label_tl}]{
                This;
                #1
              }[\l_arg_verbs_tl]\ignorespaces #3
            }
            {Ctapr} {
              \int_case:nn {\int_rand:nn {0} {1}} {
                {0} {
                  \bool_if:NT \g_debug_bool {C4.2}
                  \seq_gput_left:Nn \g_arg_mru_this {Ctapr}
                  \seq_gremove_duplicates:N \g_arg_mru_this
                  Combining~
                  \enum[\lref{\g_label_tl}]{
                    this;
                    #1
                  } \proves[\l_arg_verbs_tl]~\ignorespaces #3
                }
                {1} {
                  \bool_if:NT \g_debug_bool {C4.3}
                  \seq_gput_left:Nn \g_arg_mru_this {Ctapr}
                  \seq_gremove_duplicates:N \g_arg_mru_this
                  Combining~
                  \enum:nxnn {#1} {\lref{\g_label_tl}} {-} {\BooleanFalse}
                  \hence~
                  \proves[\l_arg_verbs_tl]~\ignorespaces #3
                }
              }
            }
          } {}
        }
      }
    } {
      \tl_if_head_eq_catcode:oNTF \l_cons_tl a {
        \seq_set_eq:NN \l_tmpa_seq \g_arg_mru_this
        \seq_remove_all:Nn \l_tmpa_seq {Ctapr}
        \seq_remove_all:Nn \l_tmpa_seq {Ahpr}
        \seq_get_right:NN \l_tmpa_seq \l_tmpa_tl
        \str_case:VnF {\l_tmpa_tl} {
          {H} {
            \bool_if:NT \g_debug_bool {C5.1}
            \seq_gput_left:Nn \g_arg_mru_this {H}
            \seq_gremove_duplicates:N \g_arg_mru_this
            Hence,~we~obtain~\ignorespaces #3
          }
          {Tapr} {
            \bool_if:NT \g_debug_bool {C5.2}
            \seq_gput_left:Nn \g_arg_mru_this {Tapr}
            \seq_gremove_duplicates:N \g_arg_mru_this
            This~\proves[\l_arg_verbs_tl]~\ignorespaces #3
          }
        } {}
      } {
        \bool_if:NT \g_debug_bool {C6.1}
        \seq_gput_left:Nn \g_arg_mru_this {Tapr}
        \seq_gremove_duplicates:N \g_arg_mru_this
        This~\proves[\l_arg_verbs_tl]~\ignorespaces #3
      }
    }
  } {
    \int_compare:nNnTF {\seq_count:N \l_arg_seq} = {0} {
      \bool_if:NTF \g_arg_start_bool {
        \bool_if:NT \g_debug_bool {C7.1}
        \Nobs\unskip
        #3
      } {
        \bool_if:NT \g_debug_bool {C7.2}
        \Moreover~
        #3
      }
    } {
      \bool_if:NTF \g_arg_start_bool {
        \bool_if:NT \g_debug_bool {C8.1}
        \tl_log:N \l_arg_verbs_tl
        \Nobs~that~
        \enum{
          #1
        }[\l_arg_verbs_tl]\ignorespaces #3
      } {
        \int_compare:nNnTF {\seq_count:N \l_arg_seq} = {1} {
          \seq_set_eq:NN \l_tmpa_seq \g_arg_mru_nothis
          \seq_remove_all:Nn \l_tmpa_seq {Nwc}
          \seq_remove_all:Nn \l_tmpa_seq {Itnswc}
          \seq_get_right:NN \l_tmpa_seq \l_tmpa_tl
        } {
          \seq_get_right:NN \g_arg_mru_nothis \l_tmpa_tl
        }
        \str_case:VnF {\l_tmpa_tl} {
          {Mo} {
            \bool_if:NT \g_debug_bool {C9.1}
            \seq_gput_left:Nn \g_arg_mru_nothis {Mo}
            \seq_gremove_duplicates:N \g_arg_mru_nothis
            Moreover,~\nobs~that~
            \enum{
              #1
            }[\l_arg_verbs_tl]\ignorespaces #3
          }
          {Fm} {
            \bool_if:NT \g_debug_bool {C9.2}
            \seq_gput_left:Nn \g_arg_mru_nothis {Fm}
            \seq_gremove_duplicates:N \g_arg_mru_nothis
            Furthermore,~\nobs~that~
            \enum{
              #1
            }[\l_arg_verbs_tl]\ignorespaces #3
          }
          {Ia} {
            \bool_if:NT \g_debug_bool {C9.3}
            \seq_gput_left:Nn \g_arg_mru_nothis {Ia}
            \seq_gremove_duplicates:N \g_arg_mru_nothis
            In~addition,~\nobs~that~
            \enum{
              #1
            }[\l_arg_verbs_tl]\ignorespaces #3
          }
          {N} {
            \bool_if:NT \g_debug_bool {C9.4}
            \seq_gput_left:Nn \g_arg_mru_nothis {N}
            \seq_gremove_duplicates:N \g_arg_mru_nothis
            Next,~\nobs~that~
            \enum{
              #1
            }[\l_arg_verbs_tl]\ignorespaces #3
          }
          {Itns} {
            \bool_if:NT \g_debug_bool {C9.5}
            \seq_gput_left:Nn \g_arg_mru_nothis {Itnswc}
            \seq_gput_left:Nn \g_arg_mru_nothis {Itns}
            \seq_gremove_duplicates:N \g_arg_mru_nothis
            In~the~next~step~we~\nobs~that~
            \enum{
              #1
            }[\l_arg_verbs_tl]\ignorespaces #3
          }
          {Nwc} {
            \bool_if:NT \g_debug_bool {C9.6}
            \seq_gput_left:Nn \g_arg_mru_nothis {Nwc}
            \seq_gremove_duplicates:N \g_arg_mru_nothis
            Next~we~combine~
            \enum{
              #1
            }to~obtain~\ignorespaces #3
          }
          {Itnswc} {
            \bool_if:NT \g_debug_bool {C9.7}
            \seq_gput_left:Nn \g_arg_mru_nothis {Itns}
            \seq_gput_left:Nn \g_arg_mru_nothis {Itnswc}
            \seq_gremove_duplicates:N \g_arg_mru_nothis
            In~the~next~step~we~combine~
            \enum{
              #1
            }to~obtain~\ignorespaces #3
          }
        } {}
      }
    }
  }
  \bool_gset_false:N \g_arg_start_bool
  \bool_gset_false:N \l_insidearg_bool
  \cfload[.]
  \color{black}
}

\tl_new:N \g_label_tl
\tl_gset:Nn \g_label_tl { }

\NewDocumentCommand{\savelabel}{m}{
  \bool_if:NTF \l_insidearg_bool {
    \tl_gset:Nn \g_label_tl {#1}
  } {
    \tl_gset:Nn \g_label_tl { }
  }
}

\ExplSyntaxOff

\ExplSyntaxOn

\NewDocumentEnvironment {athm} {m m o} {
\str_if_eq:noTF {example} {#1} {
  \bool_gset_true:N \g_example_bool
} {
  \bool_gset_false:N \g_example_bool
}
\cfclear
\IfNoValueTF{#3}{
\begin{#1}\label{#2}\global\def\loc{#2}
}{
\begin{#1}[#3]\label{#2}\global\def\loc{#2}
}
}{
\end{#1}
}

\NewDocumentEnvironment {adef} {m} {
\begin{definition}\label{#1}\global\def\loc{#1}
}{
\end{definition}
}

\NewDocumentEnvironment{aproof} {} {
\bool_if:NTF \g_example_bool {
  \bool_gset_true:N \g_arg_start_bool
  \begin{proof}[Proof~for~\cref{\loc}]
} {
  \bool_gset_true:N \g_arg_start_bool
  \begin{proof}[Proof~of~\cref{\loc}]
}
\bool_gset_false:N \g_finishproof_bool
}{
\bool_if:NTF \g_finishproof_bool {}
{\finishproofthus}
\end{proof}
}

\NewDocumentCommand{\finishproofthus} {} {
  \bool_gset_true:N \g_finishproof_bool
  \bool_if:NTF \g_example_bool {
    The~proof~for~\cref{\loc}~is~thus~complete.
  } {
    The~proof~of~\cref{\loc}~is~thus~complete.
  }
}
\NewDocumentCommand{\finishproofthis} {} {
  \bool_gset_true:N \g_finishproof_bool
  \bool_if:NTF \g_example_bool {
    This~completes~the~proof~for~\cref{\loc}.
  } {
    This~completes~the~proof~of~\cref{\loc}.
  }
}

\ExplSyntaxOff

\NewDocumentEnvironment{cproof}{m}
{\begin{proof}[Proof of \cref{#1}]}%
{\noindent The proof of \cref{#1} is thus complete.
\end{proof}}

\NewDocumentEnvironment{cproof2}{m}
{\begin{proof}[Proof of \cref{#1}]}%
{\noindent This completes the proof of \cref{#1}.
\end{proof}}

\hypersetup{
    colorlinks,
    linkcolor={blue!80!black},
    citecolor={red},
    urlcolor={blue!80!black}
}

\makeatletter
\ExplSyntaxOn
\seq_new:N \g_abbrs
\prop_new:N \g_abbr_counts
\tl_new:N \l_abbr_count_tl

\ExplSyntaxOn

\bool_new:N \g_forexample

\NewDocumentCommand{\eg}{ o }{
	\IfValueT{#1}{
		\str_if_eq:noTF {fe} {#1} {
			\bool_gset_true:N \g_forexample
		} {\bool_gset_false:N \g_forexample}
	}
	\bool_if:nTF { \g_forexample } {
		\bool_gset_false:N \g_forexample
		for~example
	}{
		\bool_gset_true:N \g_forexample
		for~instance
	}
}

\NewDocumentCommand{\abbr}{m m O{#1} m m O{#4} m}{
	\expandafter\newcommand\csname#3\endcsname[1][]{
		\seq_if_in:NnTF \g_abbrs {#1} {
			\prop_get:NnN \g_abbr_counts {#1} \l_abbr_count_tl
			\prop_gput:Nnx \g_abbr_counts {#1} {\int_eval:n {\l_abbr_count_tl + 1}}
			\hyperref[#1]{#7}
		} {
			\seq_gput_left:Nn \g_abbrs {#1}
			\prop_gput:Nnn \g_abbr_counts {#1} {1}
			\expandafter\gdef\csname#1@def\endcsname{#2}
			\phantomsection\label{#1}
			\str_if_eq:nnTF{##1}{}{\emph{#2}}{##1}~(\hyperref[#1]{#7})
		}
	}
	\expandafter\newcommand\csname#6\endcsname[1][]{
		\seq_if_in:NnTF \g_abbrs {#1} {
			\prop_get:NnN \g_abbr_counts {#1} \l_abbr_count_tl
			\prop_gput:Nnx \g_abbr_counts {#1} {\int_eval:n {\l_abbr_count_tl + 1}}
			\hyperref[#1]{#4}
		} {
			\expandafter\gdef\csname#1@def\endcsname{#5}
			\seq_gput_left:Nn \g_abbrs {#1}
			\prop_gput:Nnn \g_abbr_counts {#1} {1}
			\phantomsection\label{#1}
			\str_if_eq:nnTF{##1}{}{\emph{#5}}{##1}~(\hyperref[#1]{#4})
		}
	}
}

\ExplSyntaxOff
		\makeatother
		\abbr{Adamax}{adaptive moment estimation maximum SGD}{Adamaxs}{Adaptive moment estimation maximum SGD}{Adamax}
		\abbr{Adagrad}{adaptive gradient SGD}{Adagrads}{adaptive gradient SGD}{Adagrad}
        \abbr{HJB}{Hamilton-Jacobi-Bellman}{HJBs}{adaptive gradient SGD}{HJB}
		\abbr{Nadam}{Nesterov accelerated adaptive moment estimation SGD}{Nadams}{Nesterov accelerated adaptive moment estimation SGD}{Nadam}
		\abbr{RMSprop}{root mean square propagation SGD}{RMSprops}{the root mean square propagation}{RMSprop}
		\abbr{ANN}{artificial neural network}{ANNs}{artificial neural networks}{ANN}
\abbr{DNN}{deep artificial neural network}{DNNs}{deep artificial neural networks}{DNN}
		\abbr{SGD}{stochastic gradient descent}{SGDs}{stochastic gradient descent}{SGD}
		\abbr{PR}{Polyak--Ruppert}{PRs}{Polyak--Rupperts}{PR}
        \abbr{PINN}{physics-informed neural network}{PINNs}{physics-informed neural networks}{PINNs}
		\abbr{PADAM}{parallel averaged \Adam}{PADAMs}{parallel averaged \Adams}{PADAM}
		\abbr{DGM}{deep Galerkin method}{DGMs}{deep Galerkin methods}{DGM}
        \abbr{DSM}{deep splitting method}{DSMs}{deep splitting methods}{DSM}
		\abbr{COD}{curse of dimensionality}{CODs}{curse of dimensionalities}{COD}
		\abbr{NSFC}{National Science Foundation of China}{NSFCs}{National Science Foundation of China}{NSFC}
		\abbr{ICBS}{International Congress of Basic Science}{ICBSs}{International Congress of Basic Sciences}{ICBS} 
		\abbr{PI}{Principal Investigator}{PIs}{Principal Investigators}{PI}
		\abbr{ERC}{European Research Council}{ERCs}{European Research Councils}{ERC}
		\abbr{MC}{Monte Carlo}{MCs}{Monte Carlos}{MC}
		\abbr{MLP}{multilevel Picard}{MLPs}{multilevel Picard}{MLP}
		\abbr{SDE}{stochastic differential equation}{SDEs}{stochastic differential equations}{SDE}
		\abbr{BSDE}{backward stochastic differential equation}{BSDEs}{backward stochastic differential equations}{BSDE}
		\abbr{SNSF}{Swiss National Science Foundation}{SNSFs}{Swiss National Science Foundations}{SNSF}
		\abbr{EMS}{European Mathematical Society}{EMSs}{European Mathematical Societies}{EMS}
		\abbr{EMA}{exponential moving average}{EMAs}{exponential moving averages}{EMA}
		\abbr{GD}{gradient descent}{GDs}{gradient descents}{GDs}
        \abbr{DL}{deep learning}{DLs}{deep learnings}{DL}
        \abbr{ELI}{elementary logical implication}{ELIs}{elementary logical implications}{ELIs}
        \abbr{DRL}{deep reinforcement learning}{DRLs}{deep reinforcement learnings}{DRL}
        \abbr{Adam}{adaptive moment estimation \SGD}{Adams}{adaptive moment estimation \SGD}{Adam}
        \abbr{AdamW}{AdamW}{AdamWs}{AdamWs}{AdamW}
		\abbr{iid}{independent and identically distributed}{i.i.d}{independent and identically distributed}{i.i.d.}
        \abbr{RePUs}{rectified power unit}{RePU}{rectified power unit}{RePU}
		\abbr{AI}{artificial intelligence}{AIs}{artificial intelligence}{AI}
		\abbr{ZF}{Zermelo--Fraenkel set theory}{ZFs}{Zermelo--Fraenkel set theories}{ZF}
        \abbr{AIBAP}{\AI\ based automatic proof verification}{AIBAPs}{\AI\ based automatic proof verification}{AIBAP}
		\abbr{LLM}{large language model}{LLMs}{large language models}{LLM}
		\abbr{PDE}{partial differential equation}{PDEs}{partial differential equations}{PDE}
		\abbr{DKM}{deep Kolmogorov method}{DKMs}{deep Kolmogorov methods}{DKM}
		\abbr{as}{almost surely}{ass}{almost surely}{a.s.}
		\abbr{pas}{$\P$-almost surely}{pass}{almost surely}{$\P$-a.s.}
		\abbr{ReLUs}{rectified linear unit}{ReLU}{rectified linear unit}
		\makeatother

\allowdisplaybreaks

        \newcommand{\reli}[3]{\mathcal{N}^{#1,#3}_{#2}\cfadd{definition: ANN}}
\newcommand{\relii}[4]{\mathcal{N}^{#1,#3}_{#2,#4}\cfadd{definition: ANN}}

\newcommand{\ffd}{\mathfrak{d}\cfadd{definition: ANN}}

			\newcommand\restr[2]{{
					\left.\kern-\nulldelimiterspace 
					#1 
					\vphantom{|} 
					\right|_{#2} 
			}}

\begin{document}
\title{Error analysis for the deep Kolmogorov method}
\author{Iulian Cîmpean$^{1,2}$, Thang Do$^{3}$, Lukas Gonon$^{4}$, Arnulf Jentzen$^{5,6}$, and Ionel Popescu$^{7,8}$
	\bigskip
	\\
     \small{$^1$ Faculty of Mathematics and Computer Science,}
	\vspace{-0.1cm}\\
	\small{University of Bucharest, Romania, e-mail: \texttt{iulian.cimpean@unibuc.ro}}
	\smallskip
	\\
     \small{$^2$ Institute of Mathematics “Simion Stoilow" of the Romanian Academy,}
	\vspace{-0.1cm}\\
	\small{Romania, e-mail: \texttt{iulian.cimpean@imar.ro}}
	\smallskip
\\
    	\small{$^3$ School of Data Science, The Chinese University of Hong Kong, Shenzhen}
	\vspace{-0.1cm}\\
	\small{ (CUHK-Shenzhen), China, e-mail: \texttt{minhthangdo@link.cuhk.edu.cn}}
 \smallskip
	\\
    \small{$^4$ Faculty of Mathematics and Statistics, University of St.\ Gallen,}
	\vspace{-0.1cm}\\
	\small{Switzerland, e-mail: \texttt{lukas.gonon@unisg.ch}}
 \smallskip
	\\
	\small{$^5$ School of Data Science and School of Artificial Intelligence, The Chinese University}
	\vspace{-0.1cm}\\
	\small{of Hong Kong, Shenzhen (CUHK-Shenzhen), China, e-mail: \texttt{ajentzen@cuhk.edu.cn}}
	\smallskip
	\\
 \small{$^6$ Applied Mathematics: Institute for Analysis and Numerics, Faculty of Mathematics and}
	\vspace{-0.1cm}\\
	\small{Computer Science, University of M{\"u}nster, Germany, e-mail: \texttt{ajentzen@uni-muenster.de}}
	\smallskip
	\\
    \small{$^7$ Faculty of Mathematics and Computer Science,}
	\vspace{-0.1cm}\\
	\small{University of Bucharest, Romania, e-mail: \texttt{ionel.popescu@fmi.unibuc.ro}}
	\smallskip
	\\
     \small{$^8$ Institute of Mathematics “Simion Stoilow" of the Romanian Academy,}
	\vspace{-0.1cm}\\
	\small{Romania, e-mail: \texttt{ionel.popescu@imar.ro}}
	\smallskip
\\
}
\date{}
\maketitle
\begin{abstract}
    The deep Kolmogorov method is a simple and popular deep learning based method for approximating solutions of partial differential equations (PDEs) of the Kolmogorov type. In this work we provide an error analysis for the deep Kolmogorov method for heat PDEs. Specifically, we reveal convergence with convergence rates for the overall mean squared distance between the exact solution of the heat PDE and the realization function of the approximating deep neural network (DNN) associated with a stochastic optimization algorithm in terms of the size of the architecture (the depth/number of hidden layers and the width of the hidden layers) of the approximating DNN, in terms of the number of random sample points used in the loss function (the number of input-output data pairs used in the loss function), and in terms of the size of the optimization error made by the employed stochastic optimization method.
\end{abstract}
\thispagestyle{empty}
\tableofcontents

\section{Introduction}
It is nowadays a very active area of research to design and study deep learning methods to approximately solve \PDEs\ and related scientific computing problems (cf., for example, the review articles and monographs \cite{MR4547087,MR4356985,MR4457972,Germainetal2021arXivOverviewArticle,MR4795589,ArBePhi2024}). Popular such deep learning approximation methods for \PDEs\ include \PINNs\ (cf.\ \cite{MR3881695}) and the \DGM\ (cf.\ \cite{MR3874585}) as well as the deep \BSDE\ method (cf.\ \cite{MR3847747,MR3736669}). Another quite simple and nowadays well-known deep learning approximation scheme is the \DKM\ to approximately solve \PDEs\ of the Kolmogorov type; see \cite{MR4293960} and \cite{MR3946472}. (In the references \cite{MR4293960} and \cite{MR3946472} the method has not been named as \DKM\ but this has been done in the follow-up work \cite[Chapter 17]{ArBePhi2024}.) We also refer, \eg, to \cite{MR4310910,MR4673229} for suitable \DSMs\ that are based on the combination of the splitting method and the \DKM.

Beyond the design and the numerical investigation of such deep learning approximation methods for \PDEs, it is also an active topic of research to rigorously analyze such methods mathematically, that is, to establish rigorous upper bounds for the approximation error between the exact solution of the \PDE\ and the deep learning approximation under consideration. While for classical approximation methods for solutions of \PDEs\ such as finite difference (cf., \eg, \cite{MR3136501} and the references therein) and finite element methods (cf., \eg, \cite{MR2249024} and the references therein) such error analyses for the difference between the \PDE\ solution and the considered numerical approximation are -- to a large extent -- well understood, the situation is much more subtle in the case of deep learning approximation methods for \PDEs\ and it basically remains a fundamental open problem of research to establish a full error analysis for any reasonable deep learning approximation scheme for \PDEs.

There are, however, several partial error analyses for deep learning methods for \PDEs\ in the literature. For instance, in the case of the \PINNs\ the work Mishra \& Molinaro \cite[Theorem 3.1]{MR4565573} explicitly bounds the overall root mean squared distance between the exact solution of the \PDE\ and the realization of the approximating \ANN\ from above by the product of an error explicit constant and the sum of the loss function and the number of discretization points where the error constant also depends on the realization function of the approximating \ANN. A qualitative convergence result without a convergence speed for the \DKM\ applied to semilinear Kolmogorov \PDEs\ in terms of the expectation of the loss function can also be found in \cite[Proposition 4.7, Appendix A.3, and (21)]{MR3946472}.

In this work we push these error analyses, in the case of the \DKM\ applied to heat \PDEs, a bit further. Specifically, in \cref{main theorem: conclude} in \cref{main theorem} below and \cref{main theorem 3} in \cref{sec: deep Kolmogorov method}, respectively, we establish that there exists a strictly positive real error constant $c \in (0, \infty)$ such that the overall mean squared distance between the exact solution of the \PDE\ and the random deep \ANN\ approximation converges to zero with explicit convergence constants and convergence rates using the
error constant $c$ as the optimization error of the underlying stochastic optimization algorithm converges to zero (convergence rate $1$ in terms of the optimization error), as the number of random sample points used in the loss function (the number of input-output data pairs used in the loss function) converges to infinity (convergence rate \nicefrac{1}{2} in terms of the number of random
sample points), and as the width of the approximating \ANN\ (the number of neurons on the hidden layers of the approximating \ANN) converges to infinity (convergence rate $2/(d+5)$ in terms of the width of the \ANN). The strictly positive positive real error constant $c \in (0,\infty)$ in \cref{main theorem: conclude} in \cref{main theorem} is completely independent of the number of random sample points in the loss function (the number of input-output data pairs in the loss function), is completely independent of the architecture of the approximating deep \ANNs\ (the depth/the number of hidden layers and the widths of the hidden layers of the approximation \ANNs), and is completely independent of the random \ANN\ parameter vector computed with stochastic optimization methods.

The specific analysis in this work is restricted to simple heat \PDEs, to \ANNs\ with the \ReLU\ activation, and to the \DKM, which seems to be easier to analyze than many other deep learning approximation methods for \PDEs\ such as \PINNs\ and deep \BSDE\ approximations, but we expect that our arguments can be extended to more challenging deep learning methods and \PDE\ problems.
\subsection{Artificial neural networks (ANNs)}
To formulate our main error analysis result within this introductory section in \cref{main theorem} below, we need to describe realization functions of \ANNs\ in accurate mathematical terms. This is precisely the subject of the notion in \cref{definition: ANN} below. In the simplest form, a realization function of an \ANN\ can be considered as a function that is given by multiple compositions of affine functions and multidimensional versions of a certain fixed one-dimensional nonlinear function, which is referred to as activation function of the \ANN. In \cref{definition: ANN} the natural number $L \in \N = \{ 1, 2, 3, \dots \}$ specifies the number of affine functions in the compositions and the natural numbers $\ell_0, \ell_1, \dots, \ell_L$ specify the dimensionalities of the domains and codomains, respectively, of the affine functions in the compositions.
\begin{samepage}
\begin{definition}[\ANNs]
\label{definition: ANN}
For every $ L \in \N $,
$ \ell = ( \ell_0, \ell_1, \dots, \ell_L ) \in \N^{ L + 1 } $
we denote by $ \ffd( \ell ) \in \N $ the natural number
given by
$
  \ffd( \ell ) = \sum_{i=1}^{L}\ell_i ( \ell_{ i - 1 } + 1 )
$
and for every $ L\in \N $,
$ \ell = ( \ell_0, \ell_1, \dots, \ell_L ) \in \N^{ L + 1 } $,
$ \theta=(\theta_1,\dots,\theta_{\ffd(\ell)}) \in \R^{ \ffd( \ell ) } $ we denote by
  $\reli{ \ell }{ v }{ \theta }
  =
  (
    \relii{ \ell }{ v }{ \theta }{ 1 },
    \dots,
    \relii{ \ell }{ v }{ \theta }{ \ell_v }
  )
  \colon \R^{ \ell_0 } \to \R^{ \ell_v }
$,
$ v \in \{ 1, 2, \dots, L \} $,
the functions which satisfy for all
$ v \in \{ 0, 1, \dots, L - 1 \} $,
$ x = (x_1, \dots, x_{ \ell_0 } ) \in \R^{ \ell_0 } $,
$ i \in \{ 1, 2, \dots, \ell_{ v + 1 } \} $
that
\begin{equation}
\label{realization multi}
\begin{split}
 & \relii{ \ell }{ v + 1 }{ \theta }{ i }( x )
  =
  \theta_{ \ell_{ v + 1 } \ell_v + i + \sum_{ h = 1 }^v \ell_h ( \ell_{ h - 1 } + 1 ) }
\\ &
  +
  \textstyle\sum_{ j = 1 }^{ \ell_v }
  \theta_{ (i - 1) \ell_v + j + \sum_{ h = 1 }^v \ell_h ( \ell_{ h - 1 } + 1 ) }
  \bigl(
    x_j
    \indicator{ \{ 0 \} }( v )
    +
    \max\{ \relii{ \ell }{ \max\{ v, 1 \} }{ \theta }{ j }(x),0\} 
    \indicator{ \N }( v )
  \bigr) .
\end{split}
\end{equation}
\end{definition}
\end{samepage}
In \cref{realization multi} appears the \ReLU\ activation function $ \R \ni x \mapsto \max\{ x, 0 \} \in \R $ 
and in \cref{definition: ANN} we thus describe \ANNs\ with the \ReLU\ activation function. Moreover, we note that the components of the vector $\theta = ( \theta_1, \dots, \theta_{ \fd( \ell ) } ) \in \R^{\fd(\ell)}$ in \cref{definition: ANN} represent the weight and the bias parameters of the considered \ANN. \cref{definition: ANN} essentially coincides with the formulation of realization functions of \ANNs\ in, \eg, \cite[Setting 2.1]{HannibalJentzenThang2024}.
\subsection{Error analysis for the deep Kolmogorov method}
Having specified realization functions of \ReLU\ \ANNs\ in \cref{definition: ANN} above, 
we are now in the position to formulate in \cref{main theorem} below one of the main error analysis results of this work. The function $u \colon [0,T] \times \R^d \to \R$ in \cref{main theorem} is the solution of the heat equation in \cref{main theorem: eq2} that we intend to approximate by deep \ANNs\ by means of the \DKM; see \cref{main theorem: eq3} for the loss function of the \DKM\ in which \iid\ random variables appear.

We now present \cref{main theorem} with all mathematical details and, thereafter, we briefly explain the contribution of \cref{main theorem} and some of the mathematical objects appearing in \cref{main theorem} in words.
\begin{samepage}
\begin{athm}{theorem}{main theorem}
Let $ T \in (0,\infty) $, $ d \in \N $, let $u\in C^2([0,T]\times\R^{d},\R)$ be\footnote{Note that for every $d, n \in \N$ and every closed $D \subseteq \R^d$ it holds that $C^n( D, \R ) = \{ f \in C( D, \R ) \colon ( \exists\, g \in C^n( \R^d, \R ) \colon g|_D = f ) \}$.} at most polynomially growing\footnote{Note that for every $d \in \N$, every $A \subseteq \R^d$, and every $f \colon A \to \R$ it holds that $f$ is at most polynomially growing if and only if there exists $c \in \R$ such that for all $x = ( x_1, \dots, x_d ) \in A$ it holds that $| f(x) | \leq c ( 1 + \sum_{ i = 1 }^d | x_i | )^c$.}, assume for all $t\in (0,T)$, $x\in \R^d$ that
\begin{equation}\label{main theorem: eq2}
\textstyle
  \frac{ \partial }{ \partial t } u( t, x )
  +
  \frac{ 1 }{ 2 }
  \Delta_x u(t,x)
  = 0,
\end{equation}
 let $ ( \Omega, \mathcal{F}, \P ) $ be a probability space,
let $ \mathbb{T}_n \colon \Omega \to [0,T] $, $n\in \N$, be independent uniformly distributed random variables, let $D\subseteq\R^d$ be open and bounded,
let $ \mathbb{X}_n \colon \Omega \to D $, $n\in \N$, be independent uniformly distributed random variables,
let $ W^n \colon [0,T] \times \Omega \to \R^d $, $n\in \N$, be \iid\ standard Brownian motions, 
assume that $ (\mathbb{T}_n)_{n\in \N} $, $ (\mathbb{X}_n)_{n\in \N} $, and $ (W^n)_{n\in \N} $ are independent,
and for every $M=(M_1,M_2)\in \N^2$, $ L \in \N \backslash \{ 1 \} $,
$
  \ell  \in
  ( \{ d + 1 \} \times \N^{ L - 1 } \times \{ 1 \} )
$
 let $ \mathbb{F}_{M,\ell} \colon \R^{ \ffd( \ell ) } \to \R $
satisfy
for all
$ \theta \in \R^{ \ffd( \ell ) } $
that
\begin{equation}\label{main theorem: eq3}
\textstyle
  \mathbb{F}_{M,\ell}( \theta )=
  \frac {1}{M_1} \sum_{m=1}^{M_1}
    \big|
      \reli{ \ell }{ L }{ \theta }( \mathbb{T}_m, \mathbb{X}_m )
      - \bigl[ \frac{ 1 }{ M_2 } \sum_{ n = 1 }^{ M_2 } u(T, \bbX_{ m } + W^{ m M_2 + n }_{ T-\bbT_m } ) \bigr]
    \big|^2.
\end{equation}
Then there exists $c \in (0,\infty)$ such that for every $R\in [c,\infty)$, $M=(M_1,M_2)\in \N^2$, $L\in \N\backslash\{1\}$, $\ell=(\ell_0,\ell_1,\dots,\ell_L)\in (\{d+1\}\times\N^{L-1}\times\{1\}) $ and every random variable $\vartheta \colon \Omega \to [-R,R]^{\fd(\ell)}$ it holds that
\begin{equation}\label{main theorem: conclude}
\begin{split}
&\textstyle
  \E[
    \int_{[0,T]\times D}
      | u(y) - \reli{ \ell }{ L}{ \vartheta }( y ) |^2
    \, \d y
  ]\textstyle\leq [R+\textstyle\max\{\ell_1,\ell_2,\dots,\ell_{L-1}\}]^{cL}(M_1)^{-\nicefrac{1}{2}}\\
  &+c [\min\{\ell_1, \ell_2, \dots, \ell_{ L - 1 } \}]^{-2/(d+5)}+c\,\E\bigl[\textstyle \bbF_{M,\ell}( \vartheta ) - \inf_{\theta\in [-R,R]^{\fd(\ell)}} \bbF_{M,\ell}( \theta ) \bigr]  \ifnocf.
  \end{split}
\end{equation}
\cfout[.]
\end{athm}
\end{samepage}
\cref{main theorem} is an immediate consequence of \cref{main cor} in \cref{subsec: qualitative results} below. \cref{main cor}, in turn, is established through an application of \cref{prop: linear II relu pinn} in \cref{sec: deep Kolmogorov method} below, which is the main result of this work (see \cref{subsec: qualitative results} for details). 

The strictly positive real number $T \in (0,\infty)$ in \cref{main theorem} describes the time horizon of the heat equation under consideration and the natural number $d \in \N$ specifies the space dimension of the heat equation under consideration. In \cref{main theorem: eq2} we present the heat \PDE\ whose solutions we intend to approximate by means of the \DKM\ in \cref{main theorem: eq3} and the function $u \colon [0,T]\times \R^d \to \R$ is a solution of the \PDE\ in \cref{main theorem: eq2}. The function specified in \cref{main theorem: eq3} is precisely the loss function of the \DKM\ (cf., \eg, \cite[Framework 3.2]{MR4293960} and \cite[Subsection 4.2]{MR3946472}). The vector $\ell = ( \ell_0, \ell_1, \dots, \ell_L ) \in ( \{ d + 1 \} \times \N^{ L - 1 } \times \{ 1 \} )$ in \cref{main theorem: eq3} and \cref{main theorem: conclude} specifies the architecture of the considered approximating \ANNs\ in the sense  
\begin{itemize}
\item that $\ell_0$ specifies the number of neurons on the $1\textsuperscript{st}$ layer of the \ANN\ (the number of neurons on the input layer), 
\item that $\ell_1$ specifies the number of neurons on the $2\textsuperscript{nd}$ layer of the \ANN\ (the number of neurons on the $1\textsuperscript{st}$ hidden layer), 
\item that $\ell_2$ specifies the number of neurons on the $3\textsuperscript{rd}$ layer of the \ANN\ (the number of neurons on the $2\textsuperscript{nd}$ hidden layer), 
\item $\ldots$ ,
\item that $\ell_{ L - 1 }$ specifies the number of neurons on the $L\textsuperscript{th}$ layer of the \ANN\ (the number of neurons on the $(L-1)\textsuperscript{th}$  hidden layer), and 
\item that $\ell_L$ specifies the number of neurons on the $(L+1)\textsuperscript{th}$ layer of the \ANN\ (the number of neurons on the output layer).
\end{itemize}
We note that \cref{main theorem: conclude} in \cref{main theorem} assures that there exists a strictly positive real error constant $c \in (0,\infty)$ such that for every upper bound $R \in [c,\infty)$ for the parameters of the \ANN\ associated to the stochastic optimization method, every number of random sample points $M_1 \in \N$ in the loss function (see \cref{main theorem: eq3}), every number of Monte Carlo samples $M_2 \in \N$ in the loss function (see \cref{main theorem: eq3}), every depth $L \in \N \backslash \{1\}$ of the approximating \ANNs, every architecture vector $\ell = ( \ell_0,\dots,\ell_L ) \in ( \{d+1\}\times\N^{L-1}\times\{1\} )$ for the approximating \ANNs, and every random variable $\vartheta \colon \Omega \to[-R,R]^{\fd(\ell)}$ resulting from the stochastic optimization method we have that the overall mean squared distance 
\begin{equation}
 \textstyle \E[
    \int_{[0,T]\times D}
      | u(y) - \reli{ \ell }{ L}{ \vartheta }( y ) |^2
    \, \d y
  ]
\end{equation}
on the bounded set $[0,T] \times D$ between the exact solution $u \colon [0,T]\times\R^d \to\R$ of the \PDE\ in \cref{main theorem: eq2} and the realization function $\reli{ \ell }{ L}{ \vartheta } \colon \R^{d+1} \to\R$ of the \ANN\ associated with the random \ANN\ parameter vector $\vartheta$ is bounded by the sum

\begin{enumerate}[label=(\roman*)]
    \item \label{item 1: introduction} of the term $[R+\textstyle\max\{\ell_1,\ell_2,\dots,\ell_{L-1}\}]^{cL}(M_1)^{-\nicefrac{1}{2}}$ (that converges to zero as $M_1$ increases to infinity),
    \item \label{item 2: introduction} of the term $c [\min\{\ell_1, \ell_2, \dots, \ell_{ L - 1 } \}]^{-2/(d+5)}$ (that converges to zero as the minimum $\min\{ \ell_1,\allowbreak \ell_2, \dots, \ell_{ L-1 }\}$ increases to infinity), and
    \item \label{item 3: introduction} of -- up to the constant $c$ -- the optimization error $\E\bigl[\textstyle \bbF_{M,\ell}( \vartheta ) - \inf_{\theta\in [-R,R]^{\fd(\ell)}} \bbF_{M,\ell}( \theta ) \bigr] $.
\end{enumerate}
We also note in \cref{main theorem} that the value of the number of Monte Carlo samples $M_2 \in \N$ does not influence the first and the second error terms on the right hand side of \cref{main theorem: conclude} (see \cref{item 1: introduction,item 2: introduction}  above) but does only influence the optimization error -- the third error term on the right hand side of \cref{main theorem: conclude} (see \cref{item 3: introduction} above). In particular, we note that even if $M_2$ tends to infinity, then one can not expect that the first error term on the right hand side of \cref{main theorem: conclude} (see \cref{item 1: introduction} above) vanishes. Nonetheless, in numerical simulations it may be beneficial to choose $M_2$ much larger than 1. In this context we also note that, under suitable assumptions, it may very well be possible to improve the convergence rate in $M_1$ in \cref{main theorem: conclude} (see \cref{item 1: introduction} above) from $\nicefrac{ 1 }{ 2 }$ to $1$ by employing the techniques in \cite[Lemma 4]{MR4134774} and \cite[Appendix C]{Hiebersupplement}.
\subsection{Literature review}
In this subsection we provide a short review on selected further results in the literature that develop error analyses for deep learning based approximation schemes for \PDEs\ and related problems. For more comprehensive reviews on deep learning approximation methods for \PDEs\ in the literature we refer, for instance, to the overview and review articles \cite{Germainetal2021arXivOverviewArticle,Gononetal2024arXiv,MR4547087,MR4356985,MR4270459,MR4412280,MR4457972,MR4795589} and the references therein.

The mathematical analysis of deep learning methods for solving \PDEs\ has attracted increasing attention in recent years. Numerous works have contributed insights into the approximation capabilities and theoretical guarantees of such methods. Substantial progress has been made on bounding individual components of the total error -- such as approximation, sampling, or optimization errors in the past years. Several works in the literature focus on analysing the approximation errors for neural network approximations to \PDE\ solutions for heat equations \cite{MR4454915},
Black Scholes \PDEs\
\cite{MR4574851,DePhiArSch2022},
non-local \PDEs\ \cite{GononSch2021,GononSch2023,Valenzuela2022,ArTuanWu2022},
parametric \PDEs\ \cite{Gitta2022},
smoothness spaces associated to \PDEs \cite{Ingo2021,MR4131037,MR4376564,MR4420580,Beneventano2020},
semilinear (heat) equations \cite{MR4292849,NeufeldNguyen2024arXiv,Ackermanndeepneural2023,Ackermann2024,MaArBe2021},
non-linear waves \cite{Claudiobound2024},
zero-sum games
\cite{MR4154658},
stochastic volatility models
\cite{MR4787434},
optimal stopping problems 
\cite{MR4765842},
 or McKean-Vlasov equations
\cite{MR4776392}.
Several articles analyze approximation and sampling errors jointly by focusing on the empirical risk minimization problem (see, \eg, \cite{NEURIPS2020_c1714160,Jichangerroranalysis2023,MR4127967}).
Concerning the analysis of the generalization error we refer to 
\cite{TaoHaizhao2020}
for the analysis of the generalization error for deep learning-based approximations for certain type of second-order linear \PDEs\ and two-layer networks. For an analysis of the generalization error of several related deep learning methods for elliptic \PDEs\ we refer to 
\cite{KaiGuerroranalysis}
for deep mixed residual methods, to \cite{MR3767958,JiaoLaiLoWangYang2021,Lu21a}
for the deep Ritz method, to \cite{JiaoLaiWangYangYan2023} for a \DGM\ for weak solutions, and to \cite{ar2408.13511} for a method for computing eigenvalues of the Schrödinger operator using deep \ANN\ test functions.
Full error analysis for random feature neural network-based methods are provided in \cite{MR4633578,NeufeldSchmokerWu2024,NeufeldSchmoker2023}.
An error analysis for an \ANN-based algorithm for ergodic mean field control is provided in \cite{MR4264647}.

In \cite[Theorem 2.4]{TaoHaizhao2020} and \cite[Theorem 3.1]{MR4565573} Mishra and Molinaro provide upper estimates for the generalization/overall error of \PINNs\ for approximating \PDEs\ and inverse problems with error constants depending on the approximating (\ANN\ realization) function. The error bounds in \cite{MR4793683} provide separate analyses of the relevant error components. The work  
\cite{MR4600824} proposes and mathematically analyses Friedrichs learning as a novel deep learning methodology that can learn weak solutions of \PDEs\ via a minimax formulation. 

Moreover, there are various convergence results for minimizers of loss functions of deep learning algorithms in the literature focusing on discretization or simulation errors, without explicit approximation rates for the underlying \ANNs\ or bounds in terms of the optimization error. For example, the following articles provide such convergence results for the deep \BSDE\ method and related (forward) \BSDE\ based methods
\cite{MR4122227,MR4551869,MR4399896,Belak2021,ar2412.11010,ar2505.18297,ar2212.14372,ar2405.11392,ar2301.12895,Bichuch2024,MR4237495},
deep backward schemes \cite{MR4081911,MR4500281,MR4358471,ar2407.14566},
the deep splitting scheme \cite{MR4865059,NeufeldSchmokerWu2024},
the deep fictitious play algorithm \cite{MR4484148},
algorithms for first order \HJB\ equations \cite{MR4646437},
methods related to the \DGM\ \cite{MR4447670,MR3874585},
\PINNs\ \cite{MR4548578,ar2402.15592,ar2409.17938,HuShukla2024},
and a randomized quasi-Monte Carlo based deep learning method \cite{ar2310.18100}.

Results for related algorithms are provided, \eg, for regularized dynamical parametric approximations in \cite{ar2403.19234} for the variational Monte Carlo algorithm in \cite{MR4047010},
for physics-informed kernel learning in \cite{ar2409.13786,ar2402.07514},
for a Monte Carlo Euler scheme in \cite{ar2306.16811},
for a quasi-Regression Monte Carlo Scheme in \cite{MR4103948},
and for a policy iteration algorithm \cite{MR4779451}. 

We also refer, \eg, to \cite[Lemma 4]{MR4134774} and \cite[Appendix C]{Hiebersupplement} for abstract general error estimates for \ANN\ approximations in terms of the optimization error.

A key novelty of the error analyses in \cref{main theorem} (and \cref{prop: linear II relu pinn}, respectively) is that the overall mean squared distance between the exact solution of the \PDE\ and the deep \ANN\ approximation (the overall mean squared error of the deep learning method) is estimated in terms of the number of random sample points (input-output data pairs) used in the loss function (convergence rate $\nicefrac{ 1 }{ 2 }$ for the number of random sample points used in the loss function), in terms of the width of the \ANN\ (convergence rate $2/(d+5)$ for the width of the \ANN), and in terms of the \emph{optimization error} $\E\bigl[\textstyle \bbF_{M,\ell}( \vartheta ) - \inf_{\theta\in [-R,R]^{\fd(\ell)}} \bbF_{M,\ell}( \theta ) \bigr]$ (convergence rate $1$ for the optimization error) with the error constant $c \in (0,\infty)$ being completely independent of the number of random sample points, of the \ANN\ architecture, and of the optimization method.

It should, however, be pointed out that the error analysis in this work is restricted to the simple \DKM, to simple heat \PDEs, and to \ANNs\ with the \ReLU\ activation. Nonetheless, we expect that our arguments can be extended to more general classes of activation functions and \PDE\ problems and we plan to treat more advanced deep learning methods and \PDE\ problems in future research works.

\subsection{Structure of this article}
The remainder of this article is structured in the following way. After this introductory section, we establish in \cref{sec: priori estimate} several elementary explicit a priori estimates for realization functions of \ANNs. In \cref{sec: ANN approximation} we present in \cref{cor: Pinns for Poisson3 intro2 pre2 approximation} suitable approximation error estimates for deep \ANN\ approximations of twice continuously differentiable functions. Finally, in \cref{sec: deep Kolmogorov method} we combine the findings from \cref{sec: priori estimate} and \cref{sec: ANN approximation} to establish in \cref{prop: linear II relu pinn draft}, \cref{prop: linear II relu pinn}, \cref{cor: linear II relu pinn}, and \cref{main cor} error estimates for the \DKM. \cref{main theorem} above is a direct consequence of \cref{main cor}.
\section{Explicit a priori estimates for ANNs}\label{sec: priori estimate}
In this section we establish in \cref{lem: anngrowth}, \cref{lem: ann lipschitz}, and \cref{cor: ann lipschitz} several elementary a priori estimates for realizations of arbitrarily deep \ANNs\ with general activation functions. The arguments in our proofs of \cref{lem: anngrowth}, \cref{lem: ann lipschitz}, and \cref{cor: ann lipschitz} are based, \eg, on \cite[Subsection 2.1.5]{Beckfullerror2019} and \cite[Theorem 2.6 and Section A.3]{MR4127967}, respectively. 

In \cref{sec: deep Kolmogorov method} we employ the a priori estimates in \cref{lem: anngrowth} and \cref{cor: ann lipschitz} to establish in \cref{prop: linear II relu pinn draft} error estimates for the \DKM\ and, thereby, to prove \cref{main theorem} in the introduction.
\subsection{Explicit a priori estimates for realization functions of ANNs}\label{subsec: explicit estimates}
In this section we establish several elementary a priori estimates for realizations of \ANNs\ with a general activation function which does not necessarily coincide with the \ReLU\ activation function. In view of this, we have no use of the realizations for \ReLU\ \ANNs\ specified in \cref{realization multi} in \cref{definition: ANN} above but instead we need to specify realizations of \ANNs\ with a general activation function $\bbA \colon \R \to \R$. This is precisely the subject of the framework provided in \cref{setting: ANN} below.
\begin{setting}\label{setting: ANN} Let $d,L\in \N$,
$ \ell = ( \ell_0, \ell_1, \dots, \ell_L ) \in (\{d\}\times\N^{L-1}\times\{1\}) $, for every $k\in \{0,1,\dots,L\}$ let $\bfd_k\in \R$ satisfy $\bfd_k=\sum_{h=1}^k\ell_h(\ell_{h-1}+1)$, let $\bbA\colon \R\to\R$ be a function, and for every 
$ \theta=(\theta_1,\dots,\theta_{\ffd( \ell ) }) \in \R^{ \ffd( \ell ) } $ let
$
  \Reli{ \ell }{ v }{ \theta }{ \Act }
  =
  (
    \Relii{ \ell }{ v }{ \theta }{ \Act }{ 1 },
    \dots,
    \Relii{ \ell }{ v }{ \theta }{ \Act }{ \ell_v }
  )
  \colon \R^{ \ell_0 } \to \R^{ \ell_v }
$,
$ v \in \{ 1, 2, \dots, L \} $,
 satisfy for all
$ v \in \{ 0, 1, \dots, L - 1 \} $,
$ x = (x_1, \dots, x_{ \ell_0 } ) \in \R^{ \ell_0 } $,
$ i \in \{ 1, 2, \dots, \ell_{ v + 1 } \} $
that
\begin{equation}
\label{setting ANN: realization multi}
\begin{split}
  \Relii{ \ell }{ v + 1 }{ \theta }{ \Act }{ i }( x )
&
  =
  \theta_{ \ell_{ v + 1 } \ell_v + i +\bfd_v }
\\ &\quad
  +
  \textstyle\sum_{ j = 1 }^{ \ell_v }
  \theta_{ (i - 1) \ell_v + j +\bfd_{v} }
  \bigl(
    x_j
    \indicator{ \{ 0 \} }( v )
    +
    \Act( \Relii{ \ell }{ \max\{ v, 1 \} }{ \theta }{ \Act }{ j } (x))
    \indicator{ \N }( v )
  \bigr) \ifnocf.
\end{split}
\end{equation}
\cfload[.]
\end{setting}

Using the notation from \cref{setting: ANN} above, we now establish in \cref{lem: anngrowth} below an elementary a priori estimate for the evaluation of the realization of an \ANN\ with the general activation function $\bbA \colon \R\to \R$ from \cref{setting: ANN}.
\begin{athm}{lemma}{lem: anngrowth}
    Assume \cref{setting: ANN} and let $\polycons\in [0,\infty)$, $c,a\in \R$, $b\in (a,\infty)$ satisfy for all $x\in \R$ that $|\bbA(x)|\leq c(1+|x|^\polycons)$. Then it holds for all $k\in \{1,2,\dots,L\}$, $\theta=(\theta_1,\dots,\theta_{\fd(\ell)})\in \R^{\fd(\ell)}$ that
    \begin{equation}\llabel{need to prove}
  \begin{split}
 &\textstyle\max_{i\in \{1,2,\dots,\ell_{k}\}}\sup_{x\in [a,b]^d}|  \Relii{ \ell }{ k }{ \theta }{ \Act}{i }( x )|\\
 &\textstyle
 \leq [\max\{1,|a|,|b|\}]^{\max\{1,\polycons^{k-1}\}}\bigl[\max\{1,\max_{j\in \{1,2,\dots,\bfd_{k}\}}|\theta_j|\}\bigr]^{k\max\{1,\polycons^{k-1}\}}\\
&\textstyle\quad\cdot\bigl[1+2\max\{c,1\}\max\{\ell_0,\ell_1,\dots,\ell_{k-1}\}\bigr]^{k\max\{1,\polycons^{k-1}\}}\ifnocf.
\end{split}
    \end{equation}
    \cfout[.]
\end{athm}
\begin{aproof}
    We prove \lref{need to prove} by induction on $k\in \{1,2,\dots,L\}$. For the base case $k=1$ note that \cref{setting ANN: realization multi} shows that for all $\theta=(\theta_1,\dots,\theta_{\fd(\ell)})\in \R^{\fd(\ell)}$, $x=(x_1,\dots,x_d)\in [a,b]^d$, $i\in \{1,2,\dots,\ell_1\}$ it holds that
    \begin{equation}\llabel{base case}
    \begin{split}
        &|\Relii{\ell}{1}{\theta}{\bbA}{i}(x)|=\big|\theta_{\ell_1\ell_0+i}+\textstyle\sum_{j=0}^d\theta_{(i-1)\ell_0+j}x_j\bigr|\leq |\theta_{\ell_1\ell_0+i}|+\bigl|\textstyle\sum_{j=0}^d\theta_{(i-1)\ell_0+j}x_j\bigr|\\
        &\textstyle \leq (d+1)\bigl[\max\{1,|a|,|b|\}\max_{j\in \{1,2,\dots,\bfd_1\}}|\theta_j|\bigr] \\
        &\textstyle\leq [\max\{1,|a|,|b|\}]^{\max\{1,\polycons^0\}}\bigl[\max\{1,\max_{j\in \{1,2,\dots,\bfd_{0+1}\}}|\theta_j|\}\bigr]^{\max\{1,\polycons^0\}}\\
&\textstyle\quad\cdot\bigl[1+2\max\{c,1\}\ell_0\bigr]^{\max\{1,\polycons^k\}}.
        \end{split}
    \end{equation}
    This proves \lref{need to prove} in the base case $k=1$. For the induction step we assume that there exists $k\in \N\cap[0,L)$ which satisfies for all $\theta=(\theta_1,\dots,\theta_{\fd(\ell)})\in \R^{\fd(\ell)}$, $x\in [a,b]^d$ that
    \begin{equation}\llabel{assume}
        \begin{split}
        &\textstyle\max_{i\in\{1,2,\dots,\ell_{k}\}}|\Relii{\ell}{k}{\theta}{\bbA}{i}(x)|\\
        &\textstyle\leq [\max\{1,|a|,|b|\}]^{\max\{1,\polycons^{k-1}\}}\bigl[\max\{1,\max_{j\in \{1,2,\dots,\bfd_{k}\}}|\theta_j|\}\bigr]^{k\max\{1,\polycons^{k-1}\}}\\
&\textstyle\quad\cdot\bigl[1+2\max\{c,1\}\max\{\ell_0,\ell_1,\dots\ell_{k}\}\bigr]^{k\max\{1,\polycons^{k-1}\}}.
        \end{split}
    \end{equation}
    \startnewargseq
    \argument{\cref{setting ANN: realization multi};}{that for all $\theta=(\theta_1,\dots,\theta_{\fd(\ell)})\in \R^{\fd(\ell)}$, $x\in [a,b]^d$, $i\in\{1,2,\dots,\ell_{k+1}\}$ it holds that
    \begin{equation}\llabel{eq1}
    \begin{split}
        |\Relii{\ell}{k+1}{\theta}{\bbA}{i}(x)|&\leq \textstyle |\theta_{\ell_{k+1}\ell_{k}+i+\bfd_{k}}|+\bigl[\sum_{j=1}^{\ell_{k}}|\theta_{(i-1)\ell_{k}+j+\bfd_{k}}\bbA(\Relii{\ell}{k}{\theta}{\bbA}{j}(x))|\bigr]\\
        &\leq\textstyle \bigl[\max_{j\in\{1,2,\dots,\bfd_{k+1}\}}|\theta_j|\bigr]\\
        &+\textstyle\ell_{k}\bigl[\max_{j\in\{1,2,\dots,\bfd_{k+1}\}}|\theta_j|\bigr]\bigl[\max_{j\in\{1,2,\dots,\ell_{k}\}}|\bbA(\Relii{\ell}{k}{\theta}{\bbA}{j}(x))|\bigr]
        .
        \end{split}
    \end{equation}}
    \argument{\lref{eq1};the assumption that for all $x\in \R$ it holds that $|\bbA(x)|\leq c(1+|x|^\polycons)$}{that for all $\theta=(\theta_1,\dots,\theta_{\fd(\ell)})\in \R^{\fd(\ell)}$, $x\in [a,b]^d$, $i\in\{1,2,\dots,\ell_{k+1}\}$ it holds that
      \begin{equation}\llabel{eq2}
    \begin{split}
       & |\Relii{\ell}{k+1}{\theta}{\bbA}{i}(x)|
        \leq\textstyle \bigl[\max_{j\in\{1,2,\dots,\bfd_{k+1}\}}|\theta_j|\bigr]\\
        &+\textstyle\ell_{k}\bigl[\max_{j\in\{1,2,\dots,\bfd_{k+1}\}}|\theta_j|\bigr]\bigl[\max_{j\in\{1,2,\dots,\ell_{k}\}}\bigl(c+c|\Relii{\ell}{k}{\theta}{\bbA}{j}(x)|^\polycons\bigr)\bigr]\\
        &=\textstyle\bigl[\max_{j\in\{1,2,\dots,\bfd_{k+1}\}}|\theta_j|\bigr]\bigl(1+\ell_{k}c+\ell_{k}c\bigl[\max_{j\in\{1,2,\dots,\ell_{k}\}}|\Relii{\ell}{k}{\theta}{\bbA}{j}(x)|^\polycons\bigr]\bigr).
        \end{split}
    \end{equation}}
    \argument{\lref{eq2};\lref{assume}}{for all $\theta=(\theta_1,\dots,\theta_{\fd(\ell)})\in \R^{\fd(\ell)}$, $x\in [a,b]^d$, $i\in\{1,2,\dots,\ell_{k+1}\}$ that
    \begin{equation}\llabel{eq3}
         \begin{split}
        |\Relii{\ell}{k+1}{\theta}{\bbA}{i}(x)|
        &\leq\textstyle (1+\ell_{k}c)\bigl[\max_{j\in\{1,2,\dots,\bfd_{k+1}\}}|\theta_j|\bigr]\\
        &+\textstyle\ell_{k}c\bigl[\max_{j\in\{1,2,\dots,\bfd_{k+1}\}}|\theta_j|\bigr][\max\{1,|a|,|b|\}]^{\polycons\max\{1,\polycons^{k-1}\}}\\
        &\textstyle \quad\cdot\bigl[\max\{1,\max_{j\in \{1,2,\dots,\bfd_{k}\}}|\theta_j|\}\bigr]^{\polycons k\max\{1,\polycons^{k-1}\}}\\
&\textstyle\quad\cdot\bigl[1+2\max\{c,1\}\max\{\ell_1,\ell_2,\dots\ell_{k-1}\}\bigr]^{\polycons k\max\{1,\polycons^{k-1}\}}\\
&\leq \textstyle (1+\ell_{k}c)\bigl[\max_{j\in\{1,2,\dots,\bfd_{k+1}\}}|\theta_j|\bigr]\\
        &+\textstyle\ell_{k}c\bigl[\max_{j\in\{1,2,\dots,\bfd_{k+1}\}}|\theta_j|\bigr][\max\{1,|a|,|b|\}]^{\max\{1,\polycons^{k}\}}\\
        &\textstyle \quad\cdot\bigl[\max\{1,\max_{j\in \{1,2,\dots,\bfd_{k+1}\}}|\theta_j|\}\bigr]^{k\max\{1,\polycons^{k}\}}\\
&\textstyle\quad\cdot\bigl[1+2\max\{c,1\}\max\{\ell_1,\ell_2,\dots\ell_{k}\}\bigr]^{k\max\{1,\polycons^{k}\}}.
        \end{split}
    \end{equation}}
    \argument{\lref{eq3};}{for all $\theta=(\theta_1,\dots,\theta_{\fd(\ell)})\in \R^{\fd(\ell)}$, $x\in [a,b]^d$, $i\in\{1,2,\dots,\ell_{k+1}\}$ that
    \begin{equation}\llabel{eq4}
         \begin{split}
        |\Relii{\ell}{k+1}{\theta}{\bbA}{j}(x)|
&\leq 
        \textstyle(1+2\ell_{k}c)\bigl[\max_{j\in\{1,2,\dots,\bfd_{k+1}\}}|\theta_j|\bigr][\max\{1,|a|,|b|\}]^{\max\{1,\polycons^{k}\}}\\
        &\textstyle \quad\cdot\bigl[\max\{1,\max_{j\in \{1,2,\dots,\bfd_{k+1}\}}|\theta_j|\}\bigr]^{k\max\{1,\polycons^{k}\}}\\
&\textstyle\quad\cdot\bigl[1+2\max\{c,1\}\textstyle\max\{\ell_1,\ell_2,\dots\ell_{k}\}\bigr]^{k\max\{1,\polycons^{k}\}}\\
&\textstyle\leq [\max\{1,|a|,|b|\}]^{\max\{1,\polycons^{k}\}}\bigl[\max\{1,\max_{j\in \{1,2,\dots,\bfd_{k+1}\}}|\theta_j|\}\bigr]^{(k+1)\max\{1,\polycons^{k}\}}\\
&\quad\textstyle \cdot\bigl[1+2\max\{c,1\}\max\{\ell_1,\ell_2,\dots\ell_{k}\}\bigr]^{(k+1)\max\{1,\polycons^{k}\}}.
        \end{split}
    \end{equation}}
    \argument{\lref{eq4};\lref{base case};induction}{\lref{need to prove}\dott}
    \end{aproof}
    \subsection{Explicit a priori estimates for gradients of realization functions of ANNs}
    
    In this subsection we employ the notation from \cref{setting: ANN} above to establish in \cref{lem: ann lipschitz} below an elementary a priori estimate for the distance between the realizations of two \ANNs\ with different \ANN\ parameter vectors.
    \begin{athm}{lemma}{lem: ann lipschitz}
Assume \cref{setting: ANN} and let $c,a\in \R$, $b\in (a,\infty)$ satisfy for all $x,y\in \R$ that $|\bbA(x)-\bbA(y)|\leq c|x-y|$. Then it holds for all $\theta=(\theta_1,\dots,\theta_{\fd(\ell)})$, $\vartheta=(\vartheta_1,\dots,\vartheta_{\fd(\ell)})\in \R^{\fd(\ell)}$,  $x\in [a,b]^d$ that
    \begin{equation}\llabel{conclude}
  \begin{split}
 |  \Reli{ \ell }{ L }{ \theta }{ \Act }( x )-\Reli{ \ell }{ L }{ \vartheta }{ \Act }( x )|
 &\leq \textstyle\max\{1,|a|,|b|\}\bigl[(3+|\bbA(0)|)\max\{c,1\}\max\{\ell_0,\ell_1,\dots,\ell_{L-1}\}\bigr]^{L}\\
&\textstyle\cdot\bigl[\max_{j\in \{1,2,\dots,\fd(\ell)\}}\max\{1,|\theta_j|,|\vartheta_j|\}\bigr]^{L}\bigl[\max_{j\in \{1,2,\dots,\fd(\ell)\}}|\theta_j-\vartheta_j|\bigr]\ifnocf.
\end{split}
    \end{equation}
    \cfout[.]        
    \end{athm}
    \begin{aproof}
     Throughout this proof, for every  $k\in \{0,1,\dots,L\}$ let $\bfd_k\in \R$ satisfy $\bfd_k=\sum_{m=1}^k\ell_m(\ell_{m-1}+1)$.
    We claim that for all $k\in \{ 1,2,\dots,L\}$, $\theta=(\theta_1,\dots,\theta_{\fd(\ell)})$, $\vartheta=(\vartheta_1,\dots,\vartheta_{\fd(\ell)})\in \R^{\fd(\ell)}$,  $x\in [a,b]^d$ it holds that
    \begin{equation}\llabel{need to prove}
    \begin{split}
        &\textstyle\max_{i\in\{1,2,\dots,\ell_{k}\}}|\Relii{\ell}{k}{\theta}{\bbA}{i}(x)-\Relii{\ell}{k}{\vartheta}{\bbA}{i}(x)|\\
        &\textstyle\leq \max\{1,|a|,|b|\}\bigl[\max_{j\in \{1,2,\dots,\bfd_{k}\}}\max\{1,|\theta_j|,|\vartheta_j|\}\bigr]^{k}\\
&\textstyle\quad\cdot\bigl[1+(2+|\bbA(0)|)\max\{c,1\}\max\{\ell_0,\ell_1,\dots,\ell_{k-1}\}\bigr]^{k}\big[\max_{j\in \{1,2,\dots,\bfd_{k}\}}|\theta_j-\vartheta_j|\bigr].
        \end{split}
    \end{equation}
    We prove \lref{need to prove} by induction on $k\in \{1,2,\dots,L\}$. For the base case $k=1$ note that \cref{realization multi} shows that for all $\theta=(\theta_1,\dots,\theta_{\fd(\ell)})$, $\vartheta=(\vartheta_1,\dots,\vartheta_{\fd(\ell)})\in \R^{\fd(\ell)}$, $x=(x_1,\dots,x_d)\in [a,b]^d$, $i\in \{1,2,\dots,\ell_1\}$ it holds that
    \begin{equation}\llabel{base case}
    \begin{split}
        |\Relii{\ell}{1}{\theta}{\bbA}{i}(x)-\Relii{\ell}{1}{\vartheta}{\bbA}{i}(x)|
&=\big|\theta_{\ell_1\ell_0+i}+\textstyle\sum_{j=0}^d\theta_{(i-1)\ell_0+j}x_j-\bigl(\vartheta_{\ell_1\ell_0+i}+\textstyle\sum_{j=0}^d\vartheta_{(i-1)\ell_0+j}x_j\bigr)\big|\\
        &\leq |\theta_{\ell_1\ell_0+i}-\vartheta_{\ell_1\ell_0+i}|+\big|\textstyle\sum_{j=0}^d(\theta_{(i-1)\ell_0+j}-\theta_{(i-1)\ell_0+j})x_j\big|\\
        &\textstyle \leq (d+1)[\max\{1,|a|,|b|\}]\bigl[\max_{j\in \{1,2,\dots,\bfd_1\}}|\theta_j-\vartheta_j|\bigr] \\
        &\textstyle\leq [\max\{1,|a|,|b|\}]\bigl[\max_{j\in \{1,2,\dots,\bfd_{0+1}\}}\max\{1,|\theta_j|,|\vartheta_j|\}\bigr]^{1}\\
&\textstyle\quad\cdot[1+(2+|\bbA(0)|)\max\{c,1\}\ell_0]^{1}\big[\max_{j\in \{1,2,\dots,\bfd_1\}}|\theta_j-\vartheta_j|\bigr].
        \end{split}
    \end{equation}
    This proves \lref{need to prove} in the base case $k=1$. For the induction step we assume that there exists $k\in\N\cap[0,L)$ which satisfies for all $\theta=(\theta_1,\dots,\theta_{\fd(\ell)})$, $\vartheta=(\vartheta_1,\dots,\vartheta_{\fd(\ell)})\in \R^{\fd(\ell)}$, $x\in [a,b]^d$ that
    \begin{equation}\llabel{assume}
    \begin{split}
       &\textstyle\max_{i\in\{1,2,\dots,\ell_{k}\}}|\Relii{\ell}{k}{\theta}{\bbA}{i}(x)-\Relii{\ell}{k}{\vartheta}{\bbA}{i}(x)|\\
        &\textstyle\leq \max\{1,|a|,|b|\}\bigl[\max_{j\in \{1,2,\dots,\bfd_{k}\}}\max\{1,|\theta_j|,|\vartheta_j|\}\bigr]^{k}\\
&\textstyle\quad\cdot\bigl[1+(2+|\bbA(0)|)\max\{c,1\}\max\{\ell_0,\ell_1,\dots,\ell_{k-1}\}\bigr]^{k}\big[\max_{j\in \{1,2,\dots,\bfd_{k}\}}|\theta_j-\vartheta_j|\bigr].
\end{split}
    \end{equation}
    \startnewargseq
    \argument{\cref{setting ANN: realization multi};}{that for all $\theta=(\theta_1,\dots,\theta_{\fd(\ell)})$, $\vartheta=(\vartheta_1,\dots,\vartheta_{\fd(\ell)})\in \R^{\fd(\ell)}$, $x\in [a,b]^d$, $i\in\{1,2,\dots,\ell_{k+1}\}$ it holds that
    \begin{equation}\llabel{eq1}
    \begin{split}
        &|\Relii{\ell}{k+1}{\theta}{\bbA}{i}(x)-\Relii{\ell}{k+1}{\vartheta}{\bbA}{i}(x)|\\
        &\leq \textstyle |\theta_{\ell_{k+1}\ell_{k}+i+\bfd_{k}}-\vartheta_{\ell_{k+1}\ell_{k}+i+\bfd_{k}}|\\
        &+\textstyle\bigl[\sum_{j=1}^{\ell_{k}}|\theta_{(i-1)\ell_{k}+j+\bfd_{k}}\bbA(\Relii{\ell}{k}{\theta}{\bbA}{j}(x))-\vartheta_{(i-1)\ell_{k}+j+\bfd_{k}}\bbA(\Relii{\ell}{k}{\vartheta}{\bbA}{j}(x))|\bigr]\\
        &\leq\textstyle \bigl[\max_{j\in\{1,2,\dots,\bfd_{k+1}\}}|\theta_j-\vartheta_j|\bigr]\\
        &+\textstyle\ell_{k}\bigl[\max_{j\in\{1,2,\dots,\bfd_{k+1}\}}|\theta_j|\bigr]\bigl[\max_{j\in\{1,2,\dots,\ell_{k}\}}|\bbA(\Relii{\ell}{k}{\theta}{\bbA}{j}(x))-\bbA(\Relii{\ell}{k}{\vartheta}{\bbA}{j}(x))|\bigr]\\
        &+\textstyle \ell_{k}\bigl[\max_{j\in\{1,2,\dots,\bfd_{k+1}\}}|\theta_j-\vartheta_j|\bigr]\bigl[\max_{j\in\{1,2,\dots,\ell_{k}\}}|\bbA(\Relii{\ell}{k}{\vartheta}{\bbA}{j}(x))|\bigr]
        .
        \end{split}
    \end{equation}}
    \argument{\lref{eq1};the assumption that for all $x,y\in \R$ it holds that $|\bbA(x)-\bbA(y)|\leq c|x-y|$}{that for all $i\in\{1,2,\dots,\ell_{k+1}\}$, $\theta=(\theta_1,\dots,\theta_{\fd(\ell)})$, $\vartheta=(\vartheta_1,\dots,\vartheta_{\fd(\ell)})\in \R^{\fd(\ell)}$, $x\in [a,b]^d$ it holds that
      \begin{equation}\llabel{eq2}
   \begin{split}
        &|\Relii{\ell}{k+1}{\theta}{\bbA}{i}(x)-\Relii{\ell}{k+1}{\vartheta}{\bbA}{i}(x)|\\
        &\leq\textstyle \bigl[\max_{j\in\{1,2,\dots,\bfd_{k+1}\}}|\theta_j-\vartheta_j|\bigr]\\
        &+\textstyle\ell_{k}c\bigl[\max_{j\in\{1,2,\dots,\bfd_{k+1}\}}|\theta_j|\bigr]\bigl[\max_{j\in\{1,2,\dots,\ell_{k}\}}|\Relii{\ell}{k}{\theta}{\bbA}{j}(x)-\Relii{\ell}{k}{\vartheta}{\bbA}{j}(x)|\bigr]\\
        &+\textstyle \ell_{k}|\bbA(0)|\bigl[\max_{j\in\{1,2,\dots,\bfd_{k+1}\}}|\theta_j-\vartheta_j|\bigr]\\
        &+\textstyle \ell_{k}c\bigl[\max_{j\in\{1,2,\dots,\bfd_{k+1}\}}|\theta_j-\vartheta_j|\bigr]\bigl[\max_{j\in\{1,2,\dots,\ell_{k}\}}|\Relii{\ell}{k}{\vartheta}{\bbA}{j}(x)|\bigr]
        .
        \end{split}
    \end{equation}}
    \argument{\lref{eq2};\lref{assume};\cref{lem: anngrowth}}{for all $i\in\{1,2,\dots,\ell_{k+1}\}$, $\theta=(\theta_1,\dots,\theta_{\fd(\ell)})$, $\vartheta=(\vartheta_1,\dots,\vartheta_{\fd(\ell)})\allowbreak\in \R^{\fd(\ell)}$, $x\in [a,b]^d$ that
    \begin{equation}\llabel{eq3}
         \begin{split}
       & |\Relii{\ell}{k+1}{\theta}{\bbA}{i}(x)-\Relii{\ell}{k+1}{\vartheta}{\bbA}{i}(x)|\\
       &\textstyle\leq (1+\ell_{k}|\bbA(0)|)\bigl[\max_{j\in\{1,2,\dots,\bfd_{k+1}\}}|\theta_j-\vartheta_j|\bigr]\\
       &+\textstyle\ell_{k}c\bigl[\max_{j\in\{1,2,\dots,\bfd_{k+1}\}}|\theta_j|\bigr]\max\{1,|a|,|b|\} \textstyle\bigl[\max_{j\in \{1,2,\dots,\bfd_{k}\}}\max\{1,|\theta_j|,|\vartheta_j|\}\bigr]^{k}\\
&\textstyle\quad\cdot\bigl[\max_{l\in \{0,1,\dots,k-1\}}(1+(2+|\bbA(0)|)\max\{c,1\}\ell_l)\bigr]^{k}\big[\max_{j\in \{1,2,\dots,\bfd_{k}\}}|\theta_j-\vartheta_j|\bigr]\\
       &+\textstyle\ell_{k}c\bigl[\max_{j\in\{1,2,\dots,\bfd_{k+1}\}}|\theta_j-\vartheta_j|\bigr]\textstyle[\max\{1,|a|,|b|\}]\\
       &\textstyle\quad \cdot\bigl[\max\{1,\max_{j\in \{1,2,\dots,\fd(\ell)\}}|\vartheta_j|\}\bigr]^{k}\bigl[1+2\max\{c,1\}\max\{\ell_0,\ell_1,\dots,\ell_{L-1}\}\bigr]^{k}.
        \end{split}
    \end{equation}}
    \argument{\lref{eq3};}{for all $i\in\{1,2,\dots,\ell_{k+1}\}$, $\theta=(\theta_1,\dots,\theta_{\fd(\ell)})$, $\vartheta=(\vartheta_1,\dots,\vartheta_{\fd(\ell)})\in \R^{\fd(\ell)}$, $x\in [a,b]^d$ that
    \begin{equation}\llabel{eq4}
         \begin{split}
      &|\Relii{\ell}{k+1}{\theta}{\bbA}{i}(x)-\Relii{\ell}{k+1}{\vartheta}{\bbA}{i}(x)|
      \\
&\leq\textstyle \bigl[\max_{j\in\{1,2,\dots,\bfd_{k+1}\}}|\theta_j-\vartheta_j|\bigr][\max\{1,|a|,|b|\}] \bigl[\max_{j\in \{1,2,\dots,\bfd_{k}\}}\max\{1,|\theta_j|,|\vartheta_j|\}\bigr]^{k}\\
&\textstyle\quad \cdot \bigl[1+(2+|\bbA(0)|)\max\{c,1\}\max\{\ell_0,\ell_1,\dots,\ell_{k-1}\}\bigr]^{k}\\
&\textstyle\quad \cdot\Bigl[1+\ell_{k}|\bbA(0)|+\ell_{k}c\bigl[\max_{j\in\{1,2,\dots,\bfd_{k+1}\}}|\theta_j|\bigr]+\ell_{k}c\Bigr]\\
&\leq\textstyle \bigl[\max_{j\in\{1,2,\dots,\bfd_{k+1}\}}|\theta_j-\vartheta_j|\bigr][\max\{1,|a|,|b|\}]\bigl[\max_{j\in \{1,2,\dots,\bfd_{k}\}}\max\{1,|\theta_j|,|\vartheta_j|\}\bigr]^{k}\\
&\textstyle\quad \cdot \bigl[1+(2+|\bbA(0)|)\max\{c,1\}\max\{\ell_0,\ell_1,\dots,\ell_{k-1}\}\bigr]^{k}\\
&\textstyle\quad \cdot\bigl[1+\ell_{k}|\bbA(0)|+2\ell_{k}c\bigr]\bigl[\max\{1,\max_{j\in\{1,2,\dots,\bfd_{k+1}\}}|\theta_j|\}\bigr]\\
&\leq\textstyle \bigl[\max_{j\in\{1,2,\dots,\bfd_{k+1}\}}|\theta_j-\vartheta_j|\bigr][\max\{1,|a|,|b|\}]\bigl[\max_{j\in \{1,2,\dots,\bfd_{k}\}}\max\{1,|\theta_j|,|\vartheta_j|\}\bigr]^{k+1}\\
&\textstyle\quad \cdot \bigl[1+(2+|\bbA(0)|)\max\{c,1\}\max\{\ell_0,\ell_1,\dots,\ell_{k}\}\bigr]^{k+1}.
        \end{split}
    \end{equation}}
    \argument{\lref{eq4};\lref{base case};induction}{\lref{need to prove}\dott}
    \startnewargseq
    \argument{\lref{need to prove};}{\lref{conclude}\dott}
    \end{aproof}
    
In \cref{cor: ann lipschitz} below we employ the notation from \cref{setting: ANN} and the conclusion from \cref{lem: ann lipschitz}  to establish an elementary a priori estimate for the evaluation of the gradient (with respect to the \ANN\ parameter vector) of the realization of an \ANN\ with the general activation function $\bbA \colon \R\to \R$ from \cref{setting: ANN}.

    \begin{athm}{cor}{cor: ann lipschitz}
       Assume \cref{setting: ANN}, let $a\in \R$, $b\in (a,\infty)$, and assume $\bbA\in C^1(\R,\R)$. Then it holds\footnote{Note that for all $d \in \N$, $v = ( v_1, \dots, v_d ) \in \R^d$ it holds that $\| v \| = [\sum_{ i = 1 }^d | v_i |^2]^{ 1 / 2 }$ (standard norm).} for all $\theta=(\theta_1,\dots,\theta_{\fd(\ell)})\in \R^{\fd(\ell)}$,  $x\in [a,b]^d$ that
    \begin{equation}\llabel{conclude}
  \begin{split}
&\| \nabla_\theta( \Reli{ \ell }{ L }{ \theta }{ \Act }( x ))\|
 \leq \textstyle\max\{1,|a|,|b|\}(3+|\bbA(0)|)^L\\
 &\textstyle\cdot\bigl[\max\{\ell_0,\ell_1,\dots,\ell_{L-1}\}\max\{1,\max_{j\in \{1,2,\dots,\fd(\ell)\}}|\theta_j|\}\max\{1,\sup_{x\in \R}|\bbA'(x)|\}\bigr]^{L}\ifnocf.
\end{split}
    \end{equation}
    \cfout[.]        
    \end{athm}
    \begin{aproof}
    Throughout this proof, assume without loss of generality that $\sup_{ x \in \R } \allowbreak | \bbA'(x) | < \infty$ (otherwise note that the fact that for all $\theta \in \R^\fd$, $x \in [a,b]^d$ it holds that $ \| \nabla_\theta( \Reli{ \ell }{ L }{ \theta }{ \Act }( x ))\|\leq \infty$
    establishes \lref{conclude}).
\argument{\cref{lem: ann lipschitz};}{\lref{conclude}\dott}
    \end{aproof}
\section{Error estimates for ANN approximations}\label{sec: ANN approximation}
In this section we establish in \cref{cor: Pinns for Poisson3 intro2 pre2 approximation} below basic approximation error estimates for deep \ANN\ approximations of twice continuously differentiable functions. In \cref{sec: deep Kolmogorov method} we employ \cref{cor: Pinns for Poisson3 intro2 pre2 approximation} together with \cref{lem: anngrowth} and \cref{cor: ann lipschitz} from \cref{sec: priori estimate} to establish in \cref{prop: linear II relu pinn draft}, \cref{prop: linear II relu pinn}, \cref{cor: linear II relu pinn}, and \cref{main cor} error estimates for the \DKM\ and, thereby, prove \cref{main theorem} in the introduction.

In our proof of \cref{cor: Pinns for Poisson3 intro2 pre2 approximation} we also employ the well-known fact in \cref{lem: existence of smooth} below that a constant function on a compact (or even only closed) set can be extended to a smooth function with the support in an open set containing the compact set (cf., \eg, \cite[Theorem 2.29 in Chapter~2]{MR2954043}). Only for completeness we include here a detailed proof for \cref{lem: existence of smooth}.


Our proof of the \ANN\ approximation result with convergence rates in \cref{cor: Pinns for Poisson3 intro2 pre2 approximation} strongly builds up on the findings in \cite{MR4626415}. We also refer, \eg, to \cite{MR4319100,Beckfullerror2019,Yarotskyerrorbound,MR4131039,Petersenoptimal2018,MR4505885,Ackermann2024} and the references therein for further \ANN\ approximation results.
\subsection{Error estimates for shallow ANN approximations for 
$C^2$-functions}\label{subsec: estimate shallow C2}
\cfclear
In \cref{conj: Pinns for Poisson3 intro2 pre2 approximation shallow ann} and \cref{conj: Pinns for Poisson3 intro2 pre2 approximation shallow ann general} in this subsection we first establish variants of the \ANN\ approximation result in \cref{cor: Pinns for Poisson3 intro2 pre2 approximation} below but for shallow \ANNs\ instead of deep \ANNs\ as in \cref{cor: Pinns for Poisson3 intro2 pre2 approximation}. In \cref{subsec: error estimate C2} below we extend shallow \ANNs\ to deep \ANNs\ with the same realization functions (see \cref{lem: embedding ANN} below) and combine this with \cref{conj: Pinns for Poisson3 intro2 pre2 approximation shallow ann general} to establish \cref{cor: Pinns for Poisson3 intro2 pre2 approximation}. \cref{conj: Pinns for Poisson3 intro2 pre2 approximation shallow ann} is strongly based on \cite{MR4626415}.
\begin{athm}{lemma}{conj: Pinns for Poisson3 intro2 pre2 approximation shallow ann}
    Let $d \in \N$, $u \in C^2( [-1,1]^d, \R )$. Then there exists $ c \in (0,\infty)$ such that for all $\ell=(\ell_0,\ell_1,\ell_2)\allowbreak\in (\{d\}\times\N\times\{1\}) $ there exists $\theta \in ( - c, c )^{ \fd(\ell) }$ such that
    \begin{equation}\llabel{conclude}
\sup_{x\in (-1,1)^d}
  | u( x)
 -\textstyle
  \reli{\ell}{2}{\theta}( x )
| \leq
  c (\ell_1)^{ \frac{-1}{d+5} } \ifnocf.
    \end{equation}
    \cfout[.]
\end{athm}
\begin{aproof}
    Throughout this proof, 
     for every $\ell=(\ell_0,\ell_1,\ell_2) \in \N^3 $ let $N_\ell,\fL_\ell\in \Z$ satisfy
    \begin{equation}\llabel{def: N}
        N_\ell=\lfloor(\ell_1)^{\frac{d+2}{d(d+4)}}\rfloor\qqandqq \fL_\ell=\ceil{\log_2(N_\ell)},
    \end{equation}
   let $K_{\ell}\colon \R\to\R$ satisfy for all $t\in \R$ that
    \begin{equation}\llabel{def: K}
        K_{\ell}(t)=\biggl[\int_{-\pi}^\pi\Bigl(\frac{\sin(\frac{(N_\ell+1)s}{2})}{\sin(\frac{s}{2})}\Bigr)^2\,\d s\biggr]^{-1}\biggl(\frac{\sin(\frac{(N_\ell+1)t}{2})}{\sin(\frac{t}{2})}\biggr)^2,
    \end{equation}
    let $(\bfa_{k}^\ell)_{k\in \Z}\subseteq \R$ satisfy for all $t\in \R$ that
    \begin{equation}\llabel{def: bfa}
        K_{\ell}(t)=\sum_{k\in \Z}\bfa_{k}^\ell \exp(\imag kt),
        \end{equation}
        let $(a_{k}^\ell)_{k\in\Z^d}\subseteq\R$ satisfy for all $k=(k_1,\dots,k_d)\in \Z^d$ that
        \begin{equation}\llabel{def: a}
            a_{k}^\ell=\prod_{j=1}^d \bfa_{k_j}^\ell,
        \end{equation} 
    let $G_{\ell}\colon \R^d\to\mathbb C$ satisfy for all $x=(x_1,\dots,x_d)\in \R^d$ that
    \begin{equation}\llabel{def: G}
         G_{\ell}(x)=\sum_{k=(k_1,\dots,k_d)\in \Z^d}a_k^\ell\exp\bigl(\imag\textstyle\sum_{i=1}^dk_ix_i\bigr),
    \end{equation}
    (cf.\ (3.6) in Section 3 in \cite{MR4626415})
    and let $J_{\ell}\colon C([-\pi,\pi]^d,\R)\times\R^d\to\mathbb C$ satisfy for all $g\in C([-\pi,\pi]^d,\R)$, $x=(x_1,\dots,x_d)\in \R^d$ that
    \begin{equation}\llabel{def: J}
J_\ell(g,x)=\int_{[-\pi,\pi]^d}\textstyle G_\ell(x-y)g(y)\,\d y
    \end{equation}
      (cf.\ (3.7) in Section 3 in \cite {MR4626415}).
     \argument{\cite[Chapter 6]{Steinsingular}}{that there exist $F\in C^1([-\pi,\pi]^d,\R)$, $\fC\in \R$ which satisfy that
    \begin{enumerate}[label=(\roman*)]
        \item \llabel{item 1} it holds for all $x\in (-1,1)^d$ that $F(x)=u(x)$,
        \item \llabel{item 2} it holds for all $i\in \{1,2,\dots,d\}$, $x=(x_1,\dots,x_{d})\in [-\pi,\pi]^{d}$ with $x_i=-\pi$ that $F(x)=F(x_1,x_2,\dots,x_{i-1},x_i+2\pi,x_{i+1},\dots,x_{d})$, and
        \item \llabel{item 3} it holds that
        \begin{equation}
            \begin{split}
            &\textstyle\max\bigl\{\sup_{x\in [-\pi,\pi]^d} |(F(x)|,\sup_{x\in [-\pi,\pi]^d} \|(\nabla F)(x)\|\bigr\}\\
            &\leq \textstyle \fC\bigl[\max\bigl\{\sup_{x\in (-1,1)^d}|u(x)|,\sup_{x\in (-1,1)^d}\|(\nabla u)(x)\|\bigr\}\bigr].
            \end{split}
        \end{equation}
    \end{enumerate}}
\startnewargseq
    \argument{(3.7) in Section 3 in \cite{MR4626415};}{that for all $\ell\in \N^3$, $f\in C([-\pi,\pi]^d,\R)$ it holds that 
    \begin{equation}\llabel{eqt2}
         ([-\pi,\pi]^d \ni x\mapsto J_{\ell}(f,x)\in \mathbb{C})\in C([-\pi,\pi]^d,\mathbb{C}).
    \end{equation}}
    \startnewargseq
    In the following 
    for every $\ell\in \N^3$, $k=(k_1,\dots,k_d)\in \Z^d$ let $\fJ_\ell(k)\in \mathbb C$ satisfy
    \begin{equation}
        \fJ_\ell(k)=\frac{1}{(2\pi)^d}\int_{[-\pi,\pi]^d} J_\ell(F,y_1,y_2,\dots,y_d)\exp\textstyle\bigl(-\imag \sum_{i=1}^dk_iy_i\bigr)\,\d y_1\d y_2\dots\d y_d
    \end{equation}
     and for every $\ell\in \N^3$ let $v_\ell\in [0,\infty]$ satisfy
    \begin{equation}\llabel{def: v}
        v_\ell=\sum_{k=(k_1,\dots,k_d)\in \Z^d}|\fJ_\ell(k)|( |k_1|+|k_2|+\ldots+|k_d|)^2.
    \end{equation}
    \startnewargseq
    \argument{(3.8) in Section 3 in \cite{MR4626415};\lref{def: N};\lref{def: K};\lref{def: J}}{that there exists $\rho\in (0,\infty)$ which satisfies for all $\ell\in (\{d\}\times\N\times\{1\})  $ that
    \begin{equation}\llabel{def: C}
        \sup_{t\in [-\pi,\pi]^d}|J_\ell(F,t)-F(t)|\leq 
\rho(N_\ell)^{-1}\max\biggl\{\sup_{x\in [-\pi,\pi]^d}|F(x)|,\sup_{x\in [-\pi,\pi]^d}\|(\nabla F)(x)\|\biggr\}.
    \end{equation}}
    \startnewargseq
    \argument{\cite[Lemma 1]{MR4626415};\lref{def: v}}{that there exist $\theta_\ell=(\theta_\ell^1,\dots,\theta_\ell^{\fd(\ell)})\in\R^{\fd(\ell)}$, $\ell\in (\{d\}\times\N\times\{1\})$, and $\varrho\in (0,\infty)$ which satisfy for all $\ell=(\ell_0,\ell_1,\ell_2)\in (\{d\}\times\N\times\{1\})  $ that
    \begin{equation}\llabel{def: theta}
        \sup_{i\in \{1,2,\dots,\fd(\ell)\}}|\theta_\ell^i|\leq \max\Bigl\{1,\frac{8\pi^2v_\ell}{\ell_1}\Bigr\} 
                \end{equation}
        \begin{equation}\llabel{def: fC}
            \text{and}\qquad \sup_{x\in (-1,1)^d}|J_\ell(F,x)-\reli{\ell}{2}{\theta_\ell}(x)|\leq \varrho v_\ell\max\{\sqrt{\ln(\ell_1)},\sqrt{d}\}(\ell_1)^{-\frac 12-\frac 1d}
    \end{equation}}
    \startnewargseq
    \argument{\lref{def: N};\lref{def: C};\lref{def: fC};\lref{item 1}}{that for all $\ell=(\ell_0,\ell_1,\ell_2)\in (\{d\}\times\N\times\{1\})$ it holds that
    \begin{equation}\llabel{eq1}
        \begin{split}
      &  \textstyle\sup_{x\in (-1,1)^d}|u(x)-\reli{\ell}{2}{\theta_\ell}(x)|\\
      &\leq \textstyle
\rho(N_\ell)^{-1} \max\bigl\{\sup_{x\in [-\pi,\pi]^d}|F(x)|,\sup_{x\in[-\pi,\pi]^d}\|(\nabla F)(x)\|\bigr\}\\
&+\varrho v_\ell\max\{\sqrt{\ln(\ell_1)},\sqrt{d}\}(\ell_1)^{-\frac 12-\frac 1d}\\
      &\textstyle\leq 
2\rho (\ell_1)^{-\frac{d+2}{d(d+4)}} \max\bigl\{\sup_{x\in [-\pi,\pi]^d}|F(x)|,\sup_{x\in [-\pi,\pi]^d}\|(\nabla F)(x)\|\bigr\}\\
&+\varrho v_\ell\max\{\sqrt{\ln(\ell_1)},\sqrt{d}\}(\ell_1)^{-\frac 12-\frac 1d}.
        \end{split}
    \end{equation}}
     \argument{(4.1) in Section 4 in \cite{MR4626415}; (4.2) in Section 4 in \cite{MR4626415}; (4.6) in Section 4 in \cite{MR4626415};\lref{def: N}}{that there exist $c_1,c_2\in (0,\infty)$ such that for all $\ell\in (\{d\}\times\N\times\{1\})$ it holds that
    \begin{equation}\llabel{Eq1}
    \begin{split}
       \textstyle v_\ell \leq 2d(4\pi c_1)^dc_2\bigl[\sup_{x\in [-\pi,\pi]^d}\max\{|F(x)|,\|(\nabla F)(x)\|\}\bigr]2^{-1}3^{d/2}\sum_{i=0}^{\fL_\ell}(i+1)^d(2^i)^{1+\frac d2}.
       \end{split}
    \end{equation}}
    \argument{\lref{Eq1}; the fact that $F\in C^1([-\pi,\pi]^d,\R)$}{that there exists $c\in (0,\infty)$ such that for all $\ell\in (\{d\}\times\N\times\{1\})$ it holds that
    \begin{equation}\llabel{Eq2}
    \begin{split}
      v_\ell\leq \textstyle c\sum_{i=0}^{\fL_\ell}(i+1)^d(2^i)^{1+\frac d2}.
       \end{split}
    \end{equation}}
    \argument{\lref{def: N}}{that for all $\ell\in (\{d\}\times\N\times\{1\})$ it holds that
    \begin{equation}\llabel{Eq3}
    \begin{split}
        &\sum_{i=0}^{\fL_\ell}(i+1)^d(2^i)^{1+\frac d2}\leq (\fL_\ell+1)^d \sum_{i=0}^{\fL_\ell}(2^i)^{1+\frac d2}=(\fL_\ell+1)^d\frac{2^{1+\frac d2}}{2^{1+\frac d2}-1}\frac{\bigl(2^{1+\frac d2}\bigr)^{\fL_\ell+1}-1}{2^{1+\frac d2}}\\
        &\leq (\log_2(N_\ell)+2)^d\frac{\sqrt{2}}{\sqrt{2}-1}(2^{1+\frac d2})^{\fL_\ell}\leq \frac{\sqrt{2}}{\sqrt{2}-1}(\log_2(N_\ell)+2)^d(2N_\ell)^{1+\frac d2}.
        \end{split}
    \end{equation}}
\argument{\lref{Eq3};\lref{def: N};\lref{Eq2};}{that there exist $c_1,c_2\in (0,\infty)$ such that for all $\ell=(\ell_0,\ell_1,\ell_2)\in (\{d\}\times\N\times\{1\})$ it holds that
\begin{equation}\llabel{def: scrC}
\begin{split}
    v_\ell&\leq c_1(\log_2(N_\ell)+2)^d(2N_\ell)^{1+\frac d2}\\
    &\leq c_1\Bigl(\frac{d+2}{d(d+4)}\log_2(\ell_1)+2\Bigr)^d2^{1+\frac d2}(\ell_1)^{\frac{(d+2)^2}{2d(d+4)}}\\
    &\leq c_2 (\ell_1)^{\frac{1}{2d(d+4)}}(\ell_1)^{\frac{(d+2)^2}{2d(d+4)}}=c_2(\ell_1)^{\frac{1}{2}+\frac{5}{2d(d+4)}}.
    \end{split}
\end{equation}}
\argument{\lref{def: scrC};\lref{def: theta};\lref{eq1};the fact that $F\in C^1([-\pi,\pi]^d,\R)$}{that there exist $c_1,c_2,c_3\in (0,\infty)$ such that for all $\ell=(\ell_0,\ell_1,\ell_2)\in (\{d\}\times\N\times\{1\})$ it holds that
    \begin{align}
      & \textstyle\sup_{x\in (-1,1)^d}|u(x)-\reli{\ell}{2}{\theta_\ell}(x)|\textstyle\notag\\
      &\textstyle \leq c_1 (\ell_1)^{-\frac{d+2}{d(d+4)}} \max\bigl\{\sup_{x\in [-\pi,\pi]^d}|F(x)|,\sup_{x\in [-\pi,\pi]^d}\|(\nabla F)(x)\|\bigr\}+c_1\sqrt{\ln(\ell_1)}(\ell_1)^{-\frac{2d+3}{2d(d+4)}}\notag\\
       &\leq c_2 (\ell_1)^{-\frac{d+2}{d(d+4)}} +c_2\sqrt{\ln(\ell_1)}(\ell_1)^{-\frac{2d+3}{2d(d+4)}}\leq c_3(\ell_1)^{-\frac{1}{d+5}}\llabel{eqc}
    \end{align}
       \begin{equation}\llabel{def: fC2}
            \text{and}\qquad \sup_{i\in \{1,2,\dots,\fd(\ell)\}}|\theta_\ell^i|\leq \max\Bigl\{1,c_3(\ell_1)^{\frac{-1}{2}+\frac{5}{2d(d+4)}}\Bigr\}\leq \max\{1,c_3\}.
    \end{equation}}
    \argument{\lref{def: fC2}}{that there exist $c_1,c_2\in (0,\infty)$ such that for all $\ell=(\ell_0,\ell_1,\ell_2)\allowbreak\in (\{d\}\times\N\times\{1\}) $ there exists $\theta \in ( - c_2, c_2 )^{ \fd(\ell) }$ such that
    \begin{equation}\llabel{eq9}
\sup_{x\in (-1,1)^d}
  | u( x)
 -\textstyle
  \reli{\ell}{2}{\theta}( x )
| \leq
  c_1 (\ell_1)^{ \frac{-1}{d+5} } \ifnocf.
    \end{equation}}
    \argument{\lref{eq9};}{\lref{conclude}\dott}
\end{aproof}
\cfclear
\begin{athm}{lemma}{conj: Pinns for Poisson3 intro2 pre2 approximation shallow ann general}
    Let $d \in \N$, $a\in \R$, $b\in (a,\infty)$, $u \in C^2( [a,b]^d, \R )$. Then there exists $ c \in (0,\infty)$ such that for all $\ell=(\ell_0,\ell_1,\ell_2)\allowbreak\in (\{d\}\times\N\times\{1\}) $ there exists $\theta \in ( - c, c )^{ \fd(\ell) }$ such that
    \begin{equation}\llabel{conclude}
\sup_{x\in (a,b)^d}
  | u( x)
 -\textstyle
  \reli{\ell}{2}{\theta}( x )
| \leq
  c (\ell_1)^{ \frac{-1}{d+5} } \ifnocf.
    \end{equation}
    \cfout[.]
\end{athm}
\begin{aproof}
    Throughout this proof, let $v\colon [-1,1]^d\to\R$ satisfy for all $x\in [-1,1]^d$ that
    \begin{equation}\llabel{def: v}
        v(x)=\textstyle u\bigl(\frac{b-a}{2}x+(\frac{a+b}{2},\frac{a+b}{2},\dots,\frac{a+b}{2})\bigr).
    \end{equation}
    \argument{\cref{conj: Pinns for Poisson3 intro2 pre2 approximation shallow ann}}{that there exist $c\in (0,\infty)$ and $\theta_\ell=(\theta_\ell^1,\dots,\theta_\ell^{\fd(\ell)}) \in ( - c, c )^{ \fd(\ell) }$, $\ell\in \N^3 $, which satisfy for all $\ell=(\ell_0,\ell_1,\ell_2)\allowbreak\in (\{d\}\times\N\times\{1\}) $ that
    \begin{equation}\llabel{def: c}
        \sup_{x\in (-1,1)^d}
  | v( x)
 -\textstyle
  \reli{\ell}{2}{\theta_\ell}( x )
| \leq
  c (\ell_1)^{ \frac{-1}{d+5} } \ifnocf.
    \end{equation}}
    In the following for every $\ell=(\ell_0,\ell_1,\ell_2)\in \N^3$ let $\vartheta_\ell=(\vartheta_\ell^1,\dots,\vartheta_\ell^{\fd(\ell)})\in \R^{\fd(\ell)}$ satisfy that
    \begin{enumerate}[label=(\roman*)]
        \item \llabel{item 1} it holds for all $i\in \N\cap [1,\ell_1\ell_0]$ that $\vartheta_\ell^i=\frac{2}{b-a}\theta_\ell^i$,
        \item \llabel{item 2} it holds for all $i\in \N\cap(\ell_1\ell_0,\ell_1\ell_0+\ell_1] $ that $\vartheta_\ell^i=\theta_\ell^i-\frac{a+b}{2}$, and
        \item \llabel{item 3} it holds for all $i\in \N\cap(\ell_1\ell_0+\ell_1,\fd(\ell)] $ that $\vartheta_\ell^i=\theta_\ell^i$.
    \end{enumerate}
    \startnewargseq
    \argument{\lref{item 1};\lref{item 2};\lref{item 3}}{that for all $\ell\in (\{d\}\times\N\times\{1\})$ it holds that
    \begin{equation}\llabel{eq1}
         \reli{\ell}{2}{\theta_\ell}( x )= \reli{\ell}{2}{\vartheta_\ell}\bigl(\textstyle\frac{b-a}{2}x+(\frac{a+b}{2},\frac{a+b}{2},\dots,\frac{a+b}{2})\bigr).
    \end{equation}}
    \argument{\lref{eq1};\lref{def: v};\lref{def: c}}{for all $\ell\in (\{d\}\times\N\times\{1\})$ that
    \begin{equation}\llabel{eq2}
        \sup_{x\in (-1,1)^d}\textstyle\bigl|u\bigl(\frac{b-a}{2}x+(\frac{a+b}{2},\frac{a+b}{2},\dots,\frac{a+b}{2})\bigr)-\reli{\ell}{2}{\vartheta_\ell}\bigl(\textstyle\frac{b-a}{2}x+(\frac{a+b}{2},\frac{a+b}{2},\dots,\frac{a+b}{2})\bigr)\bigr|\leq c(\ell_1)^{ \frac{-1}{d+5} } . 
    \end{equation}}
    \argument{\lref{eq2}}{for all $\ell\in (\{d\}\times\N\times\{1\})$ that
    \begin{equation}\llabel{eq3}
        \sup_{x\in (a,b)^d}
  | u( x)
 -\textstyle
  \reli{\ell}{2}{\vartheta_\ell}( x )
| \leq
  c (\ell_1)^{ \frac{-1}{d+5} }.
    \end{equation}}
    \argument{\lref{def: c};\lref{item 1};\lref{item 2};\lref{item 3}}{for all $\ell\in (\{d\}\times\N\times\{1\})$, $i\in \{1,2,\dots,\fd(\ell)\}$ that
    \begin{equation}\llabel{eq4}
        \textstyle|\vartheta_\ell^i|\leq c(1+|\frac{b-a}{2}|)+|\frac{a+b}{2}|.
    \end{equation}}
    \argument{\lref{eq3};\lref{eq4}}{\lref{conclude}\dott}
\end{aproof}
\subsection{Error estimates for deep ANN approximations for 
$C^2$-functions}\label{subsec: error estimate C2}
\cref{main theorem} provides error estimates for the \DKM\ for the case of, both, shallow ($L=2$ in \cref{main theorem: conclude}) as well as deep ($L>2$ in \cref{main theorem: conclude}) \ANN\ approximations. In the situation of shallow \ANN\ approximations ($L=2$ in \cref{main theorem: conclude} in \cref{main theorem}) \cref{main theorem} can be deduced by employing the \ANN\ approximation result from \cref{conj: Pinns for Poisson3 intro2 pre2 approximation shallow ann general} above together with the arguments from Sections \ref{sec: priori estimate} and \ref{sec: deep Kolmogorov method}. For the situation of deep \ANN\ approximations ($L>2$ in \cref{main theorem: conclude} in \cref{main theorem}) we extend in this subsection \cref{conj: Pinns for Poisson3 intro2 pre2 approximation shallow ann general} from \cref{subsec: estimate shallow C2} above to deep \ANNs\ in \cref{prop: Pinns for Poisson3 intro2 pre2 approximation}, \cref{cor: Pinns for Poisson3 intro2 pre2 approximation pre}, and \cref{cor: Pinns for Poisson3 intro2 pre2 approximation} below. Our proofs of \cref{prop: Pinns for Poisson3 intro2 pre2 approximation}, \cref{cor: Pinns for Poisson3 intro2 pre2 approximation pre}, and \cref{cor: Pinns for Poisson3 intro2 pre2 approximation}, respectively, make use of the elementary \ANN\ extension result in \cref{lem: embedding ANN} below.
\begin{athm}{lemma}{lem: existence of smooth}
    Let $d \in \N$, let $O \subseteq \R^d$ be open, and let $K \subseteq O$ be compact. Then there exists $\phi \in C^{ \infty }( \R^d, \R )$ such that
\begin{enumerate}[label=(\roman*)]
    \item \label{item 1: existence of smooth function} it holds for all $x\in K$ that $\phi(x)=1$ and
    \item \label{item 2: existence of smooth function}  it holds for all $x\in \R^d\backslash O$ that $\phi(x)=0$.
\end{enumerate}
\end{athm}
\begin{aproof}
Throughout this proof, assume without loss of generality that $K \neq \emptyset$ (otherwise we can choose $\phi = ( \R^d \ni x \mapsto 0 \in \R)$), for every $r\in \R$, $x\in \R^d$ let $B_r(x)\subseteq\R^d$ satisfy $B_r(x)=\{y\in \R^d\colon \|x-y\|< r\}$, let $R \in (0,\infty)$ satisfy $K \subseteq B_R( 0 )$, 
let $\bbO \subseteq \R^d$ satisfy $\bbO = O \cap B_R( 0 )$, let $r\in \R$ satisfy
\begin{equation}\llabel{def: r}
    r=\tfrac {1}{3}\textstyle\bigl[\inf_{ x \in K } \inf_{ y \in (\R^d \backslash \bbO) } \| x - y \|\bigr],
\end{equation}
for every $x \in \R^d$ let
$f_x \colon \R^d \to \R$ 
satisfy for all $y \in \R^d$ that
\begin{equation}\llabel{def: fv}
    f_x(y)=\begin{cases}
        \exp( [ \| x - y \|^2 - r^2 ]^{ - 1 } )&\colon \| x - y \| < r,\\
        0 &\colon \| x - y \| \geq r,
    \end{cases}
\end{equation}
let $N \in \N \cap ( r^{ - 1 } 2 d^{ 1 / 2 }, \infty )$, let $I\subseteq\R^d$ satisfy 
 \begin{equation}\llabel{def: I}
     I = \{ v \in \R^d \colon N v \in \Z^d \}, 
 \end{equation} 
and let $\phi\colon \R^d\to\R$ satisfy for all $x\in \R^d$ that
\begin{equation}\llabel{def: phi}
   \textstyle \phi(x)=\bigl( \sum_{ v \in I,\, \operatorname{supp}( f_v ) \cap K \neq \emptyset } f_v(x) \bigr) \bigl( \sum_{ v\in I } f_v( x ) \bigr)^{-1}
\end{equation}
\argument{\lref{def: fv};the fact that $\bigl([0,\infty)\ni x\mapsto \exp\bigl(\frac{1}{x-r^2+\mathbbm 1_{[0,r^2]}(x)}\bigr)\mathbbm 1_{[0,r^2)}(x)\in \R\bigr)\in C^{\infty}([0,\infty),\R)$; the fact that for all $v\in \R^\d$ it holds that $(\R^d\ni w\mapsto \|v-w\|^2\in \R)\in C^{\infty}(\R^d,\R)$}{that for all $v\in \R^\d$ it holds that 
\begin{equation}\llabel{arg1}
    f_v\in C^{\infty}(\R^d,\R).
\end{equation}}
\argument{\lref{def: fv};\lref{def: I}}{for all $s\in (0,\infty)$, $x\in B_s(0)$ that
\begin{equation}\llabel{eq1}
    \sum_{v\in I}f_v(x)=\sum_{v\in I\cap B_{r+s}(0)}f_v(x)
\end{equation}
\begin{equation}\llabel{eq1'}
   \text{and}\qquad \sum_{ v \in I,\, \operatorname{supp}( f_v ) \cap K \neq \emptyset } f_v(x)=\sum_{ v \in I\cap B_{r+s}(0),  \,\operatorname{supp}( f_v ) \cap K \neq \emptyset } f_v(x).
\end{equation}}
\argument{\lref{def: I};}{for all $s\in (0,\infty)$ that 
\begin{equation}\llabel{eq2}
    \#\{I\cap B_s(0)\}<\infty.
\end{equation}}
\argument{\lref{eq1};\lref{eq2}}{for all $s\in (0,\infty)$ that
\begin{equation}\llabel{eq3}
  \bigl( \textstyle B_s(0)\ni x \mapsto \sum_{v\in I}f_v(x)\in \R\bigr)\in C^\infty(B_s(0),\R)
\end{equation}
\begin{equation}\llabel{eq3'}
    \text{and}\qquad \bigl( \textstyle B_s(0)\ni x \mapsto \sum_{ v \in I,\, \operatorname{supp}( f_v ) \cap K \neq \emptyset }f_v(x)\in \R\bigr)\in C^\infty(B_s(0),\R).
\end{equation}
}
\argument{\lref{eq3'}}{that
\begin{equation}\llabel{eq4}
  \bigl( \textstyle \R^d\ni x \mapsto \sum_{v\in I}f_v(x)\in \R\bigr)\in C^\infty(\R^d,\R)
\end{equation}
\begin{equation}\llabel{eq4'}
   \text{and}\qquad  \bigl( \textstyle \R^d\ni x \mapsto \sum_{ v \in I,\, \operatorname{supp}( f_v ) \cap K \neq \emptyset }f_v(x)\in \R\bigr)\in C^\infty(\R^d,\R).
\end{equation}
}
\argument{the fact that $N\geq \frac{2\sqrt{d}}{r}$}{that for all $x=(x_1,\dots,x_d)\in \R^d$ there exists $y=(y_1,y_2,\dots,y_d)\in \R^d$ such that for all $i\in \{1,2,\dots,d\}$ it holds that 
\begin{equation}\llabel{arg2}
    |x_i-y_i|< \frac{r}{\sqrt{d}} \qqandqq Ny_i\in \Z.
\end{equation}}
\argument{\lref{arg2};\lref{def: I}}{for all \llabel{arg2'} $x\in \R^d$ that there exists $y\in I$ such that $\|x-y\|<r$\dott}
\argument{\lref{arg2'};\lref{def: fv}}{for all $x\in \R^d$ it holds that
\begin{equation}\llabel{eq5}
    \sum_{v\in I}f_v(x)>0.
\end{equation}}
\argument{\lref{def: phi};\lref{eq5};\lref{eq4};\lref{eq4'}}{that
\begin{equation}\llabel{evidence 1}
    \phi\in C^\infty(\R^d,\R).
\end{equation}}
\argument{\lref{def: fv};}{that for all $v\in \R^d$ it holds that
\begin{equation}\llabel{eq6}
    \operatorname{supp}(f_v)=\overline{B_r(v)}.
\end{equation}}
\argument{\lref{eq6}}{for all $x\in K$, $v\in I\cap \overline{B_r(x)}$ that
\begin{equation}\llabel{eq7}
   \operatorname{supp}(f_v)\cap K\neq \emptyset.
\end{equation}}
\argument{\lref{eq7}; \lref{def: phi};the fact that for all $x\in \R^d$ it holds that \begin{equation}
    \sum_{v\in I}f_v(x)=\sum_{v\in I\cap \overline{B_r(x)}}f_v(x)
\end{equation}}{that for all $x\in K$ it holds that
\begin{equation}\llabel{evidence 2}
    \phi(x)=\textstyle\bigl( \sum_{ v \in I\cap \overline{B_r(x)}} f_v(x) \bigr) \bigl( \sum_{ v\in I \cap \overline{B_r(x)}} f_v( x ) \bigr)^{-1}=1.
\end{equation}}
\argument{\lref{def: r};\lref{eq6}}{for all $x\in \R^d\backslash \bbO$, $v\in \overline{B_r(x)}$ it holds that
\begin{equation}\llabel{eq8}
    \operatorname{supp}(f_v)\cap K=\overline{B_r(v)}\cap K\subseteq \overline{B_{2r}(x)}\cap K=\emptyset.
\end{equation}}
\argument{\lref{eq8};the fact that for all $x\in \R^d$ it holds that \begin{equation}
    \sum_{v\in I,\,  \operatorname{supp}(f_v)\cap K\neq \emptyset}f_v(x)=\sum_{v\in I\cap \overline{B_r(x)},\, \operatorname{supp}(f_v)\cap K\neq \emptyset}f_v(x)
\end{equation}}{that for all $x\in \R^d\backslash \bbO$ it holds that
\begin{equation}\llabel{evidence 3}
    \phi(x)=\textstyle\bigl( \sum_{ v \in I\cap \overline{B_r(x)},\,\operatorname{supp}(f_v)\cap K\neq \emptyset} f_v(x) \bigr) \bigl( \sum_{ v\in I} f_v( x ) \bigr)^{-1}=0.
\end{equation}}
\argument{\lref{evidence 1};\lref{evidence 2};\lref{evidence 3}}{\cref{item 1: existence of smooth function,item 2: existence of smooth function} \dott}
\end{aproof}
\renewcommand{\scrl}{\mathbf{l}}
\begin{athm}{lemma}{lem: embedding ANN}
    Let $d\in\N$, let $D\subseteq\R^d$ be bounded, let $L\in \N\backslash\{1\}$, $\ell=(\ell_0,\ell_1,\dots,\ell_L)\allowbreak\in (\{d\}\times\N^{L-1}\times\{1\}) $, $\scrl\in (\{d\}\times\N\times\{1\})$, $\vartheta\in \R^{\fd(\scrl)}$, assume $\ell_1\geq \bfl_1$, and assume for all $i\in \N\cap(1,L)$ that $\ell_i\geq 2$. Then there exists $\theta\in [-\|\vartheta\|-1,\|\vartheta\|+1]^{\fd(\ell)}$ which satisfies for all $x\in D$ that
    \begin{equation}\llabel{conclude}
        \reli{\ell}{L}{\theta}(x)=\reli{\scrl}{2}{\vartheta}(x).
    \end{equation}
\end{athm}
\begin{aproof}
    Throughout this proof, let $\xi=(\xi_1,\dots,\xi_{\fd(\scrl)})\in \R^{\fd(\scrl)}$ satisfy $\xi=\vartheta$.
    \argument{the assumption that $\ell_1\geq \bfl_1$;the assumption that for all $i\in \N\cap(1,L)$ it holds that $\ell_i\geq 2$;\unskip, \eg, \cite[Lemma 4.4.7]{ArBePhi2024} (applied with $L\curvearrowleft 2$, $\fL\curvearrowleft L$, $l\curvearrowleft \scrl$, $\fl\curvearrowleft \ell$ in the notation of \cite[Lemma 4.4.7]{ArBePhi2024}) }{that there exists $\theta=(\theta_1,\dots,\theta_{\fd(\ell)})\in \R^{\fd(\ell)}$ which satisfies that
    \begin{enumerate}[label=(\roman*)]
        \item \llabel{item 1} it holds for all $x\in \R^{\fd}$ that
    $\reli{\ell}{L}{\theta}(x)=\reli{\scrl}{2}{\xi}(x)$ and
    \item \llabel{item 2} it holds that $\max\{|\theta_1|,|\theta_2|,\dots,|\theta_{\fd(\ell)}|\}\leq \max\{1,\max\{|\xi_1|,|\xi_2|,\dots,|\xi_{\fd(\scrl)}|\}\}$.
    \end{enumerate}}
    \startnewargseq
    \argument{\lref{item 1};\lref{item 2};the fact that $\max\{1,\max\{|\xi_1|,|\xi_2|,\allowbreak\dots,\allowbreak |\xi_{\fd(\scrl)}|\}\}\leq \|\xi\|+1$ }{\lref{conclude}\dott}
\end{aproof}
\cfclear
\renewcommand{\scrl}{\mathscr{l}}
\begin{athm}{lemma}{prop: Pinns for Poisson3 intro2 pre2 approximation}
    Let $d \in \N$, $a\in \R$, $b\in (a,\infty)$, $u \in C^2( [a,b]^d, \R )$. Then there exists $c \in (0,\infty)$ such that for all $L\in \N\backslash\{1\}$, $\ell=(\ell_0,\ell_1,\dots,\ell_L)\allowbreak\in (\{d\}\times\N^{L-1}\times\{1\}) $ there exists $\theta \in ( - c,c )^{ \fd(\ell) }$ such that
    \begin{equation}\llabel{conclude}
\displaystyle\int_{(a,b)^d}|
   u( x)
 -\textstyle
  \reli{ \ell}{L}{\theta } ( x )
|^2\,\d x\leq 
  c [ \min\{ \ell_1, \ell_2, \dots, \ell_{ L - 1 } \} ]^{ 
\frac{-2}{d+5}}
 \ifnocf.
    \end{equation}
    \cfout[.]
\end{athm}
\begin{aproof}
Throughout this proof, let $\bfl\in \N^3$, $\vartheta\in \R^{\fd(\bfl)}$ satisfy $\bfl=(1,2,1)$ and $\vartheta=(1,-1,0,0,1,-1,0)$.
    \argument{\unskip, \eg, \cite[Lemma 2.2.7]{ArBePhi2024}}{that for all $x\in \R$ it holds that 
\begin{equation}\llabel{def: bfl}
    \reli{\bfl}{2}{\vartheta}(x)=x.
\end{equation}}
\startnewargseq
    \argument{\cref{conj: Pinns for Poisson3 intro2 pre2 approximation shallow ann};}{that there exist $ c \in (0,\infty)$ and $\eta \colon  (\{d\}\times (\cup_{L=0}^\infty\N^L)\times \{1\})\to (-c,c)^{ \fd(\ell) }$ which satisfy for all $\ell=(\ell_0,\ell_1,\ell_2)\allowbreak\in (\{d\}\times\N\times\{1\}) $ that 
    \begin{equation}\llabel{def: rc}
\displaystyle\int_{(a,b)^d}|
   u( x)
 -\textstyle
  \reli{ \ell}{2}{\eta(\ell) } ( x )
|^2\,\d x\leq 
  c (  \ell_1)^{ 
\frac{-2}{d+5}}.
    \end{equation}}
    \startnewargseq
    \argument{\cref{realization multi};the fact that $u$ is bounded;}{that there exists $\fC\in\R$ such that for all  $L\in \N\backslash\{1\}$, $\ell=(\ell_0,\ell_1,\dots,\ell_L)\allowbreak\in (\{d\}\times(\N^{L-1}\times\{1\}) $ with $\min\{\ell_1,\ell_2,\dots,\ell_{L-1}\}\allowbreak= 1$ it holds that 
    \begin{equation}\llabel{evidence 1}
    \begin{split}
\displaystyle\int_{(a,b)^d}|
   u( x)
 -\textstyle
  \reli{ \ell}{L}{0 } ( x )
|^2\,\d x=\displaystyle\int_{(a,b)^d}|
   u( x)
|^2\,\d x\leq \fC = 
  \fC [\min\{  \ell_1, \ell_2, \dots, \ell_{ L - 1 } \} ]^{  
\frac{-2}{d+5} }.
  \end{split}
    \end{equation}}
  \argument{\cref{lem: embedding ANN};\lref{def: bfl}}{that for all $\mathscr{l}=(\scrl_0,\scrl_1,\scrl_2)\in (\{d\}\times\N\times\{1\})$, $L\in \N\backslash\{1\}$, $\ell=(\ell_0,\ell_1,\dots,\ell_L)\allowbreak\in (\{d\}\times(\N\cap[\max\{2,\scrl_1\},\infty))^{L-1}\times\{1\}) $ there exists $\theta \in [ -\|\vartheta\|-1, \|\vartheta\|+1 ]^{ \fd(\ell) }$ such that for all $x\in \R^d$ it holds that
    \begin{equation}\llabel{eq1}
\reli{\ell}{L}{\theta}(x)=\reli{\scrl}{2}{\eta(\scrl)}(x).
    \end{equation}}
      \argument{\lref{eq1};}{that for all $\mathscr{l}=(\scrl_0,\scrl_1,\scrl_2)\in (\{d\}\times(\N\cap[2,\infty))\times\{1\})$, $L\in \N\backslash\{1\}$, $\ell=(\ell_0,\ell_1,\dots,\ell_L)\allowbreak\in (\{d\}\times(\N\cap[\scrl_1,\infty))^{L-1}\times\{1\}) $ there exists $\theta \in ( -\max\{\|\vartheta\|+1, c\}, \max\{\|\vartheta\|+1, c\})^{ \fd(\ell) }$ such that for all $x\in \R^d$ it holds that
    \begin{equation}\llabel{eq2}
\reli{\ell}{L}{\theta}(x)=\reli{\scrl}{2}{\eta(\scrl)}(x).
    \end{equation}}
    \argument{\lref{eq2};\lref{def: rc}}{that for all $L\in \N\backslash\{1\}$, $\ell=(\ell_0,\ell_1,\dots\ell_L)\allowbreak\in (\{d\}\times\N^{L-1}\times\{1\}) $ with $\min\{\ell_1,\ell_2,\dots,\ell_{L-1}\}\allowbreak\geq2$ there exists $\theta\in ( -\max\{\|\vartheta\|+1, c\}, \max\{\|\vartheta\|+1, c\})^{ \fd(\ell) }$ such that 
    \begin{equation}\llabel{evidence 3}
\displaystyle\int_{(a,b)^d}|
   u( x)
 -\textstyle
  \reli{ \ell}{L}{\theta } ( x )
|^2\,\d x\leq 
  c [ \min\{  \ell_1, \ell_2, \dots, \ell_{ L - 1 } \} ]^{  
\frac{-2}{d+5} } .
    \end{equation}}
    \argument{\lref{evidence 3};\lref{evidence 1}}{\lref{conclude}\dott}
\end{aproof}
\begin{athm}{cor}{cor: Pinns for Poisson3 intro2 pre2 approximation pre}
    Let $d \in \N$, let $D \subseteq \R^d$ be measurable and bounded, and let $u \in C^2( \bar{D}, \R )$. Then there exists $ c \in (0,\infty)$ such that for all $L\in \N\backslash\{1\}$, $\ell=(\ell_0,\ell_1,\dots,\ell_L)\allowbreak\in (\{d\}\times\N^{L-1}\times\{1\}) $ there exists $\theta \in ( - c,c )^{ \fd(\ell) }$ such that
    \begin{equation}\llabel{conclude}
\displaystyle\int_D|
   u( x)
 -\textstyle
  \reli{ \ell}{L}{\theta } ( x )
|^2\,\d x\leq 
   c [ \min\{ \ell_1, \ell_2, \dots, \ell_{ L - 1 } \} ]^{  
\frac{-2}{d+5} } \ifnocf.
    \end{equation}
    \cfout[.]
\end{athm}
\begin{aproof}
    \argument{the assumption that $D$ is bounded;}{that there exist $a\in \R$, $b\in (a,\infty)$ which satisfy
    \begin{equation}\llabel{def: ab}
        \bar{D}\subseteq (a,b)^d.
    \end{equation}}
    \startnewargseq
    \argument{the assumption that $u\in C^2(\bar{D},\R)$}{that there exists an open set $V\subseteq\R^d$ and $\bar{u}\in C^2(V,\R)$ which satisfy for all $x\in \bar{D}$ that
    \begin{equation}\llabel{def: V}
        \bar{D}\subseteq V \qqandqq  u(x)=\bar{u}(x).
    \end{equation}}
    \startnewargseq
    \argument{\lref{def: V};\cref{lem: existence of smooth}; the assumption that $D$ is bounded}{that there exists $\phi\in C^{\infty}(\R^d,\R)$ which satisfies for all $x\in \bar{D}$, $y\in \R^d\backslash V$ that
\begin{equation}\llabel{def: phi}
    \phi(x)=1\qqandqq \phi(y)=0.
\end{equation}}
\startnewargseq
In the following let $v\colon \R^d\to\R$ satisfy for all $x\in \R^d$ that
\begin{equation}\llabel{def: v}
    v(x)=\begin{cases}
        \phi(x)\bar{u}(x) &\colon x\in V\\
        0 &\colon x\in \R^d\backslash V.
    \end{cases}
\end{equation}
\startnewargseq
\argument{\lref{def: V};\lref{def: phi};\lref{def: v}}{that for all $x\in \bar{D}$ it holds that 
\begin{equation}\llabel{def: v1}
   v(x)=\phi(x)\bar{u}(x)=\phi(x)u(x)=u(x).
\end{equation}}
\argument{\lref{def: phi};\lref{def: v}; the fact that $\bar{u}\in C^2(V,\R)$}{that \llabel{arg1} $v\in C^1(\R^d,\R)$\dott}
    \argument{\lref{def: ab};\lref{def: v1};\lref{arg1};\cref{prop: Pinns for Poisson3 intro2 pre2 approximation}}{that there exists $c\in (0,\infty)$ such that for all $L\in \N\backslash\{1\}$, $\ell=(\ell_0,\ell_1,\dots,\ell_L)\allowbreak\in (\{d\}\times\N^{L-1}\times\{1\}) $ there exists $\theta \in ( -c,c )^{ \fd(\ell) }$ such that
    \begin{equation}\llabel{eq1}
    \begin{split}
        \displaystyle\int_D|
   u( x)
 -\textstyle
  \reli{ \ell}{L}{\theta } ( x )
|^2\,\d x&\leq \displaystyle\int_{(a,b)^d}|
   v( x)
 -\textstyle
  \reli{ \ell}{L}{\theta } ( x )
|^2\,\d x\\
&\leq 
  c [ \min\{ \ell_1, \ell_2, \dots, \ell_{ L - 1 } \} ]^{  \frac{-2}{d+5} } 
  \end{split}
    \end{equation}}
    \argument{\lref{eq1};}{\lref{conclude}\dott}
\end{aproof}
\begin{athm}{cor}{cor: Pinns for Poisson3 intro2 pre2 approximation}
    Let $d \in \N$, let $D \subseteq \R^d$ be measurable and bounded, and let $u \in C^2( \bar{D}, \R )$. Then there exists $ c \in (0,\infty)$ such that for all $R\in (0,\infty]$, $L\in \N\backslash\{1\}$, $\ell=(\ell_0,\ell_1,\dots,\ell_L)\allowbreak\in (\{d\}\times\N^{L-1}\times\{1\}) $ there exists $\theta \in ( - \min\{c,R\}, \min\{c,R\} )^{ \fd(\ell) }$ such that
    \begin{equation}\llabel{conclude}
\displaystyle\int_D|
   u( x)
 -\textstyle
  \reli{ \ell}{L}{\theta } ( x )
|^2\,\d x\leq 
   c\mathbbm 1_{(-\infty,c)}(R)+ c [ \min\{ \ell_1, \ell_2, \dots, \ell_{ L - 1 } \} ]^{  
\frac{-2}{d+5} } \ifnocf.
    \end{equation}
    \cfout[.]
\end{athm}
\begin{aproof}
  \argument{\cref{cor: Pinns for Poisson3 intro2 pre2 approximation pre};}{that there exists $ c \in (0,\infty)$ such that for all $L\in \N\backslash\{1\}$, $\ell=(\ell_0,\ell_1,\dots,\ell_L)\allowbreak\in (\{d\}\times\N^{L-1}\times\{1\}) $ there exists $\theta \in ( - c,c )^{ \fd(\ell) }$ such that
    \begin{equation}\llabel{def: c}
\displaystyle\int_D|
   u( x)
 -\textstyle
  \reli{ \ell}{L}{\theta } ( x )
|^2\,\d x\leq 
   c [ \min\{ \ell_1, \ell_2, \dots, \ell_{ L - 1 } \} ]^{ 
\frac{-2}{d+5} } \ifnocf.
    \end{equation}}
    \startnewargseq
    \argument{\lref{def: c};}{that for all $R\in [c,\infty]$, $L\in \N\backslash\{1\}$, $\ell=(\ell_0,\ell_1,\dots,\ell_L)\allowbreak\in (\{d\}\times\N^{L-1}\times\{1\}) $ there exists $\theta \in ( - \min\{c,R\}, \min\{c,R\} )^{ \fd(\ell) }$ such that
    \begin{equation}\llabel{evidence 1}
    \begin{split}
\displaystyle\int_D|
   u( x)
 -\textstyle
  \reli{ \ell}{L}{\theta } ( x )
|^2\,\d x&\leq c [ \min\{ \ell_1, \ell_2, \dots, \ell_{ L - 1 } \} ]^{  
\frac{-2}{d+5} } .
\end{split}
    \end{equation}}
    \argument{\cref{realization multi};the fact that $u$ is bounded;the assumption that $f$ is non-decreasing;}{that there exists $\fC\in (0,\infty)$ such that for all $R\in (0,c)$, $L\in \N\backslash\{1\}$, $\ell=(\ell_0,\ell_1,\dots,\ell_L)\allowbreak\in (\{d\}\times\N^{L-1}\times\{1\}) $ there exists $\theta \in ( - \min\{c,R\},\allowbreak \min\{c,R\} )^{ \fd(\ell) }$ such that
    \begin{equation}\llabel{evidence 2}
    \begin{split}
&\displaystyle\int_D|
   u( x)
 -\textstyle
  \reli{ \ell}{L}{0 } ( x )
|^2\,\d x=\displaystyle\int_D|
   u( x)
|^2\,\d x\leq \fC= \fC\mathbbm 1_{(-\infty,c)}(R).
\end{split}
    \end{equation}}
    \argument{\lref{evidence 1};\lref{evidence 2}}{\lref{conclude}\dott}
\end{aproof}
\renewcommand{\Reli}[4]{\mathbf{N}^{#1,#3,#4}_{#2}}
\renewcommand{\Relii}[5]{\mathbf{N}^{#1,#3,#4}_{#2,#5}}
\section{Error estimates for the deep Kolmogorov method}\label{sec: deep Kolmogorov method}
In this section we employ \cref{lem: anngrowth} and \cref{cor: ann lipschitz} from \cref{sec: priori estimate} and \cref{cor: Pinns for Poisson3 intro2 pre2 approximation} from \cref{sec: ANN approximation} to establish in \cref{prop: linear II relu pinn draft}, \cref{prop: linear II relu pinn}, \cref{cor: linear II relu pinn}, and \cref{main cor} below error estimates for the \DKM. \cref{main theorem} in the introduction is an immediate consequence of \cref{main cor}. Among other things, our proof of \cref{prop: linear II relu pinn draft} employs
\begin{itemize}
    \item the well-known Sobolev estimates with explicit embedding constants in \cref{lem: montecarlo 4} and \cref{cor: montecarlo 4} below,
    \item the Feynman-Kac representation for solutions of heat \PDEs\ that we recall in \cref{Feynman-Kac formula} below, and 
    \item the well-known error estimates for the Monte Carlo method in \cref{lem: montecarlo2} below.
\end{itemize}
 We also refer, \eg, to \cite[Lemma 4]{MR4134774} and \cite[Appendix C]{Hiebersupplement} for abstract general error estimates for \ANN\ approximations in terms of the optimization error. 
\subsection{Sobolev embedding estimates with explicit constants}
\begin{athm}{lemma}{lem: montecarlo 4}
  Let $d\in \N$, $\Phi\in C^1((0,1)^d,\R)$. Then
  \begin{equation}\llabel{conclude}
      \sup_{x\in (0,1)^d}|\Phi(x)|\leq 8\sqrt{e}\biggl[\int_{(0,1)^d}\bigl(|\Phi(x)|^{\max\{2,d^2\}}+\|(\nabla\Phi)(x)\|^{\max\{2,d^2\}}\bigr)\,\d x\biggr ]^{\nicefrac{1}{\max\{2,d^2\}}}.
\end{equation}
\end{athm}
\begin{aproof}
    \argument{\unskip, \eg, \cite[Corollary 2.15]{MR4454915}}{ \lref{conclude}\dott}
\end{aproof}
\begin{athm}{cor}{cor: montecarlo 4}
    Let $d\in \N$, $a\in \R$, $b\in (a,\infty)$, $p\in [\max\{2,d^2\},\infty)$, $\Phi\in C^1((a,b)^d,\R)$. Then
  \begin{equation}\llabel{conclude}
  \begin{split}
      \sup_{x\in (a,b)^d}|\Phi(x)|\leq 16\sqrt{e}\max\{(b-a)^{\nicefrac{-d}{p}},b-a\}\biggl[\int_{(a,b)^d}\bigl(|\Phi(x)|^{p}+\|(\nabla\Phi)(x)\|^{p}\bigr)\,\d x\biggr ]^{\nicefrac{1}{p}}.
         \end{split}
      \end{equation}
\end{athm}
\begin{aproof}
Throughout this proof, let $\bfe\in \R^d$ satisfy $\bfe=(1,1,\dots,1)$ and let $\Psi\colon (0,1)^d\to\R$ satisfy for all $x=(x_1,\dots,x_d)\in (0,1)^d$ that
\begin{equation}\llabel{def: Psi}
\Psi(x)=\Phi(a\bfe+(b-a)x).
\end{equation}
\startnewargseq
\argument{\lref{def: Psi}; the assumption that $\Phi\in C^1((a,b)^d,\R)$;}{that \llabel{arg1} $\Psi\in C^1((0,1)^d,\R)$\dott}
\argument{\lref{arg1};\cref{lem: montecarlo 4};the assumption that $p\geq \max\{2,d^2\}$; Holder inequality}{that
\begin{equation}\llabel{eq1}
\begin{split}
     \sup_{x\in (0,1)^d}|\Psi(x)|
     &\leq 8\sqrt{e}\biggl[\int_{(0,1)^d}\bigl(|\Psi(x)|^{\max\{2,d^2\}}+\|(\nabla\Psi)(x)\|^{\max\{2,d^2\}}\bigr)\,\d x\biggr ]^{\nicefrac{1}{\max\{2,d^2\}}}\\
     &\leq 8\sqrt{e}\biggl[\int_{(0,1)^d}\bigl(|\Psi(x)|^{\max\{2,d^2\}}+\|(\nabla\Psi)(x)\|^{\max\{2,d^2\}}\bigr)^{\nicefrac{p}{\max\{2,d^2\}}}\,\d x\biggr ]^{\nicefrac{1}{p}}.
     \end{split}
\end{equation}}
\argument{\lref{eq1};the fact that for all $x,y\in \R^d$, $p\in [1,\infty)$ it holds that $\|x+y\|^p\leq 2^{p}\|x\|^p+2^{p}\|y\|^p$; the fact that $\max\{2,d^2\}\geq 1$}{that 
\begin{equation}\llabel{eq1'}
\begin{split}
     \sup_{x\in (0,1)^d}|\Psi(x)|
    \leq 16 \sqrt{e}\biggl[\int_{(0,1)^d}\bigl(|\Psi(x)|^{p}+\|(\nabla\Psi)(x)\|^{p}\bigr)\,\d x\biggr ]^{\nicefrac{1}{p}}.
     \end{split}
\end{equation}}
\argument{\lref{def: Psi}}{that for all $x\in (0,1)^d$ it holds that
\begin{equation}\llabel{eq2}
    (\nabla \Psi)(x)= (b-a)(\nabla\Phi)(a\bfe+(b-a)x)
\end{equation}}
\argument{\lref{eq2};\lref{def: Psi}}{that
\begin{align}
   & \int_{(0,1)^d}\bigl(|\Psi(x)|^{p}+\|(\nabla\Psi)(x)\|^{p}\bigr)\,\d x \notag\\
   &\leq \max\{1,(b-a)^{p}\}\biggl[\int_{(0,1)^d}\bigl(|\Psi(a\bfe+(b-a)x)|^{p}+\|(\nabla\Psi)(a\bfe+(b-a)x)\|^{p}\bigr)\,\d x\biggr]\notag\\
   &=\frac{\max\{1,(b-a)^{p}\}}{(b-a)^d}\biggl[\int_{(0,1)^d}\bigl(|\Psi(a\bfe+(b-a)x)|^{p}\notag\\
   &\quad\quad+\|(\nabla\Psi)(a\bfe+(b-a)x)\|^{p}\bigr)\,\d (a\bfe+(b-a)x)\biggr]\notag\\
   &=\frac{\max\{1,(b-a)^{p}\}}{(b-a)^d}\biggl[\int_{(a,b)^d}\bigl(|\Phi(x)|^{p}+\|(\nabla\Phi)(x)\|^{p}\bigr)\,\d x\biggr ]\llabel{eq3}\\
   &= \max\{(b-a)^{-d}, (b-a)^{p}\}\biggl[\int_{(a,b)^d}\bigl(|\Phi(x)|^{p}+\|(\nabla\Phi)(x)\|^{p}\bigr)\,\d x\biggr ]\notag.
\end{align}}
\argument{\lref{eq3};\lref{def: Psi};\lref{eq1'}}{\lref{conclude}\dott}
\end{aproof}
\subsection{Feynman-Kac representations for solutions of heat PDEs}
\begin{athm}{lemma}{Feynman-Kac formula}
Let $ T\in (0,\infty) $, $ d \in \N $, let $u\in C([0,T]\times\R^{d},\R)$ be at most polynomially growing, assume $u|_{ (0,T) \times \R^d } \in C^{ 1, 2 }( (0,T)\times\R^d, \R )$, assume for all $t\in (0,T)$, $x\in \R^d$ that
\begin{equation}\llabel{pde}
\textstyle
  \frac{ \partial }{ \partial t } u( t, x )
  +
 \frac{1}{2}
  \Delta_x u(t,x)
  = 0,
\end{equation}
let $ ( \Omega, \mathcal{F}, \P ) $ be a probability space,
let $ W \colon [0,T] \times \Omega \to \R^d $ be a standard Brownian motion, let $\bbT\colon\Omega\to[0,T]$ and $\bbX\colon\Omega\to\R^d$ be random variables, let $\polycons\in [0,\infty)$ satisfy $\E[\|\bbX\|^{\polycons}]+\sup_{x\in \R^d}((1+\|x\|)^{-\polycons}|u(T,x)|)<\infty$, and assume that $ (\bbT, \mathbb{X}) $ and $ W$ are independent. Then it holds $\P$-a.s.\ that
\begin{equation}\llabel{conclude}
    \E\bigl[u(T,\bbX+W_{T-\bbT})\bigl|(\bbT,\bbX)\bigr]=u(\bbT,\bbX).
\end{equation}
\end{athm}
\begin{aproof}
    Throughout this proof, let $g\colon \R^d\to\R$ satisfy for all $x\in \R^d$ that $g(x)=u(T,x)$, let $v\colon [0,T]\times \R^d\to\R$ satisfy for all $t\in [0,T]$, $x\in \R^d$ that
    \begin{equation}\llabel{def: v}
        v(t,x)=u(T-t,x).
    \end{equation}
    \argument{\lref{pde};\lref{def: v}}{that for all $t\in (0,T)$, $x\in \R^d$ it holds that
    \begin{equation}\llabel{eq1}
        \textstyle
  \frac{ \partial }{ \partial t } v( t, x )
  -
 \frac{1}{2}
  \Delta_x v(t,x)
  = 0
\qquad
  \text{and}
\qquad
  v(0,x) = g(x),
\end{equation}
}
\argument{\lref{eq1};\cite[Corollary 4.17]{HairerHutzenthalerJentzen2015}}{for all $t\in [0,T]$, $x\in \R^d$ that
\begin{equation}\llabel{eq2}
    v(t,x)=\E[g(x+W_t)].
\end{equation}}
\argument{\lref{eq2};\lref{def: v}}{for all $t\in [0,T]$, $x\in \R^d$ that
\begin{equation}\llabel{eq3}
    u(t,x)=\E[g(x+W_{T-t})].
\end{equation}}
In the following let $D=[0,T]\times\R^d$, $E=C([0,T],\R^d)$, let $\Phi\colon D\times E\to [-\infty,\infty]$ satisfy for all $(t,x)\in [0,T]\times \R^d$, $\bfW\in C([0,T]\times\R^d)$ that
\begin{equation}\llabel{def: Phi}
    \Phi(t,x,\bfW)=g(x,\bfW_{T-t}),
\end{equation}
and let $\phi\colon D\to[-\infty,\infty]$ satisfy for all $(t,x)\in D$ that
\begin{equation}\llabel{def: phi}
    \phi(t,x)=\E[\Phi(t,x,W)].
\end{equation}
\argument{\lref{def: Phi};the fact that $\E[\|\bbX\|^\polycons]<\infty$; the fact that $\sup_{x\in \R^d} ((1\allowbreak+\|x\|)^{-\polycons}|g(x)|)<\infty$; the fact that $\P(|\bbT|\leq T)=1$; the assumption that $W$ is a standard Brownian motion;}{that there exists $c\in \R$ such that
\begin{equation}\llabel{evidence 1}
    \textstyle\E\bigl[|\Phi(\bbT,\bbX,W)|\bigr]=\E\bigl[|g(\bbX,W_{T-\bbT})|\bigr]\leq c\bigl(1+\E\bigl[\|\bbX\|^\polycons\bigr]+\E\bigl[\sup_{t\in [0,T]}\|W_t\|^\polycons\bigr]\bigr)<\infty.
\end{equation}}
\argument{\lref{def: Phi}; the fact that $\sup_{x\in \R^d} ((1+\|x\|)^{-\polycons}|g(x)|)<\infty$ ;the assumption that $W$ is a standard Brownian motion}{that for all $(t,x)\in D$ there exists $c\in \R$ such that
\begin{equation}\llabel{evidence 2}
    \textstyle \E\bigl[|\Phi(t,x,W)|\bigr]=\E\bigl[|g(x,W_{T-t})|\bigr]\leq c\bigl(1+\E\bigl[\sup_{t\in [0,T]}\|W_t\|^\polycons\bigr]\bigr)<\infty.
\end{equation}}
\argument{\lref{evidence 2};\lref{def: Phi};\lref{def: phi};\lref{evidence 1}; the assumption that $(\bbT,\bbX)$ and $W$ are independent;\cite[Proposition 3.15]{Dereichnonconvergence2024} (applied\footnote{Note that for every set $\Omega$, every measurable space $(S, \cS)$, and every function $X \colon \Omega \to S$ it holds that $\sigma(X) = \{ A \subseteq \Omega \colon ( \exists\, B \in \cS \colon A = X^{ - 1 }( B ))\}$ (the smallest sigma-algebra on $\Omega$ with respect to which $X$ is measurable).} with $\Omega\curvearrowleft \Omega$, $\cF\curvearrowleft\cF$, $\P\curvearrowleft\P$, $D\curvearrowleft D$, $E\curvearrowleft E$, $\cG\curvearrowleft \sigma((\bbT,\bbX))$, $X\curvearrowleft (\bbT,\bbX)$, $Y\curvearrowleft W$, $\Phi \curvearrowleft\Phi$, $\phi\curvearrowleft\phi$ in the notation of \cite[Proposition 3.15]{Dereichnonconvergence2024})}{that it holds $\P$-a.s.\ that
\begin{equation}\llabel{eq4}
\begin{split}
     \E\bigl[g(\bbX+W_{T-\bbT})\big|(\bbT,\bbX)\bigr]=\E\bigl[\Phi(\bbX,\bbT,W)\big|(\bbT,\bbX)\bigr]=\phi(\bbT,\bbX).
     \end{split}
\end{equation}}
\argument{\lref{eq4};\lref{eq3};}{that it holds $\P$-a.s.\ that
\begin{equation}\llabel{eq5}
    \E\bigl[g(\bbX+W_{T-\bbT})\bigl|(\bbT,\bbX)\bigr]=u(\bbT,\bbX).
\end{equation}}
\argument{\lref{eq5}}{\lref{conclude}\dott}
\end{aproof}
\subsection{Error estimates for the Monte Carlo method}
\cfclear
\begin{athm}{lemma}{lem: montecarlo2}
    Let $p \in [2, \infty)$, $d, n \in \N$, let $(\Omega, \cF, \P)$ be a probability space, and let $X_i \colon \Omega\to
\R^d$, $i \in \{1, 2,\dots,n\}$, be \iid\ random variables with $\E[\|X_1\|]<\infty$. Then
\begin{equation}\llabel{conclude}
    \Bigl(\textstyle\E\bigl[\bigl\|\E[X_1]-\frac 1n\bigl(\textstyle\sum_{i=1}^nX_i\bigr)\bigr\|^p\bigr]\Bigr)^{1/p}\leq 2\Bigl[\displaystyle\frac{p-1}{n}\Bigr]^{1/2}\bigl(\E\bigl[\|X_1-\E[X_1]\|^p\bigr]\bigr)^{1/p}.
\end{equation}
\end{athm}
\cfclear
\begin{aproof}
    \argument{\unskip, \eg, \cite[Corollary 2.6]{MR4574851} }{\lref{conclude}\dott}
\end{aproof}
\subsection{Error estimates for the deep Kolmogorov method}\label{subsec: qualitative results}\label{susec: qualitative result}

In \cref{prop: linear II relu pinn draft} below we employ the a priori estimates for realizations of \ANNs\ in \cref{lem: anngrowth} and \cref{cor: ann lipschitz} from \cref{sec: priori estimate} together with \cref{cor: montecarlo 4}, \cref{Feynman-Kac formula}, and \cref{lem: montecarlo2} from this section to establish in \cref{conclude: prop: linear II relu pinn draft} below an abstract error estimate for the \DKM\ applied to the head \PDE\ in \cref{pde: prop: linear II relu pinn draft} without estimating the \ANN\ approximation error (see the first summand on the right hand side of \cref{conclude: prop: linear II relu pinn draft}).
\renewcommand{\bfF}{\mathbb{S}}
\renewcommand{\fF}{\mathbf{S}}
\renewcommand{\bfS}{\mathbf{F}}
\begin{athm}{prop}{prop: linear II relu pinn draft}
Let $ T,\kappa \in (0,\infty) $, $ d \in \N $, let $u\in C([0,T]\times\R^{d},\R)$ be at most polynomially growing, assume $u|_{ (0,T) \times \R^d } \in C^{ 1, 2 }( (0,T)\times\R^d, \R )$, assume for all $t\in (0,T)$, $x\in \R^d$ that
\begin{equation}\label{pde: prop: linear II relu pinn draft}
\textstyle
  \frac{ \partial }{ \partial t } u( t, x )
  +
 \frac{\kappa^2}{2}
  \Delta_x u(t,x)
  = 0,
\end{equation}
let $ ( \Omega, \mathcal{F}, \P ) $ be a probability space,
let $ \mathbb{T}_n \colon \Omega \to [0,T] $, $n\in \N$, be \iid\ random variables, 
let $ \mathbb{X}_n \colon \Omega \to \R^d $, $n\in \N$, be \iid\ random variables,
let $ W^n \colon [0,T] \times \Omega \to \R^d $, $n\in \N$, be \iid\ standard Brownian motions, let $\polycons\in [1,\infty)$ satisfy $\E[\|\bbX_1\|^{2\polycons}]+\sup_{x\in \R^d}((1+\|x\|)^{-\polycons}|u(T,x)|)<\infty$,
assume that $ (\mathbb{T}_n)_{n\in \N} $, $ (\mathbb{X}_n)_{n\in \N} $, and $ (W^n)_{n\in \N} $ are independent,
 for every $M=(M_1,M_2)\in \N^2$, $ L \in \N \backslash \{ 1 \} $,
$
  \ell  \in
  ( \{ d + 1 \} \times \N^{ L - 1 } \times \{ 1 \} )
$
 let $ \mathbb{F}_{M,\ell} \colon \R^{ \ffd( \ell ) } \to \R $
satisfy
for all
$ \theta \in \R^{ \ffd( \ell ) } $
that
\begin{equation}\llabel{def: bbF}
\textstyle
  \mathbb{F}_{M,\ell}( \theta )=
  \frac {1}{M_1} \sum_{m=1}^{M_1}
    \big|
      \reli{ \ell }{ L }{ \theta }( \mathbb{T}_m, \mathbb{X}_m )
      - \bigl[ \frac{ 1 }{ M_2 } \sum_{ n = 1 }^{ M_2 } u(T, \bbX_{ m } + \kappa W^{ m M_2 + n }_{ T-\bbT_m } ) \bigr]
    \big|^2, 
\end{equation} and for every $L\in \N\backslash\{1\}$, $\ell\in (\{d+1\}\times\N^{L-1}\times\{1\})$ let $\bfS_\ell\colon \R^{\fd(\ell)}\to [0,\infty]$ satisfy for all $\theta\in \R^{\fd(\ell)}$ that
\begin{equation}\llabel{def: bfS}
    \bfS_\ell(\theta)=\E\bigl[
    |
      \reli{ \ell }{ L }{ \theta }( \mathbb{T}_1, \mathbb{X}_1 )-u(\bbT_1,\bbX_1)|^2\bigr].
\end{equation}
Then there exists $c \in \R$ such that for every $R_1,R_2\in (1,\infty)$, $M=(M_1,M_2)\in \N^2$, $L\in \N\backslash\{1\}$, $\ell=(\ell_0,\ell_1,\dots,\ell_L)\in (\{d+1\}\times\N^{L-1}\times\{1\}) $ and every random variable $\vartheta \colon \Omega \to [-R_2,R_2]^{\fd(\ell)}$ it holds that
\begin{align}
&\textstyle
  \textstyle\E[\bfS_\ell(\vartheta)]\textstyle\leq \bigl[\inf_{\theta\in (-R_1,R_1)^{\fd(\ell)}}\bfS_\ell(\theta)\bigr]+\E\bigl[\textstyle \bbF_{M,\ell}( \vartheta ) - \inf_{\theta\in (-R_1,R_1)^{\fd(\ell)}} \bbF_{M,\ell}( \theta ) \bigr]\label{conclude: prop: linear II relu pinn draft}\\
  &+c^L(\fd(\ell))^{4\polycons} \bigl[\max\{R_1,R_2\}\textstyle\max\{\ell_1,\ell_2,\dots,\ell_{L-1}\}\bigr]^{2L+2}\bigl(\E\bigl[(1+\|\bbX_1\|)^{2\polycons(\fd(\ell))^2}\bigr]\bigr)^{(\fd(\ell))^{-2}}(M_1)^{-\nicefrac{1}{2}}\notag\ifnocf.  
\end{align}
\cfout[.]
\end{athm}
\begin{aproof}
Throughout this proof, let $U\colon [0,T]\times\R^d\to\R$ satisfy for all $t\in [0,T]$, $x\in \R^d$ that $U(t,x)=u(t,\kappa x)$, let $G\colon \R^d\to\R$ satisfy for all $x\in \R^d$ that $G(x)=u(T,\kappa x)$, for every $n\in \N$ let $\bfX_n\colon \Omega\to\R^d$ satisfy $\bfX_n=\kappa^{-1}\bbX_n$, for every $ L\in \N $,
$ \ell = ( \ell_0, \ell_1, \dots, \ell_L ) \in \N^{ L + 1 } $, $ \theta \in \R^{ \ffd( \ell ) } $, $\bbA\in C( \R,\R)$ let
  $\Reli{ \ell} { v}{ \theta }{\bbA}
  =
  (
    \Relii{ \ell}{ v }{ \theta }{ \bbA}{ 1 },
    \dots,
    \Relii{ \ell}{ v }{ \theta }{ \bbA}{\ell_v }
  )
  \colon \R^{ \ell_0 } \to \R^{ \ell_v }
$,
$ v \in \{ 1, 2, \dots, L \} $,
satisfy for all
$ v \in \{ 0, 1, \dots, L - 1 \} $,
$ x = (x_1, \dots, x_{ \ell_0 } ) \in \R^{ \ell_0 } $,
$ i \in \{ 1, 2, \dots, \ell_{ v + 1 } \} $
that
\begin{equation}
\llabel{realization multi general ac}
\begin{split}
 &\Relii{ \ell }{ v + 1 }{ \theta }{ \bbA}{i }( x )
  =
  \theta_{ \ell_{ v + 1 } \ell_v + i + \sum_{ h = 1 }^v \ell_h ( \ell_{ h - 1 } + 1 ) }
\\ &
  +
  \textstyle\sum_{ j = 1 }^{ \ell_v }
  \theta_{ (i - 1) \ell_v + j + \sum_{ h = 1 }^v \ell_h ( \ell_{ h - 1 } + 1 ) }
  \bigl(
    x_j
    \indicator{ \{ 0 \} }( v )
    +
 \bbA(\Relii{ \ell }{ \max\{ v, 1 \} }{ \theta }{ \bbA}{j }(x))
    \indicator{ \N }( v )
  \bigr) ,
\end{split}
\end{equation}
let $\cG\subseteq\cF$ satisfy
\begin{equation}\llabel{def: cG}
    \cG=\sigma((\bbT_1,\bfX_1)),
\end{equation}
let $\scrA_n\in C^{\min\{n,1\}}(\R,\R)$, $n\in \N_0$, satisfy that
\begin{enumerate}[label=(\roman*)]
    \item \llabel{item 1} it holds for all $x\in \R$ that $\limsup_{n\to\infty}\big[|\scrA_n(x)-\max\{x,0\}| + |(\scrA_n)'(x) - \mathbbm 1_{(0,\infty)}(x)|\big]=0$,
    \item \llabel{item 2} it holds for all $x\in \R$ that $\scrA_0(x)=\max\{x,0\}$,
    and 
    \item \llabel{item 3} it holds that $\sup_{n \in \N}\sup_{x\in\R}\allowbreak|(\scrA_n)'(x)|<\infty$
\end{enumerate}
(cf., \eg, \cite[Lemma 2.4]{MaArconvergenceproof}),
  for every $n\in \N_0$, $L\in \N\backslash\{1\}$, $M\in \N$, $\ell=(\ell_0,\ell_1,\dots,\ell_L)\in (\{d+1\}\times\N^{L-1}\times\{1\}) $ let $\fF_{M,\ell}^n\colon \R^{\fd(\ell)}\to\R$ satisfy for all $\theta\in \R^{\fd(\ell)}$ that
\begin{equation}\llabel{def: fF} 
\fF_{M,\ell}^n(\theta)= \textstyle\E\bigl[\big|
      \Reli{ \ell }{L }{ \theta }{ \scrA_n }( \mathbb{T}_1, \kappa \bfX_1 )
      - \bigl[ \frac{ 1 }{ M } \sum_{ m = 1 }^{ M } G( \bfX_{ 1 } + W^{ m }_{ T-\bbT_1 } ) \bigr]
    \big|^2\bigr],
\end{equation}
for every $n\in\N_0$, $L\in \N\backslash\{1\}$, $M=(M_1,M_2)\in \N^2$, $\ell=(\ell_0,\ell_1,\dots,\ell_L)\in (\{d+1\}\times\N^{L-1}\times\{1\}) $ let $\bbS^n_\ell\colon \R^{\fd(\ell)}\to\R$ satisfy for all $\theta\in \R^{\fd(\ell)}$ that
\begin{equation}\llabel{def: bfF} 
\bfF_{M,\ell}^n(\theta)= \textstyle
  \frac {1}{M_1} \sum_{m=1}^{M_1}
    \big|
      \Reli{ \ell }{ L }{\theta }{ \scrA_n}( \mathbb{T}_m, \kappa\mathbf{X}_m )
      - \bigl[ \frac{ 1 }{ M_2 } \sum_{ n = 1 }^{ M_2 } G( \bfX_{ m } + W^{ m M_2 + n }_{ T-\bbT_m } ) \bigr]
    \big|^2,
\end{equation}
and let $\bfW=(\bfW^{(1)},\dots,\bfW^{(d)})\colon[0,T]\times\Omega\to\R^d$ satisfy for all $t\in [0,T]$ that $\bfW_t=W_t^1$.
\argument{\cref{lem: anngrowth}; the fact that $\sup_{x\in \R^d} ((1+\|x\|)^{-\polycons}|G(x)|)<\infty$; the assumption that $\E[\|\bfX_1\|^{2\polycons}]<\infty$; the fact that $\bfW$ is a standard Brownian motion}{that for every $M\in \N$, $L\in \N\backslash\{1\}$, $\ell=(\ell_0,\ell_1,\dots,\allowbreak \ell_L)\allowbreak\in (\{d+1\}\times\N^{L-1}\times\{1\}) $ and every bounded set $D\subseteq \R^{\fd(\ell)}$ there exists $c\in \R$ such that for all $\theta\in D$ it holds that 
\begin{equation}\llabel{eq2p}
    \begin{split}
        \E\bigl[ |\reli{ \ell }{ L }{ \theta }( \mathbb{T}_1, \kappa\mathbf{X}_1 )|^2\bigr]\leq c(1+\E[\|\bbT_1\|^2]+\E[\|\bfX_1\|^2])<\infty 
    \end{split}
\end{equation}
and 
\begin{equation}\llabel{eq1p}
    \begin{split}
       \textstyle \E\Bigl[\bigl|\frac{ 1 }{ M } \sum_{ m = 1 }^{ M } G( \bfX_{ 1 } + W^{ m }_{ T-\bbT_1 } )\bigr|^2\Bigr]\leq c\bigl(1+\E[\|\bfX_1\|^{2\polycons}]+\E\bigl[\sup_{t\in [0,T]}\|\bfW_t\|^{2\polycons}\bigr]\bigr)<\infty \ifnocf.
    \end{split}
\end{equation}}
\argument{\lref{def: fF};\lref{eq1p}}{that for all $M\in \N$, $L\in \N\backslash\{1\}$, $\ell=(\ell_0,\ell_1,\dots,\allowbreak \ell_L)\allowbreak\in (\{d+1\}\times\N^{L-1}\times\{1\}) $, $\theta\in \R^{\fd(\ell)}$ it holds that
\begin{equation}\llabel{eq1}
\begin{split}
    \fF_{M,\ell}^0(\theta)&=\textstyle\E\!\Bigl[
    \big|
      \reli{ \ell }{ L }{ \theta }( \mathbb{T}_1, \kappa\mathbf{X}_1 )-\E\bigl[\frac{ 1 }{ M } \sum_{ m = 1 }^{ M } G( \bfX_{ 1 } + W^{ m }_{ T-\bbT_1 } )|\cG\bigr]\\
      &\textstyle\quad+\E\bigl[\frac{ 1 }{ M } \sum_{ m = 1 }^{ M } G( \bfX_{ 1 } + W^{ m }_{ T-\bbT_1 } )|\cG\bigr]
      -
     \bigl[\frac{ 1 }{ M } \sum_{ m = 1 }^{ M } G( \bfX_{ 1 } + W^{ m }_{ T-\bbT_1 } )\bigr]
    \big|^2
  \Bigr]\\
  &=\textstyle\E\bigl[
    \big|
      \reli{ \ell }{ L }{ \theta }( \mathbb{T}_1, \kappa\mathbf{X}_1 )-\E\bigl[\frac{ 1 }{ M } \sum_{ m = 1 }^{ M } G( \bfX_{ 1 } + W^{ m }_{ T-\bbT_1 } )|\cG\bigr]\big|^2\bigr]\\
      &+\textstyle\E\bigl[\big|\E\bigl[\frac{ 1 }{ M } \sum_{ m = 1 }^{ M } G( \bfX_{ 1 } + W^{ m }_{ T-\bbT_1 } )|\cG\bigr]
      -
     \bigl[\frac{ 1 }{ M } \sum_{ m = 1 }^{ M } G( \bfX_{ 1 } + W^{ m }_{ T-\bbT_1 } )\bigr]\big|^2\bigr].
  \end{split}
\end{equation}}
\argument{\cref{pde: prop: linear II relu pinn draft};the assumption that for all $t\in [0,T]$, $x\in \R^d$ it holds that $U(t,x)=u(t,\kappa x)$; the assumption that for all $x\in \R^d$ it holds that $G(x)=u(T,\kappa x)$}{that for all $t\in (0,T)$, $x\in \R^d$ it holds that
\begin{equation}\llabel{pde2}
    \textstyle
  \frac{ \partial }{ \partial t } U( t, x )
  +
 \frac{1}{2}
  \Delta_x U(t,x)
  = 0.
\end{equation}}
\argument{\lref{pde2};\lref{def: cG}; the fact that $ (\mathbb{T}_n)_{n\in \N} $, $ (\mathbf{X}_n)_{n\in \N} $, and $ (W^n)_{n\in \N} $ are independent; the fact that $\E[\|\bfX_1\|^{\polycons}]<\infty$; \cref{Feynman-Kac formula}}{that for all $M\in\N$ it holds $\P$-a.s.\ that
\begin{equation}\llabel{eq2}
\begin{split}
   \textstyle\E\bigl[\frac{ 1 }{ M } \sum_{ m = 1 }^{ M } G( \bfX_{ 1 } + W^{ m }_{ T-\bbT_1} )|\cG\bigr]&=\textstyle\frac{1}{M}\textstyle \sum_{m=1}^M\textstyle\E[G( \bfX_{ 1 } + W^{ m }_{ T-\bbT_1} )|(\bbT_1,\bfX_1)\bigr]\\
 &=\textstyle\frac{1}{M}\sum_{m=1}^M U(\bbT_1,\bfX_1)=U(\bbT_1,\bfX_1).
    \end{split}
\end{equation}}
\argument{\lref{eq1};\lref{eq2};}{for all $M\in \N$, $L\in \N\backslash\{1\}$, $\ell=(\ell_0,\ell_1,\dots,\allowbreak \ell_L)\allowbreak\in (\{d+1\}\times\N^{L-1}\times\{1\}) $, $\theta\in \R^{\fd(\ell)}$ that
\begin{equation}\llabel{eq3}
\begin{split}
    \fF_{M,\ell}^0(\theta)&=\textstyle\E\bigl[
    |
      \reli{ \ell }{ L }{ \theta }( \mathbb{T}_1, \kappa\mathbf{X}_1 )-U(\bbT_1,\bfX_1)|^2\bigr]\\
      &+\textstyle\E\Bigl[\big|U(\bbT_1,\bfX_1) 
      -
     \bigl[\frac{ 1 }{ M } \sum_{ m = 1 }^{ M } G( \bfX_{ 1 } + W^{ m }_{ T-\bbT_1 } )\bigr]\big|^2\Bigr].
     \end{split}
\end{equation}}
\argument{\lref{def: bfS};the fact that $\bbX_1=\kappa\bfX_1$; the assumption that for all $t\in [0,T]$, $x\in \R^d$ it holds that $U(t,x)=u(t,\kappa x)$}{that for all $L\in \N\backslash\{1\}$, $\ell\in (\{d+1\}\times\N^{L-1}\times\{1\})$, $\theta\in \R^{\fd(\ell)}$ it holds that
\begin{equation}\llabel{def: bfS2}
    \bfS_\ell(\theta)=\E\bigl[
    |
      \reli{ \ell }{ L }{ \theta }( \mathbb{T}_1, \kappa \mathbf{X}_1 )-U(\bbT_1,\bfX_1)|^2\bigr].
\end{equation} }
\argument{\lref{eq3};\lref{def: bfS2}}{that for all $M\in \N$, $R\in (0,\infty]$, $L\in \N\backslash\{1\}$, $\ell=(\ell_0,\ell_1,\dots,\allowbreak \ell_L)\allowbreak\in (\{d+1\}\times\N^{L-1}\times\{1\}) $ it holds that
\begin{equation}\llabel{eq4}
\begin{split}
\textstyle    \inf_{\theta\in (-R,R)^{\fd(\ell)}}\fF_{M,\ell}^0(\theta)&= \textstyle \inf_{\theta\in (-R,R)^{\fd(\ell)}}\bfS_\ell(\theta) \\
    &+ \textstyle\E\Bigl[\big|U(\bbT_1,\bfX_1) 
      -
     \bigl[\frac{ 1 }{ M } \sum_{ m = 1 }^{ M } G( \bfX_{ 1 } + W^{ m }_{ T-\bbT_1 } )\bigr]\big|^2\Bigr].
  \end{split}
\end{equation}}
\argument{\lref{eq4};}{for all $M\in \N$, $R\in (0,\infty]$, $L\in \N\backslash\{1\}$, $\ell=(\ell_0,\ell_1,\dots,\allowbreak \ell_L)\allowbreak\in (\{d+1\}\times\N^{L-1}\times\{1\}) $ that
\begin{equation}\llabel{eq5}
\begin{split}
   &\textstyle\E\Bigl[\big|U(\bbT_1,\bfX_1) 
      -
     \bigl[\frac{ 1 }{ M } \sum_{ m = 1 }^{ M } G( \bfX_{ 1 } + W^{ m }_{ T-\bbT_1 } )\bigr]\big|^2\Bigr]\\
     &\textstyle=\inf_{\theta\in (-R,R)^{\fd(\ell)}}\fF_{M,\ell}^0(\theta)- \inf_{\theta\in (-R,R)^{\fd(\ell)}}\bfS_\ell(\theta).
     \end{split}
\end{equation}}
\argument{\lref{eq5};\lref{eq3}}{for all $M\in \N$, $R\in (0,\infty]$, $L\in \N\backslash\{1\}$, $\ell=(\ell_0,\ell_1,\dots,\allowbreak \ell_L)\allowbreak\in (\{d+1\}\times\N^{L-1}\times\{1\}) $, $\vartheta\in \R^{\fd(\ell)}$ that
\begin{equation}\llabel{eq6.1}
\begin{split}
     \fF_{M,\ell}^0(\vartheta)&= \textstyle\E\bigl[
    |
      \reli{ \ell }{ L }{ \vartheta }( \mathbb{T}_1, \kappa\mathbf{X}_1 )-U(\bbT_1,\bfX_1)|^2\bigr]\\
      &+\bigl[\textstyle\inf_{\theta\in (-R,R)^{\fd(\ell)}}\fF_{M,\ell}^0(\theta)\bigr]- \bigl[\inf_{\theta\in (-R,R)^{\fd(\ell)}}\bfS_\ell(\theta)\bigr].
\end{split}
\end{equation}
}
\argument{\lref{eq6.1};}{for all $M\in \N$, $R\in (0,\infty]$, $L\in \N\backslash\{1\}$, $\ell=(\ell_0,\ell_1,\dots,\allowbreak \ell_L)\allowbreak\in (\{d+1\}\times\N^{L-1}\times\{1\}) $, $\vartheta\in \R^{\fd(\ell)}$ that
\begin{equation}\llabel{eq6.2}
\begin{split}
     \E\bigl[
    |
      \reli{ \ell }{ L }{ \vartheta }( \mathbb{T}_1, \kappa\mathbf{X}_1 )-U(\bbT_1,\bfX_1)|^2\bigr]&\leq \fF_{M,\ell}^0(\vartheta)\textstyle-\textstyle\inf_{\theta\in (-R,R)^{\fd(\ell)}}\fF_{M,\ell}^0(\theta)\\
      &+ \textstyle\inf_{\theta\in (-R,R)^{\fd(\ell)}}\bfS_\ell(\theta).
\end{split}
\end{equation}}
\argument{\lref{eq6.2};\lref{def: bfS}}{for all $M\in \N$, $R\in (0,\infty]$, $L\in \N\backslash\{1\}$, $\ell=(\ell_0,\ell_1,\dots,\allowbreak \ell_L)\allowbreak\in (\{d+1\}\times\N^{L-1}\times\{1\}) $, $\vartheta\in \R^{\fd(\ell)}$ that
\begin{equation}\llabel{eq6}
\begin{split}
     \bfS_\ell(\vartheta)\leq\textstyle  \bigl[\inf_{\theta\in (-R,R)^{\fd(\ell)}}\bfS_\ell(\theta)\bigr]+\textstyle \bigl(\fF_{M,\ell}^0(\vartheta)-\inf_{\theta\in (-R,R)^{\fd(\ell)}}\fF_{M,\ell}^0(\theta)\bigr).
\end{split}
\end{equation}}
\argument{\lref{eq6};the fact that for all $M=(M_1,M_2)\in \N^2$, $L\in \N\backslash\{1\}$, $\ell=(\ell_0,\ell_1,\dots,\allowbreak \ell_L)\allowbreak\in (\{d+1\}\times\N^{L-1}\times\{1\}) $ it holds that $\fF^0_{M_2,\ell}$ and $\bbF_{M,\ell}$ are continuous}{that for all $M\in \N^2$, $R_1\in (1,\infty]$, $R_2\in (1,\infty)$, $L\in \N\backslash\{1\}$, $\ell=(\ell_0,\ell_1,\dots,\allowbreak \ell_L)\allowbreak\in (\{d+1\}\times\N^{L-1}\times\{1\}) $, $\vartheta\in [-R_2,R_2]^{\fd(\ell)}$ it holds that
\begin{equation}\llabel{eq6main1}
\begin{split}
     &\bfS_\ell(\vartheta)\leq\textstyle \bigl[\inf_{\theta\in (-R_1,R_1)^{\fd(\ell)}}\bfS_\ell(\theta)\bigr]+\textstyle\bigl(\bbF_{M,\ell}(\vartheta)-\inf_{\theta\in (-R_1,R_1)^{\fd(\ell)}}\bbF_{M,\ell}(\theta)\bigr)\\
     &\textstyle+ \bigl[\sup_{\theta\in (-R_2,R_2)^{\fd(\ell)}}| \fF^0_{ M_2,\ell }( \theta ) - \bbF_{ M, \ell }( \theta ) |\bigr]+\bigl[\sup_{\theta\in (-R_1,R_1)^{\fd(\ell)}}| \fF^0_{ M_2,\ell }( \theta ) - \bbF_{ M, \ell }( \theta ) |\bigr].
\end{split}
\end{equation}}
\argument{\lref{def: bbF};\lref{eq2p};\lref{eq1p};the assumption that $\sup_{x\in \R^d}((1+\|x\|)^{-\polycons}|u(T,x)|)<\infty$; the fact that $\bfW$ is a standard Brownian motion;the assumption that $\P(\bbT_1\leq T)=1$}{that for all $R\in (0,\infty]$, $M\in \N^2$, $L\in \N\backslash\{1\}$, $\ell=(\ell_0,\ell_1,\dots,\allowbreak \ell_L)\allowbreak\in (\{d+1\}\times\N^{L-1}\times\{1\}) $ there exists $c\in (0,\infty)$ such that
\begin{equation}\llabel{eq6.1.1}
\begin{split}
&\textstyle\E\bigl[\inf_{\theta\in (-R,R)^{\fd(\ell)}}\bbF_{M,\ell}(\theta)\bigl]
  \\
  & \leq \textstyle \E[\bbF_{M,\ell}(0)]\leq c\bigl(1+\E[\|\bbX_1\|^{2\polycons}]+\E\bigl[\sup_{t\in [0,T]}\|\bfW_t\|^{2\polycons}\bigr]\bigr)<\infty.
   \end{split}
\end{equation}}
\argument{\lref{eq6.1.1};\lref{eq6main1};}{that for every $M\in \N^2$, $R_1\in (1,\infty]$, $R_2\in (1,\infty)$, $L\in \N\backslash\{1\}$, $\ell=(\ell_0,\ell_1,\dots,\allowbreak \ell_L)\allowbreak\in (\{d+1\}\times\N^{L-1}\times\{1\}) $ and every random variable $\vartheta\colon \Omega\to [-R_2,R_2]^{\fd(\ell)}$ it holds that 
\begin{equation}\llabel{eq6main}
\begin{split}
      &\E[\bfS_\ell(\vartheta)]\leq\textstyle \E\bigl[\inf_{\theta\in (-R_1,R_1)^{\fd(\ell)}}\bfS_\ell(\theta)\bigr]+\textstyle\E\bigl[\bbF_{M,\ell}(\vartheta)-\inf_{\theta\in (-R_1,R_1)^{\fd(\ell)}}\bbF_{M,\ell}(\theta)\bigr]\\
     &\textstyle+ \E\bigl[\sup_{\theta\in (-R_2,R_2)^{\fd(\ell)}}| \fF^0_{ M_2,\ell }( \theta ) - \bbF_{ M, \ell }( \theta ) |\bigr]+\E\bigl[\sup_{\theta\in (-R_1,R_1)^{\fd(\ell)}}| \fF^0_{ M_2,\ell }( \theta ) - \bbF_{ M, \ell }( \theta ) |\bigr].
     \end{split}
\end{equation}}
   \argument{\lref{def: fF};\lref{def: bfF};\cref{cor: montecarlo 4}}{that there exists $\fC\in\R$ which satisfies for all $n\in \N$, $M=(M_1,M_2)\in \N^2$, $L\in \N\backslash\{1\}$, $\ell\in (\{d+1\}\times\N^{L-1}\times\{1\})$, $p\in (|\fd(\ell)|^2,\infty)$, $R\in [\nicefrac{1}{2},\infty]$ that
    \begin{align}
        &\E\bigl[ \textstyle\sup_{\theta\in (-R,R)^{\fd(\ell)}}| \fF_{ M_2,\ell }^n( \theta ) - \bfF_{ M, \ell }^n( \theta ) |^p \bigr]\notag\\
        &=\E\bigl[ \textstyle\sup_{\theta\in (-R,R)^{\fd(\ell)}}| \bfF^n_{ M, \ell }( \theta ) - \E[ \bfF^n_{ M, \ell }( \theta ) ] |^p \bigr]\llabel{eq-1}\\
       & \textstyle\leq |R\fC|^{p} \E\!\left[
    \int_{ [-R,R]^{\fd(\ell) } }
     \bigl(| \bfF^n_{ M, \ell }( \theta ) - \E[ \bfF^n_{ M, \ell }( \theta ) ] |^{p}
    +
    \| ( \nabla \bfF^n_{M,\ell})( \theta ) -\E[( \nabla \bfF^n_{M,\ell})( \theta )]\|^{p}\bigr)
    \,\d \theta 
  \right]\notag\\
  &=\textstyle |R\fC|^{p}
    \int_{ [-R,R]^{\fd(\ell)} } \bigl(
     \E\bigl[| \bfF^n_{ M, \ell }( \theta ) - \E[ \bfF^n_{ M, \ell }( \theta ) ] |^{p}\bigr]
    +
    \E\bigl[\| ( \nabla \bfF^n_{M,\ell})( \theta ) -\E[( \nabla \bfF^n_{M,\ell})( \theta )]\|^{p}\bigr]\bigr)
    \,\d \theta \notag.
    \end{align}}
        \startnewargseq
     \argument{\lref{eq-1};}{that for all $n\in \N$, $M=(M_1,M_2)\in \N^2$, $L\in \N\backslash\{1\}$, $\ell\in (\{d+1\}\times\N^{L-1}\times\{1\})$, $p\in (|\fd(\ell)|^2,\infty)$, $R\in [\nicefrac{1}{2},\infty]$ it holds that
    \begin{equation}\llabel{argg1''}
         \begin{split}
        &\E\bigl[ \textstyle\sup_{\theta\in (-R,R)^{\fd(\ell)}}| \fF^n_{ M_2,\ell }( \theta ) - \bfF^n_{ M, \ell }( \theta ) |^p \bigr]\\
  &\textstyle\leq |R\fC|^{p}(2R)^{\fd(\ell) }\sup_{\theta\in (-R,R)^{\fd(\ell) }}\E\bigl[\bigl| \bfF^n_{ M, \ell }( \theta ) - \E[ \bfF^n_{ M, \ell }( \theta ) ] \bigr|^{p}\bigr]\\
  &+\textstyle|R\fC|^{p}(2R)^{\fd(\ell) }\sup_{\theta\in (-R
  ,R)^{\fd(\ell) }}\E\bigl[\bigl\| ( \nabla \bfF^n_{M,\ell})( \theta ) -\E[( \nabla \bfF^n_{M,\ell})( \theta )]\bigr\|^{p}\bigr].
        \end{split}
    \end{equation}}
       \argument{\cref{lem: montecarlo2};\lref{def: bfF}; the Cauchy--Schwarz inequality;the fact that for all $x,y\in \R^d$, $p\in [1,\infty)$ it holds that $\|x+y\|^p\leq 2^{p-1}\|x\|^p+2^{p-1}\|y\|^p$;the fact that for every $p\in [2,\infty)$, $\fd\in \N$ and every random variable $X\colon \Omega\to\R^\fd$ it holds that $\E\bigl[\|X-\E[X]\|^{p}\bigr]\leq 2^p\E\bigl[\|X\|^p\bigr]$}{that for all $n\in \N_0$, $M=(M_1,M_2)\in \N^2$, $L\in \N\backslash\{1\}$, $\ell\in (\{d+1\}\times\N^{L-1}\times\{1\})$, $\vartheta \in \R^{ \fd( \ell ) }$, $p\in [2,\infty)$ it holds that
    \begin{equation}\llabel{evidence 1}
    \begin{split}
        &\E\bigl[ |  \bfF^n_{ M, \ell }( \vartheta )-\E[\bfF^n_{ M, \ell }( \vartheta )] |^p \bigr]\\
&\leq 2^{p}\Bigl[\displaystyle\frac{p-1}{M_1}\Bigr]^{\nicefrac{p}{2}}\textstyle\textstyle\E\Bigl[\textstyle 
  \big|\bigl( \Reli{ \ell }{ L }{ \vartheta }{ \scrA_n }( \mathbb{T}_1, \kappa \mathbf{X}_1 )
      - \bigl[ \frac{ 1 }{ M_2 } \sum_{ m = 1 }^{ M_2 } G( \bfX_{ 1 } + W^{ m }_{ T-\bbT_1 } ) \bigr]\bigr)^2\\
      &\textstyle\quad-\E\bigl[\bigl(\Reli{ \ell }{ L }{ \vartheta }{ \scrA_n }( \mathbb{T}_1, \kappa \mathbf{X}_1 )
      - \bigl[ \frac{ 1 }{ M_2 } \sum_{ m = 1 }^{ M_2 } G( \bfX_{ 1 } + W^{ m }_{ T-\bbT_1 } ) \bigr]\bigr)^2\bigr]\big|^{p}\Bigr]  \\
&\leq 2^{2p}\Bigl[\displaystyle\frac{p-1}{M_1}\Bigr]^{\nicefrac{p}{2}}\textstyle\textstyle\E\Bigl[\textstyle 
  \bigl|\Reli{ \ell }{ L }{ \vartheta }{ \scrA_n }( \mathbb{T}_1, \kappa\mathbf{X}_1 )
      - \bigl[ \frac{ 1 }{ M_2 } \sum_{ m = 1 }^{ M_2 } G( \bfX_{ 1 } + W^{ m }_{ T-\bbT_1 } ) \bigr]\bigr|^{2p}\Bigr] \\
      &\leq 2^{4p-1} \Bigl[\displaystyle\frac{p-1}{M_1}\Bigr]^{\nicefrac{p}{2}}\Bigl(\E\bigl[\textstyle 
  |\Reli{ \ell }{ L }{ \vartheta }{ \scrA_n }( \mathbb{T}_1, \kappa\mathbf{X}_1 )|^{2p}\bigr]+\textstyle\E\Bigl[
       \bigl| \frac{ 1 }{ M_2 } \sum_{ m = 1 }^{ M_2 } G( \bfX_{ 1 } + W^{ m }_{ T-\bbT_1 } ) \bigr|^{2p}\Bigr] \Bigr).
\end{split}
    \end{equation}}
        \argument{\lref{def: bfF}}{that for all $n\in \N$, $L\in \N\backslash\{1\}$, $\ell=(\ell_0,\ell_1,\dots,\ell_L)\in (\{d+1\}\times\N^{L-1}\times\{1\}) $, $N=(N_1,N_2)\in \N^2$, $\theta\in \R^{\fd(\ell)}$ it holds that
 \begin{equation}\llabel{def: F1'}
 \begin{split}
       (\nabla \bfF^n_{N,\ell})(\theta)&=\textstyle\frac{ 1 }{ N_1 } \bigl(\sum_{ m = 1 }^{ N_1} 2 \bigl( \Reli{ \ell }{ L }{ \theta }{ \scrA_n }( \mathbb{T}_m, \kappa \mathbf{X}_m )
      - \bigl[ \frac{ 1 }{ N_2 } \sum_{ n = 1 }^{ N_2 } G( \bfX_{ m } + W^{ m N_2 + n }_{ T-\bbT_m } ) \bigr]\bigr)\\
      &\quad\cdot\bigl[\nabla_{ \theta } ( \Reli{ \ell}{L}{\theta }{ \scrA_n } (\bbT_m, \kappa\bfX_m ) ) \bigr]\bigr).
      \end{split}
\end{equation}}
     \argument{\cref{lem: montecarlo2};\lref{def: F1'};the Cauchy--Schwarz inequality;the fact that for all $x,y\in \R^d$, $p\in [1,\infty)$ it holds that $\|x+y\|^p\leq 2^{p-1}\|x\|^p+2^{p-1}\|y\|^p$; the fact that for every $p\in [2,\infty)$, $\fd\in \N$ and every random variable $X\colon \Omega\to\R^\fd$ it holds that $\E\bigl[\|X-\E[X]\|^{p}\bigr]\leq 2^p\E\bigl[\|X\|^p\bigr]$}{that for all $n\in \N$, $M=(M_1,M_2)\in \N^2$, $L\in \N\backslash\{1\}$, $\ell\in (\{d+1\}\times\N^{L-1}\times\{1\})$, $\vartheta \in \R^{ \fd( \ell ) }$, $p\in  [2,\infty)$ it holds that
    \begin{align}
        &\E\bigl[ \|  (\nabla\bfF^n_{ M, \ell })( \vartheta )-\E[(\nabla\bfF^n_{ M, \ell })( \vartheta )] \|^p \bigr]\notag\\
&\leq 2^{p}\Bigl[\displaystyle\frac{p-1}{M_1}\Bigr]^{\nicefrac{p}{2}}\textstyle\textstyle\E\bigl[\textstyle 
  \big\|2 \bigl( \Reli{ \ell }{ L }{ \vartheta }{ \scrA_n }( \mathbb{T}_1, \kappa\mathbf{X}_1 )
      - \bigl[ \frac{ 1 }{ M_2 } \sum_{ m = 1 }^{ M_2 } G( \bfX_{ 1 } + W^{m }_{ T-\bbT_1 } ) \bigr]\bigr)\bigl[\nabla_{ \vartheta }  (\Reli{ \ell}{L}{\vartheta }{ \scrA_n } (\bbT_1, \kappa\bfX_1 ) ) \bigr]\bigr)\notag \\
  &\quad-\textstyle\E\bigl[2 \bigl( \Reli{ \ell }{ L }{ \vartheta }{ \scrA_n }( \mathbb{T}_1, \kappa\mathbf{X}_1 )
      - \bigl[ \frac{ 1 }{ M_2 } \sum_{ m= 1 }^{ M_2 } G( \bfX_{ 1 } + W^{m }_{ T-\bbT_1 } ) \bigr]\bigr)\bigl[\nabla_{ \vartheta }  (\Reli{ \ell}{L}{\vartheta }{ \scrA_n } (\bbT_1, \kappa\bfX_1 ))  \bigr] \big\|^{p}\bigr]\notag  \\
&\leq 2^{3p}\Bigl[\displaystyle\frac{p-1}{M_1}\Bigr]^{\nicefrac{p}{2}}\textstyle\textstyle\E\bigl[\big\|  \bigl(\Reli{ \ell }{ L }{ \vartheta }{ \scrA_n }( \mathbb{T}_1, \kappa\mathbf{X}_1 )
      - \bigl[ \frac{ 1 }{ M_2 } \sum_{ m = 1 }^{ M_2 } G( \bfX_{ 1 } + W^{m }_{ T-\bbT_1 } ) \bigr]\bigr)\bigl[\nabla_{ \vartheta } ( \Reli{ \ell}{L}{\vartheta }{ \scrA_n } (\bbT_1, \kappa\bfX_1 )  )\bigr]\big\|^{p}\bigr]\notag\\
      &\leq 2^{3p+1} \Bigl[\displaystyle\frac{p-1}{M_1}\Bigr]^{\nicefrac{p}{2}}\textstyle\textstyle\E\bigl[\big| \Reli{ \ell }{ L }{ \vartheta }{ \scrA_n }( \mathbb{T}_1, \kappa\mathbf{X}_1 )
      - \bigl[ \frac{ 1 }{ M_2 } \sum_{ m = 1 }^{ M_2 } G( \bfX_{ 1 } + W^{m }_{ T-\bbT_1 } ) \bigr]\big|^{2p}\bigr]\notag\\
      &+2^{3p+1}\Bigl[\displaystyle\frac{p-1}{M_1}\Bigr]^{\nicefrac{p}{2}}\textstyle\textstyle\E\bigl[\big\| \nabla_{ \vartheta } ( \Reli{ \ell}{L}{\vartheta }{ \scrA_n } (\bbT_1, \kappa\bfX_1 )  )\big\|^{2p}\bigr] \llabel{evidence 2}\\
      &\leq 2^{5p} \Bigl[\displaystyle\frac{p-1}{M_1}\Bigr]^{\nicefrac{p}{2}}\textstyle\textstyle\E\bigl[| \Reli{ \ell }{ L }{ \vartheta }{ \scrA_n }( \mathbb{T}_1, \kappa\mathbf{X}_1 )|^{2p}\bigr]+2^{5p} \Bigl[\displaystyle\frac{p-1}{M_1}\Bigr]^{\nicefrac{p}{2}}\textstyle\textstyle\E\bigl[\big| \frac{ 1 }{ M_2 } \sum_{ m = 1 }^{ M_2 } G( \bfX_{ 1 } + W^{m }_{ T-\bbT_1 } ) \big|^{2p}\bigr]\notag\\
      &+2^{3p+1}\Bigl[\displaystyle\frac{p-1}{M_1}\Bigr]^{\nicefrac{p}{2}}\textstyle\textstyle\E\bigl[\big\| \nabla_{ \vartheta } ( \Reli{ \ell}{L}{\vartheta }{ \scrA_n } (\bbT_1, \kappa\bfX_1 )  )\big\|^{2p}\bigr]\notag.
    \end{align}}
\argument{the fact that $\sup_{x\in \R^d} (1+\|x\|)^{-\polycons}|G(x)|<\infty$; the assumption that $W^{m}$, $m\in \N$, are \iid\ standard Brownian motions;Doob's inequality;Cauchy--Schwarz inequality}{that there exist $c_1,c_2,c_3,c_4\in \R$ such that for all $M\in \N$, $p\in [2,\infty)$ it holds that
    \begin{equation}\llabel{evd1}
    \begin{split}
       &\textstyle \E\bigl[\big| \frac{ 1 }{ M } \sum_{ n = 1 }^{ M } G( \bfX_{ 1 } + W^{n }_{ T-\bbT_1 } ) \big|^{2p}\bigr]\\
       &\textstyle\leq \frac 1M\sum_{n=1}^M \E[|G( \bfX_{ 1 } + W^{n }_{ T-\bbT_1 } ) |^{2p}]=\E[|G( \bfX_{ 1 } + W^{1 }_{ T-\bbT_1 } ) |^{2p}]\\
       &\textstyle\leq |c_1|^{2p\polycons}+|c_1|^{2p\polycons}\E[\|\bfX_1\|^{2p\polycons}]+|c_1|^{2p\polycons}\E\bigl[\sup_{t\in [0,T]}\|\bfW_t\|^{2p\polycons}\bigr]\\
&=\textstyle|c_1|^{2p\polycons}+|c_1|^{2p\polycons}\E[\|\bfX_1\|^{2p\polycons}]+|c_1|^{2p\polycons}\E\bigl[\sup_{t\in [0,T]}\bigl[\bigl(\sum_{i=1}^d[\bfW^{(i)}_t]^2\bigr)^{p\polycons}\bigr]\bigr]\\
&\leq \textstyle|c_1|^{2p\polycons}+|c_1|^{2p\polycons}\E[\|\bfX_1\|^{2p\polycons}]+|c_1|^{2p\polycons}d^{pq-1}\E\bigl[\sup_{t\in [0,T]}\bigl[\sum_{i=1}^d[\bfW^{(i)}_t]^{2p\polycons}\bigr]\bigr]\\
       &\textstyle\leq |c_1|^{2p\polycons}+|c_1|^{2p\polycons}\E[\|\bfX_1\|^{2p\polycons}]+|c_1 |^{2p\polycons}d^{pq-1}\sum_{i=1}^d \E\bigl[\sup_{t\in [0,T]}|\bfW^{(i)}_t|^{2p\polycons}\bigr]\\
&=\textstyle|c_1|^{2p\polycons}+|c_1|^{2p\polycons}\E[\|\bfX_1\|^{2p\polycons}]+|c_1 |^{2p\polycons} d^{pq} \,\E\bigl[\sup_{t\in [0,T]}|\bfW^{(1)}_t|^{2p\polycons}\bigr]\\
       &\textstyle\leq |c_2|^{2p\polycons}+|c_2|^{2p\polycons} \E[\|\bfX_1\|^{2p\polycons}]+|c_2|^{2p\polycons}\Bigl(\frac{2p\polycons}{2p\polycons-1}\Bigr)^{2p\polycons}T^{p\polycons}(\ceil{2p\polycons})!!\\
       &\leq |2c_3|^{2p\polycons} \bigl(\E[\|\bfX_1\|^{2p\polycons}]+(2p\polycons)^{p\polycons}\bigr)\leq |c_4|^{2p\polycons} (p\polycons)^{p\polycons}\E\bigl[(1+\|\bfX_1\|)^{2p\polycons}\bigr].
       \end{split}
    \end{equation}}
    \argument{\cref{cor: ann lipschitz};\lref{item 3};}{that there exist $c_1,c_2\in(0,\infty)$ such that for all $R\in (1,\infty]$, $n\in \N$, $M=(M_1,M_2)\in \N^2$, $L\in \N\backslash\{1\}$, $\ell=(\ell_0,\ell_1,\dots,\ell_L)\in (\{d+1\}\times\N^{L-1}\times\{1\})$, $\vartheta=(\vartheta_1,\dots,\vartheta_{\fd(\ell)}) \in [-R,R]^{ \fd( \ell ) }$, $p\in  [2,\infty)$ it holds that
    \begin{align}
       & \E\bigl[\big\| \nabla_{ \vartheta } ( \Reli{ \ell}{L}{\vartheta }{ \scrA_n } (\bbT_1, \kappa\bfX_1 )  )\big\|^{2p}\bigr]\notag\\
       &\leq\E \biggl[\max\{1,\|(\bbT_1,\kappa\bfX_1)\|^{2p}\}\max\biggl\{1,\max_{j\in \{1,2,\dots,\fd(\ell)\}}|\vartheta_j|^{2pL}\biggr\}\biggl[3c_1\max\{\ell_0,\ell_1,\dots,\ell_{L-1}\}\biggr]^{2pL}\biggr]\notag\\
       &\textstyle\leq (c_2)^{2pL}R^{2pL}[\max\{\ell_0,\ell_1,\dots,\ell_{L-1}\}]^{2pL}\E\bigl[(1+\|\bfX_1\|)^{2p}\bigr]\notag\\
       &\textstyle = (c_2)^{2pL}[R\max\{\ell_0,\ell_1,\dots,\ell_{L-1}\}]^{2pL}\E\bigl[(1+\|\bfX_1\|)^{2p}\bigr]\llabel{evd2}.
    \end{align}}
 \argument{\cref{lem: anngrowth};\lref{item 3}}{that there exist $c_1,c_2\in(0,\infty)$ such that for all $R\in (1,\infty]$, $n\in \N$, $M=(M_1,M_2)\in \N^2$, $L\in \N\backslash\{1\}$, $\ell=(\ell_0,\ell_1,\dots,\ell_L)\in (\{d+1\}\times\N^{L-1}\times\{1\})$, $\vartheta=(\vartheta_1,\dots,\vartheta_{\fd(\ell)}) \in [-R,R]^{ \fd( \ell ) }$, $p\in  [2,\infty)$ it holds that
    \begin{align}
       & \E\bigl[\big| \Reli{ \ell}{L}{\vartheta }{ \scrA_n } (\bbT_1, \kappa\bfX_1 )  \big|^{2p}\bigr]\notag\\
       &\leq\E\biggl[ \max\{1,\|(\bbT_1,\kappa\bfX_1)\|^{2p}\} \max\biggl\{1,\max_{j\in \{1,2,\dots,\fd(\ell)\}}|\vartheta_j|^{2pL}\biggr\}\biggl[3c_1\max\{\ell_0,\ell_1,\dots,\ell_{L-1}\}\biggr]^{2pL}\biggr]\notag\\
       &\textstyle \leq (c_2)^{2pL}R^{2pL} [\max\{\ell_0,\ell_1,\dots,\ell_{L-1}\}]^{2pL}\E\bigl[(1+\|\bfX_1\|)^{2p}\bigr]\notag\\
       &\textstyle\leq (c_3)^{2pL}[R\max\{\ell_0,\ell_1,\dots,\ell_{L-1}\}]^{2pL}\E\bigl[(1+\|\bfX_1\|)^{2p}\bigr]\llabel{evd3}.
    \end{align}}
    \argument{\lref{evidence 1};\lref{evidence 2};\lref{evd1};\lref{evd2};\lref{evd3};the fact that $\polycons\geq 1$}{that there exists $c\in (0,\infty)$ such that for all $R\in (1,\infty]$, $n\in \N$, $M=(M_1,M_2)\in \N^2$, $L\in \N\backslash\{1\}$, $\ell=(\ell_0,\ell_1,\dots,\ell_L)\in (\{d+1\}\times\N^{L-1}\times\{1\})$, $\vartheta=(\vartheta_1,\dots,\vartheta_{\fd(\ell)}) \in [-R,R]^{ \fd( \ell ) }$, $p\in  [2,\infty)$ it holds that
    \begin{align}
       & \sup_{\vartheta \in [-R,R]^{ \fd( \ell ) }}\E\bigl[ |  \bfF^n_{ M, \ell }( \vartheta )-\E[\bfF^n_{ M, \ell }( \vartheta )] |^p \bigr]+\sup_{\vartheta \in [-R,R]^{ \fd( \ell ) }}\E\bigl[ \|  (\nabla\bfF^n_{ M, \ell })( \vartheta )-\E[(\nabla\bfF^n_{ M, \ell })( \vartheta )] \|^p \bigr]\notag\\
       &\leq c^{2pL}(p\polycons)^{p\polycons} [R\textstyle\max\{\ell_0,\ell_1,\dots,\ell_{L-1}\}]^{2pL}\E\bigl[(1+\|\bfX_1\|)^{2p\polycons}\bigr](M_1)^{\nicefrac{-p}{2}}\llabel{eqc1}.
    \end{align}}
    \argument{\lref{eqc1};\lref{argg1''}}{that there exists $c\in \R$ such that for all $n\in \N$, $M=(M_1,M_2)\in \N^2$, $L\in \N\backslash\{1\}$, $\ell=(\ell_0,\ell_1,\dots,\ell_L)\in (\{d+1\}\times\N^{L-1}\times\{1\})$, $R\in (1,\infty]$, $p\in [(\fd(\ell))^2,\infty)$ it holds that
    \begin{equation}\llabel{eqc2}
         \begin{split}
       & \E\bigl[ \textstyle\sup_{\theta\in (-R,R)^{\fd(\ell)}}| \fF^n_{ M_2,\ell }( \theta ) - \bfF^n_{ M, \ell }( \theta ) |^{p} \bigr]\\
       &\leq c^{2pLc}(p^2\polycons)^{p\polycons}[R\textstyle\max\{\ell_0,\ell_1,\dots,\ell_{L-1}\}]^{2p(L+1)}\E\bigl[(1+\|\bfX_1\|)^{2p\polycons}\bigr](M_1)^{-\nicefrac{p}{2}}.
        \end{split}
    \end{equation}}
     \argument{\lref{eqc2};\lref{item 1};\lref{item 2};\lref{item 3};the fact that for all $M\in \N^2$, $L\in \N\backslash\{1\}$, $\ell\in (\{d\}\times\N^{L-1}\times\{1\})$, $\theta\in \R^{\fd(\ell)}$ it holds that $\bbF_{M,\ell}(\theta)=\bfF_{M,\ell}^0(\theta)$}{that there exists $c\in \R$ such that for all $M=(M_1,M_2)\in \N^2$, $L\in \N\backslash\{1\}$, $\ell=(\ell_0,\ell_1,\dots,\ell_L)\in (\{d+1\}\times\N^{L-1}\times\{1\})$, $R\in (1,\infty]$, $p\in [(\fd(\ell))^2,\infty)$ it holds that
    \begin{equation}\llabel{eqc3}
         \begin{split}
     & \E\bigl[ \textstyle\sup_{\theta\in (-R,R)^{\fd(\ell)}}| \fF^0_{ M_2,\ell }( \theta ) - \bbF_{ M, \ell }( \theta ) |^{p} \bigr]\\
     &\leq 
       \E\bigl[ \textstyle\sup_{\theta\in (-R,R)^{\fd(\ell)}}| \fF^0_{ M_2,\ell }( \theta ) - \bfF^0_{ M, \ell }( \theta ) |^{p} \bigr]\\
       &\leq c^{2pLc}(p^2\polycons)^{p\polycons}[R\textstyle\max\{\ell_0,\ell_1,\dots,\ell_{L-1}\}]^{2p(L+1)}\E\bigl[(1+\|\bfX_1\|)^{2p\polycons}\bigr](M_1)^{-\nicefrac{p}{2}}.
        \end{split}
    \end{equation}}
     \argument{\lref{eqc3};Holder's inequality}{that there exists $c\in \R$ such that for all $M=(M_1,M_2)\in \N^2$, $L\in \N\backslash\{1\}$, $\ell=(\ell_0,\ell_1,\dots,\ell_L)\in (\{d+1\}\times\N^{L-1}\times\{1\})$, $R\in (1,\infty]$ it holds that
    \begin{equation}\llabel{eqc4}
         \begin{split}
     & \E\bigl[ \textstyle\sup_{\theta\in (-R,R)^{\fd(\ell)}}| \fF^0_{ M_2,\ell }( \theta ) - \bbF_{ M, \ell }( \theta ) | \bigr] \\
     &
     \leq c^{2L}(\fd(\ell))^{4\polycons}[R\textstyle\max\{\ell_0,\ell_1,\dots,\ell_{L-1}\}]^{2L+2}\bigl(\E\bigl[(1+\|\bfX_1\|)^{2\polycons(\fd(\ell))^2}\bigr]\bigr)^{(\fd(\ell))^{-2}}(M_1)^{-\nicefrac{1}{2}}.
        \end{split}
    \end{equation}}
    \argument{\lref{eqc4};the fact that for all $L\in \N\backslash\{1\}$, $\ell=(\ell_0,\ell_1,\dots,\ell_L)\in (\{d+1\}\times\N^{L-1}\times\{1\})$ it holds that $\max\{\ell_0,\ell_1,\dots,\ell_{L-1}\}\leq (d+1)\max\{\ell_1,\ell_2,\dots,\ell_{L-1}\}$}{that there exists $c\in \R$ which satisfies for all $M=(M_1,M_2)\in \N^2$, $L\in \N\backslash\{1\}$, $\ell=(\ell_0,\ell_1,\dots,\ell_L)\in (\{d+1\}\times\N^{L-1}\times\{1\})$, $R\in (1,\infty]$ that
    \begin{equation}\llabel{eqc5}
         \begin{split}
     & \E\bigl[ \textstyle\sup_{\theta\in (-R,R)^{\fd(\ell)}}| \fF^0_{ M_2,\ell }( \theta ) - \bbF_{ M, \ell }( \theta ) | \bigr] \\
     &
     \leq c^{2L}(\fd(\ell))^{4\polycons}[R\textstyle\max\{\ell_1,\ell_2,\dots,\ell_{L-1}\}]^{2L+2}\bigl(\E\bigl[(1+\|\bfX_1\|)^{2\polycons(\fd(\ell))^2}\bigr]\bigr)^{(\fd(\ell))^{-2}}(M_1)^{-\nicefrac{1}{2}}.
        \end{split}
    \end{equation}}
    \argument{\lref{eqc5};\lref{eq6main};}{\cref{conclude: prop: linear II relu pinn draft}\dott}
\end{aproof}

In \cref{prop: linear II relu pinn} below, which is the main result of this work, we combine the conclusion of \cref{prop: linear II relu pinn draft} above with the \ANN\ approximation result in \cref{cor: Pinns for Poisson3 intro2 pre2 approximation} from \cref{sec: ANN approximation} to establish in \cref{conclude: prop: linear II relu pinn} an error analysis for the \DKM\ when applied to the heat \PDE\ in \cref{pde: prop: linear II relu pinn}.
\begin{athm}{theorem}{prop: linear II relu pinn}
Let $ T,\kappa \in (0,\infty) $, $ d \in \N $, let $u\in C^2([0,T]\times\R^{d},\R)$ be at most polynomially growing, assume for all $t\in (0,T)$, $x\in \R^d$ that
\begin{equation}\label{pde: prop: linear II relu pinn}
\textstyle
  \frac{ \partial }{ \partial t } u( t, x )
  +
 \frac{\kappa^2}{2}
  \Delta_x u(t,x)
  = 0,
\end{equation}
let $\polycons\in [1,\infty)$ satisfy $\sup_{x\in \R^d}((1+\|x\|)^{-\polycons}|u(T,x)|)<\infty$,
 let $ ( \Omega, \mathcal{F}, \P ) $ be a probability space,
 let $ \mathbb{T}_n \colon \Omega \to [0,T] $, $n\in \N$, be independent uniformly distributed random variables, let $D\subseteq\R^d$ be open and bounded,
let $ \mathbb{X}_n \colon \Omega \to D $, $n\in \N$, be independent uniformly distributed random variables,
let $ W^n \colon [0,T] \times \Omega \to \R^d $, $n\in \N$, be \iid\ standard Brownian motions, 
assume that $ (\mathbb{T}_n)_{n\in \N} $, $ (\mathbb{X}_n)_{n\in \N} $, and $ (W^n)_{n\in \N} $ are independent,
and for every $M=(M_1,M_2)\in \N^2$, $ L \in \N \backslash \{ 1 \} $,
$
  \ell  \in
  ( \{ d + 1 \} \times \N^{ L - 1 } \times \{ 1 \} )
$
 let $ \mathbb{F}_{M,\ell} \colon \R^{ \ffd( \ell ) } \to \R $
satisfy
for all
$ \theta \in \R^{ \ffd( \ell ) } $
that
\begin{equation}\llabel{def: bbF}
\textstyle
  \mathbb{F}_{M,\ell}( \theta )=
  \frac {1}{M_1} \sum_{m=1}^{M_1}
    \big|
      \reli{ \ell }{ L }{ \theta }( \mathbb{T}_m,  \mathbb{X}_m )
      - \bigl[ \frac{ 1 }{ M_2 } \sum_{ n = 1 }^{ M_2 } u(T, \bbX_{ m } + \kappa W^{ m M_2 + n }_{ T-\bbT_m } ) \bigr]
    \big|^2\ifnocf.
\end{equation}
\cfout[.]
Then there exists $c \in \R$ such that for every $R_1,R_2\in (1,\infty)$, $M=(M_1,M_2)\in \N^2$, $L\in \N\backslash\{1\}$, $\ell=(\ell_0,\ell_1,\dots,\ell_L)\in (\{d+1\}\times\N^{L-1}\times\{1\}) $ and every random variable $\vartheta \colon \Omega \to [-R_2,R_2]^{\fd(\ell)}$ it holds\footnote{Note that for every $n \in \N$ and every open $O \subseteq \R^n$ it holds that $\lambda( O ) \in [0,\infty]$ is the Lebesgue measure of $O$.} that
\begin{align}
&\textstyle
  \E\bigl[
    \int_{[0,T] \times D}
      | u(y) - \reli{ \ell }{ L}{ \vartheta }( y ) |^2
    \, \d y
  \bigr]\notag\\
  &\leq c\mathbbm 1_{(R_1,\infty)}(c)+ T\lambda( D)\,\E\bigl[\textstyle \bbF_{M,\ell}( \vartheta ) - \inf_{\theta\in (-R_1,R_1)^{\fd(\ell)}} \bbF_{M,\ell}( \theta ) \bigr]\label{conclude: prop: linear II relu pinn}
  \\&+\textstyle c^L(\fd(\ell))^{4\polycons} \bigl[\max\{R_1,R_2\}\textstyle\max\{\ell_1,\ell_2,\dots,\ell_{L-1}\}\bigr]^{2L+2}(M_1)^{-\nicefrac{1}{2}} +c\bigl[\min\{\ell_1, \ell_2, \dots, \ell_{ L - 1 }\}\bigr]^{-\frac{2}{d+5}}\notag\ifnocf.  
\end{align}
\cfout[.]
\end{athm}
\begin{aproof}
Throughout this proof, for every $L\in \N\backslash\{1\}$, $\ell\in (\{d+1\}\times\N^{L-1}\times\{1\})$ let $\bfS_\ell\colon \R^{\fd(\ell)}\to [0,\infty]$ satisfy for all $\theta\in \R^{\fd(\ell)}$ that
\begin{equation}\llabel{def: bfS}
    \bfS_\ell(\theta)=\E\bigl[
    |
      \reli{ \ell }{ L }{ \theta }( \mathbb{T}_1, \mathbb{X}_1 )-u(\bbT_1,\bbX_1)|^2\bigr].
\end{equation}
\argument{\lref{def: bfS};\cref{prop: linear II relu pinn draft};the fact that $\bbX_1$ is a bounded random variable}{there exists $c \in (0,\infty)$ such that for every $R_1,R_2\in (1,\infty)$, $M=(M_1,M_2)\in \N^2$, $L\in \N\backslash\{1\}$, $\ell=(\ell_0,\ell_1,\dots,\ell_L)\in (\{d+1\}\times\N^{L-1}\times\{1\}) $ and every random variable $\vartheta \colon \Omega \to [-R_2,R_2]^{\fd(\ell)}$ it holds that
\begin{equation}\llabel{def: c}
\begin{split}
\textstyle
  \textstyle\E[\bfS_\ell(\vartheta)]\textstyle &\textstyle\leq \bigl[\inf_{\theta\in (-R_1,R_1)^{\fd(\ell)}}\bfS_\ell(\theta)\bigr]+\E\bigl[\textstyle \bbF_{M,\ell}( \vartheta ) - \inf_{\theta\in (-R_1,R_1)^{\fd(\ell)}} \bbF_{M,\ell}( \theta ) \bigr]\\
  &+c^L(\fd(\ell))^{4\polycons} \bigl[\max\{R_1,R_2\}\textstyle\max\{\ell_1,\ell_2,\dots,\ell_{L-1}\}\bigr]^{2L+2}(M_1)^{-\nicefrac{1}{2}}
\end{split}
\end{equation}}
\argument{\lref{def: bfS};the fact that $(\bbT_1,\bbX_1)$ is uniformly distributed on $[0,T]\times D$;}{that for all $L\in \N\backslash\{1\}$, $\ell=(\ell_0,\ell_1,\dots,\allowbreak \ell_L)\allowbreak\in (\{d+1\}\times\N^{L-1}\times\{1\}) $, $\theta\in \R^{\fd(\ell)}$ it holds that
\begin{equation}\llabel{eq1}
\begin{split}
   \bfS_\ell(\theta)= \frac{1}{T\lambda(D)}\int_{[0,T] \times D} |u(y)-\reli{\ell}{L}{\theta}(y)|^2
\,\d y.
\end{split}
\end{equation}}
\argument{\lref{eq1};\lref{def: c}}{that there exists $c\in(0,\infty)$ such that for every $R_1,R_2\in (1,\infty)$, $M=(M_1,M_2)\in \N^2$, $L\in \N\backslash\{1\}$, $\ell=(\ell_0,\ell_1,\dots,\ell_L)\in (\{d+1\}\times\N^{L-1}\times\{1\}) $ and every random variable $\vartheta \colon \Omega \to [-R_2,R_2]^{\fd(\ell)}$ it holds that
\begin{align}
\textstyle
 & \textstyle \E\bigl[
    \int_{[0,T] \times D}
      | u(y) - \reli{ \ell }{ L}{ \vartheta }( y ) |^2
    \, \d y
  \bigr]\textstyle \notag\\
  &\textstyle\leq \bigl[\inf_{\theta\in (-R_1,R_1)^{\fd(\ell)}}\bigl(T\lambda(D)\bfS_\ell(\theta)\bigr)\bigr]+T\lambda(D)\E\bigl[\textstyle \bbF_{M,\ell}( \vartheta ) - \inf_{\theta\in (-R_1,R_1)^{\fd(\ell)}} \bbF_{M,\ell}( \theta ) \bigr]\notag\\
  &+T\lambda( D)c^L(\fd(\ell))^{4\polycons} \bigl[\max\{R_1,R_2\}\textstyle\max\{\ell_1,\ell_2,\dots,\ell_{L-1}\}\bigr]^{2L+2}(M_1)^{-\nicefrac{1}{2}} \llabel{def: c2}. 
\end{align}} 
\argument{the assumption that $D$ is open and bounded;\cref{cor: Pinns for Poisson3 intro2 pre2 approximation} (applied with $d\curvearrowleft d$, $D\curvearrowleft [0,T] \times D$, $u\curvearrowleft u|_{[0,T]\times\bar{D}}$ in the notation of \cref{cor: Pinns for Poisson3 intro2 pre2 approximation});}{that there exists $c\in (0,\infty)$ which satisfies for all $R\in (0,\infty]$, $L\in \N\backslash\{1\}$, $\ell=(\ell_0,\ell_1,\dots,\allowbreak \ell_L)\allowbreak\in (\{d+1\}\times\N^{L-1}\times\{1\}) $ that
\begin{equation}\llabel{eq1'}
\begin{split}
\textstyle    \inf_{\theta\in (-R,R)^{\fd(\ell)}}\int_{[0,T] \times D} |u(y)-\reli{\ell}{L}{\theta}(y)|^2&\leq  c\mathbbm 1_{(-\infty,c)}(R)+ c  [\min\{  \ell_1, \ell_2, \dots, \ell_{ L - 1 } \}]^{\frac{-2}{d+5}}.
  \end{split}
\end{equation}}
\startnewargseq
\argument{\lref{eq1};\lref{eq1'}}{that for all $R\in (0,\infty]$, $L\in \N\backslash\{1\}$, $\ell=(\ell_0,\ell_1,\dots,\allowbreak \ell_L)\allowbreak\in (\{d+1\}\times\N^{L-1}\times\{1\}) $ it holds that
\begin{equation}\llabel{eq2}
\begin{split}
\textstyle    \inf_{\theta\in (-R,R)^{\fd(\ell)}}\bigl(T\lambda( D)\bfS_{\ell}(\theta)\bigr)&\leq  c \mathbbm 1_{(R,\infty)}(c)+ c [\min\{  \ell_1, \ell_2, \dots, \ell_{ L - 1 } \}]^{\frac{-2}{d+5}}.
  \end{split}
\end{equation}}
\argument{\lref{def: c2};\lref{eq2}}{\cref{conclude: prop: linear II relu pinn}\dott}
\end{aproof}
In \cref{cor: linear II relu pinn} below we bound the indicator term $c \mathbbm 1_{ (R_1, \infty )}( c )$ in \cref{conclude: prop: linear II relu pinn} (see the first summand on the right hand side of \cref{conclude: prop: linear II relu pinn}) from above through the term $c \exp( -R_1 ) \mathbbm 1_{ (-\infty,c) }( R_1 )$ with exponential decay in $R_1$ (see the second summand on the right hand side of \cref{conclude: cor: linear II relu pinn}).
\begin{athm}{cor}{cor: linear II relu pinn}
Let $ T,\kappa \in (0,\infty) $, $ d \in \N $, let $u\in C^2([0,T]\times\R^{d},\R)$ be at most polynomially growing, assume for all $t\in (0,T)$, $x\in \R^d$ that
\begin{equation}\llabel{pde}
\textstyle
  \frac{ \partial }{ \partial t } u( t, x )
  +
  \frac{ \kappa^2}{ 2 }
  \Delta_x u(t,x)
  = 0,
\end{equation}
let $\polycons\in [1,\infty)$ satisfy $\sup_{x\in \R^d}((1+\|x\|)^{-\polycons}|u(T,x)|)<\infty$,
let $ ( \Omega, \mathcal{F}, \P ) $ be a probability space,
let $ \mathbb{T}_n \colon \Omega \to [0,T] $, $n\in \N$, be independent uniformly distributed random variables, let $D\subseteq\R^d$ be open and bounded,
let $ \mathbb{X}_n \colon \Omega \to D $, $n\in \N$, be independent uniformly distributed random variables,
let $ W^n \colon [0,T] \times \Omega \to \R^d $, $n\in \N$, be \iid\ standard Brownian motions,  
assume that $ (\mathbb{T}_n)_{n\in \N} $, $ (\mathbb{X}_n)_{n\in \N} $, and $ (W^n)_{n\in \N} $ are independent,
 and for every $M=(M_1,M_2)\in \N^2$, $ L \in \N \backslash \{ 1 \} $,
$
  \ell  \in
  ( \{ d + 1 \} \times \N^{ L - 1 } \times \{ 1 \} )
$
 let $ \mathbb{F}_{M,\ell} \colon \R^{ \ffd( \ell ) } \to \R $
satisfy
for all
$ \theta \in \R^{ \ffd( \ell ) } $
that
\begin{equation}\llabel{def: bbF}
\textstyle
  \mathbb{F}_{M,\ell}( \theta )=
  \frac {1}{M_1} \sum_{m=1}^{M_1}
    \big|
      \reli{ \ell }{ L }{ \theta }( \mathbb{T}_m,  \mathbb{X}_m )
      - \bigl[ \frac{ 1 }{ M_2 } \sum_{ n = 1 }^{ M_2 } u(T,  \bbX_{ m } +\kappa W^{ m M_2 + n }_{ T-\bbT_m } ) \bigr]
    \big|^2\ifnocf.
\end{equation}
\cfout[.]
Then there exists $c \in \R$ such that for every $R_1,R_2\in (1,\infty)$, $M=(M_1,M_2)\in \N^2$, $L\in \N\backslash\{1\}$, $\ell=(\ell_0,\ell_1,\dots,\ell_L)\in (\{d+1\}\times\N^{L-1}\times\{1\}) $ and every random variable $\vartheta \colon \Omega \to [-R_2,R_2]^{\fd(\ell)}$ it holds that
\begin{align}
&\textstyle
  \E\bigl[
    \int_{[0,T] \times D}
      | u(y) - \reli{ \ell }{ L}{ \vartheta }( y ) |^2
    \, \d y
  \bigr]\notag\\
  &\textstyle\leq c^L(\fd(\ell))^{4\polycons} \bigl[\max\{R_1,R_2\}\textstyle\max\{\ell_1,\ell_2,\dots,\ell_{L-1}\}\bigr]^{2L+2}(M_1)^{-\nicefrac{1}{2}}+c\exp(-R_1)\mathbbm 1_{(-\infty,c)}(R_1)\notag\\
  &+c [\min\{\ell_1, \ell_2, \dots, \ell_{ L - 1 } \}]^{\frac{-2}{d+5}}+T\lambda(D)\,\E\bigl[\textstyle \bbF_{M,\ell}( \vartheta ) - \inf_{\theta\in (-R_1,R_1)^{\fd(\ell)}} \bbF_{M,\ell}( \theta ) \bigr] \label{conclude: cor: linear II relu pinn}  \ifnocf.
\end{align}
\cfout[.]
\end{athm}
\begin{aproof}
    \argument{\cref{prop: linear II relu pinn};}{that there exists $c \in (0,\infty)$ which satisfies for every $R_1,R_2\in (1,\infty)$, $M=(M_1,M_2)\in \N^2$, $L\in \N\backslash\{1\}$, $\ell=(\ell_0,\ell_1,\dots,\ell_L)\in (\{d+1\}\times\N^{L-1}\times\{1\}) $ and every random variable $\vartheta \colon \Omega \to [-R_2,R_2]^{\fd(\ell)}$ that
    \begin{align}
&\textstyle
  \E\bigl[
    \int_{[0,T]\times D}
      | u(y) - \reli{ \ell }{ L}{ \vartheta }( y ) |^2
    \, \d y
  \bigr]\notag\\
  &\leq c\mathbbm 1_{(R_1,\infty)}(c)+ T\lambda( D)\,\E\bigl[\textstyle \bbF_{M,\ell}( \vartheta ) - \inf_{\theta\in (-R_1,R_1)^{\fd(\ell)}} \bbF_{M,\ell}( \theta ) \bigr]\llabel{eq1}
  \\&+\textstyle c^L(\fd(\ell))^{4\polycons} \bigl[\max\{R_1,R_2\}\textstyle\max\{\ell_1,\ell_2,\dots,\ell_{L-1}\}\bigr]^{2L+2}(M_1)^{-\nicefrac{1}{2}} +c\bigl[\min\{\ell_1, \ell_2, \dots, \ell_{ L - 1 }\}\bigr]^{-\frac{2}{d+5}}\notag  
\end{align}}
\startnewargseq
\argument{the fact that for all $x\in (-\infty,c)$ it holds that $\exp(-x)\geq \exp(-c)$}{that for all $R\in (0,\infty)$ it holds that
\begin{equation}\llabel{eq2}
    \mathbbm 1_{(R,\infty)}(c)=\mathbbm 1_{(-\infty,c)}(R)\leq \exp(c)\exp(-R)\mathbbm 1_{(-\infty,c)}(R).
\end{equation}}
\argument{\lref{eq2};\lref{eq1}}{\cref{conclude: cor: linear II relu pinn}\dott}
\end{aproof}
In \cref{main cor} below we entirely avoid the parameter $R_1$ in \cref{conclude: cor: linear II relu pinn} in \cref{cor: linear II relu pinn} and can then also replace $R_2$ by $R$ (cf.\ \cref{conclude: cor: linear II relu pinn} in \cref{cor: linear II relu pinn} with \cref{conclude: main cor} in \cref{main cor}). \cref{main theorem} in the introduction is an immediate consequence of \cref{main cor}.
\begin{athm}{cor}{main cor}
Let $ T,\kappa \in (0,\infty) $, $ d \in \N $, let $u\in C^2([0,T]\times\R^{d},\R)$ be at most polynomially growing, assume for all $t\in (0,T)$, $x\in \R^d$ that
\begin{equation}\llabel{pde}
\textstyle
  \frac{ \partial }{ \partial t } u( t, x )
  +
  \frac{ \kappa^2}{ 2 }
  \Delta_x u(t,x)
  = 0,
\end{equation}
 let $ ( \Omega, \mathcal{F}, \P ) $ be a probability space,
let $ \mathbb{T}_n \colon \Omega \to [0,T] $, $n\in \N$, be independent uniformly distributed random variables, let $D\subseteq\R^d$ be open and bounded,
let $ \mathbb{X}_n \colon \Omega \to D $, $n\in \N$, be independent uniformly distributed random variables,
let $ W^n \colon [0,T] \times \Omega \to \R^d $, $n\in \N$, be \iid\ standard Brownian motions,  
assume that $ (\mathbb{T}_n)_{n\in \N} $, $ (\mathbb{X}_n)_{n\in \N} $, and $ (W^n)_{n\in \N} $ are independent, 
and for every $M=(M_1,M_2)\in \N^2$, $ L \in \N \backslash \{ 1 \} $,
$
  \ell  \in
  ( \{ d + 1 \} \times \N^{ L - 1 } \times \{ 1 \} )
$
 let $ \mathbb{F}_{M,\ell} \colon \R^{ \ffd( \ell ) } \to \R $
satisfy
for all
$ \theta \in \R^{ \ffd( \ell ) } $
that
\begin{equation}\llabel{def: bbF}
\textstyle
  \mathbb{F}_{M,\ell}( \theta )=
  \frac {1}{M_1} \sum_{m=1}^{M_1}
    \big|
      \reli{ \ell }{ L }{ \theta }( \mathbb{T}_m,  \mathbb{X}_m )
      - \bigl[ \frac{ 1 }{ M_2 } \sum_{ n = 1 }^{ M_2 } u(T, \bbX_{ m } +\kappa W^{ m M_2 + n }_{ T-\bbT_m } ) \bigr]
    \big|^2\ifnocf.
\end{equation}
\cfout[.]
Then there exists $c \in (0,\infty)$ such that for every $R\in (1,\infty)$, $M=(M_1,M_2)\in \N^2$, $L\in \N\backslash\{1\}$, $\ell=(\ell_0,\ell_1,\dots,\ell_L)\in (\{d+1\}\times\N^{L-1}\times\{1\}) $ and every random variable $\vartheta \colon \Omega \to [-R,R]^{\fd(\ell)}$ it holds that
\begin{equation}\label{conclude: main cor}
\begin{split}
&\textstyle
  \E\bigl[
    \int_{[0,T] \times D}
      | u(y) - \reli{ \ell }{ L}{ \vartheta }( y ) |^2
    \, \d y
  \bigr]\textstyle\leq [R+\textstyle\max\{\ell_1,\ell_2,\dots,\ell_{L-1}\}]^{cL}(M_1)^{-\nicefrac{1}{2}}\\
  &+c [\min\{\ell_1, \ell_2, \dots, \ell_{ L - 1 } \}]^{-2/(d+5)}+T\lambda(D)\,\E\bigl[\textstyle \bbF_{M,\ell}( \vartheta ) - \inf_{\theta\in (-c,c)^{\fd(\ell)}} \bbF_{M,\ell}( \theta ) \bigr]  \ifnocf.
  \end{split}
\end{equation}
\cfout[.]
\end{athm}
\begin{aproof}
Throughout this proof, let $\polycons\in [1,\infty)$ satisfy 
\begin{equation}\llabel{def: q}
     \textstyle\sup_{x\in \R^d}((1+\|x\|)^{-\polycons}|u(T,x)|)<\infty.
\end{equation}
    \argument{\lref{def: q};\cref{cor: linear II relu pinn};}{that there exists $c\in (0,\infty)$ such that for every $R_1,R_2\in (1,\infty)$, $M=(M_1,M_2)\in \N^2$, $L\in \N\backslash\{1\}$, $\ell=(\ell_0,\ell_1,\dots,\ell_L)\in (\{d+1\}\times\N^{L-1}\times\{1\}) $ and every random variable $\vartheta \colon \Omega \to [-R_2,R_2]^{\fd(\ell)}$ it holds that
\begin{align}
&\textstyle
  \E\bigl[
    \int_{[0,T] \times D}
      | u(y) - \reli{ \ell }{ L}{ \vartheta }( y ) |^2
    \, \d y
  \bigr]\notag \\
  &\textstyle\leq c^L(\fd(\ell))^{4\polycons} \bigl[\max\{R_1,R_2\}\textstyle\max\{\ell_1,\ell_2,\dots,\ell_{L-1}\}\bigr]^{2L+2}(M_1)^{-\nicefrac{1}{2}}+c\exp(-R_1)\mathbbm 1_{(-\infty,c)}(R_1) \notag\\
  &+c [\min\{\ell_1, \ell_2, \dots, \ell_{ L - 1 } \}]^{\frac{-2}{d+5}}+T\lambda( D)\,\E\bigl[\textstyle \bbF_{M,\ell}( \vartheta ) - \inf_{\theta\in (-R_1,R_1)^{\fd(\ell)}} \bbF_{M,\ell}( \theta ) \bigr]\llabel{eq1'}  \ifnocf.
\end{align}}
    \argument{\lref{eq1'}}{that there exists $c\in (0,\infty)$ such that for every $R\in (1,\infty)$, $M=(M_1,M_2)\in \N^2$, $L\in \N\backslash\{1\}$, $\ell=(\ell_0,\ell_1,\dots,\ell_L)\in (\{d+1\}\times\N^{L-1}\times\{1\}) $ and every random variable $\vartheta \colon \Omega \to [-R,R]^{\fd(\ell)}$ it holds that
\begin{equation}\llabel{eq1}
\begin{split}
&\textstyle
  \E\bigl[
    \int_{[0,T] \times D}
      | u(y) - \reli{ \ell }{ L}{ \vartheta }( y ) |^2
    \, \d y
  \bigr]\\
  &\textstyle\leq c^L(\fd(\ell))^{4\polycons} \bigl[\max\{c+1,R\}\textstyle\max\{\ell_1,\ell_2,\dots,\ell_{L-1}\}\bigr]^{2L+2}(M_1)^{-\nicefrac{1}{2}} \\
  &+c [\min\{\ell_1, \ell_2, \dots, \ell_{ L - 1 } \}]^{\frac{-2}{d+5}}+T\lambda( D)\,\E\bigl[\textstyle \bbF_{M,\ell}( \vartheta ) - \inf_{\theta\in (-c-1,c+1)^{\fd(\ell)}} \bbF_{M,\ell}( \theta ) \bigr]  \ifnocf.
  \end{split}
\end{equation}}
\argument{\lref{eq1};}{that there exists $c\in (0,\infty)$ which satisfies for every $R\in (1,\infty)$, $M=(M_1,M_2)\in \N^2$, $L\in \N\backslash\{1\}$, $\ell=(\ell_0,\ell_1,\dots,\ell_L)\in (\{d+1\}\times\N^{L-1}\times\{1\}) $ and every random variable $\vartheta \colon \Omega \to [-R,R]^{\fd(\ell)}$ that
\begin{equation}\llabel{eq2}
\begin{split}
&\textstyle
  \E\bigl[
    \int_{[0,T] \times D}
      | u(y) - \reli{ \ell }{ L}{ \vartheta }( y ) |^2
    \, \d y
  \bigr]\textstyle\leq c^L(\fd(\ell))^{4\polycons} \bigl[R\textstyle\max\{\ell_1,\ell_2,\dots,\ell_{L-1}\}\bigr]^{2L+2}(M_1)^{-\nicefrac{1}{2}}\\
  &+c \bigl[\min\{\ell_1, \ell_2, \dots, \ell_{ L - 1 } \}\bigr]^{\frac{-2}{d+5}}+T\lambda( D)\,\E\bigl[\textstyle \bbF_{M,\ell}( \vartheta ) - \inf_{\theta\in (-c,c)^{\fd(\ell)}} \bbF_{M,\ell}( \theta ) \bigr]  \ifnocf.
  \end{split}
\end{equation}}
\startnewargseq
\argument{the fact that for all $L\in \N\backslash\{1\}$, $\ell=(\ell_0,\ell_1,\dots,\ell_L)\in (\{d+1\}\times\N^{L-1}\times\{1\}) $ it holds that $\fd(\ell)=\sum_{i=1}^L\ell_i(\ell_{i-1}+1)$}{that for all $L\in \N\backslash\{1\}$, $\ell=(\ell_0,\ell_1,\dots,\ell_L)\in (\{d+1\}\times\N^{L-1}\times\{1\}) $ it holds that
\begin{equation}\llabel{eq3}
\begin{split}
  \textstyle  \fd(\ell)&=\textstyle(d+1)\ell_1+\biggl[\sum\limits_{i=2}^{L-1}\ell_i(\ell_{i-1}+1)\biggr]+\ell_{L-1}+1\\
  &\leq (d+1)[\max\{\ell_1,\ell_2,\dots,\ell_{L-1}\}]^2+2(L-2)[\max\{\ell_1,\ell_2,\dots,\ell_{L-1}\}]^2\\
  &\quad+2[\max\{\ell_1,\ell_2,\dots,\ell_{L-1}\}]^2\\
  &\leq (2L+d)[\max\{\ell_1,\ell_2,\dots,\ell_{L-1}\}]^2.
  \end{split}
\end{equation}}
\argument{\lref{eq3};\lref{eq2}}{that for every $R\in (1,\infty)$, $M=(M_1,M_2)\in \N^2$, $L\in \N\backslash\{1\}$, $\ell=(\ell_0,\ell_1,\dots,\ell_L)\in (\{d+1\}\times\N^{L-1}\times\{1\}) $ and every random variable $\vartheta \colon \Omega \to [-R,R]^{\fd(\ell)}$ it holds that
\begin{align}
&\textstyle
  \E\bigl[
    \int_{[0,T] \times D}
      | u(y) - \reli{ \ell }{ L}{ \vartheta }( y ) |^2
    \, \d y
  \bigr]\textstyle\leq c^L(2L+d)^{4\polycons} [R\textstyle\max\{\ell_1,\ell_2,\dots,\ell_{L-1}\}]^{2L+8q+2}(M_1)^{-\nicefrac{1}{2}} \notag\\
  &+c [\min\{\ell_1, \ell_2, \dots, \ell_{ L - 1 } \}]^{\frac{-2}{d+5}}+T\lambda( D)\,\E\bigl[\textstyle \bbF_{M,\ell}( \vartheta ) - \inf_{\theta\in (-c,c)^{\fd(\ell)}} \bbF_{M,\ell}( \theta ) \bigr] \llabel{eq4}  \ifnocf.
\end{align}}
\argument{\lref{eq4};}{for every $R\in (1,\infty)$, $M=(M_1,M_2)\in \N^2$, $L\in \N\backslash\{1\}$, $\ell=(\ell_0,\ell_1,\dots,\ell_L)\in (\{d+1\}\times\N^{L-1}\times\{1\}) $ and every random variable $\vartheta \colon \Omega \to [-R,R]^{\fd(\ell)}$ that
\begin{align}
&\textstyle
  \E\bigl[
    \int_{[0,T] \times D}
      | u(y) - \reli{ \ell }{ L}{ \vartheta }( y ) |^2
    \, \d y
  \bigr]\textstyle\leq c^L(2L+d)^{4\polycons} [R+\textstyle\max\{\ell_1,\ell_2,\dots,\ell_{L-1}\}]^{4L+16q+4}(M_1)^{-\nicefrac{1}{2}}\notag\\
  &+c [\min\{\ell_1, \ell_2, \dots, \ell_{ L - 1 } \}]^{\frac{-2}{d+5}}+T\lambda(D)\,\E\bigl[\textstyle \bbF_{M,\ell}( \vartheta ) - \inf_{\theta\in (-c,c)^{\fd(\ell)}} \bbF_{M,\ell}( \theta ) \bigr] \llabel{eq5} \ifnocf.
\end{align}}
\argument{the fact that for all $x\in (0,\infty)$ it holds that $2^x\geq x$}{that there exist $\scrc_1,\scrc_2\in (0,\infty)$ such that for all $R\in (1,\infty)$, $L\in \N\backslash\{1\}$, $\ell=(\ell_0,\ell_1,\dots,\ell_L)\in (\{d+1\}\times\N^{L-1}\times\{1\}) $ it holds that
\begin{equation}\llabel{eq6}
    \begin{split}
        &c^L(2L+d)^{4\polycons} [R+\textstyle\max\{\ell_1,\ell_2,\dots,\ell_{L-1}\}]^{4L+16\polycons+4}\\
       &\leq 2^{Lc}2^{L\scrc_1}[R+\textstyle\max\{\ell_1,\ell_2,\dots,\ell_{L-1}\}]^{4L+16\polycons+4}\leq [R+\textstyle\max\{\ell_1,\ell_2,\dots,\ell_{L-1}\}]^{L\scrc_2}.
    \end{split}
\end{equation}}
\argument{\lref{eq6};\lref{eq5}}{\cref{conclude: main cor}\dott}
\end{aproof}
\subsection{Convergence results for the deep Kolmogorov method}

In this subsection we employ \emph{the error analysis with convergence rates} in \cref{cor: linear II relu pinn} above to conclude in \cref{main theorem 3} and \cref{main theorem 2} below \emph{convergence to zero} of the overall mean squared error of the \DKM\ without any information on the speed of convergence.

To conclude convergence to zero of the overall mean squared error, \cref{main theorem 3} and \cref{main theorem 2} both assume that the minimal width of the approximating \ANN\ converges to infinity, that the optimization error converges to zero, and that the number of random sample points converges sufficiently quickly to infinity.

\cref{main theorem 3} and \cref{main theorem 2} differ in the sense that \cref{main theorem 2} is a bit shorter but in \cref{main theorem 2} the number of random sample points is assumed to converge a bit quicker to infinity than in \cref{main theorem 3} (cf.\ \cref{def: bbF: main theorem 3}--\cref{conclude: main theorem 3} in \cref{main theorem 3} with \cref{def: bbF: main theorem 2}--\cref{conclude: main theorem 2} in \cref{main theorem 2}). \cref{main theorem 3} is established through an application of \cref{cor: linear II relu pinn} above. \cref{main theorem 2}, in turn, is proved through an application of \cref{main theorem 3}.

\begin{athm}{cor}{main theorem 3}
Let $ T \in (0,\infty) $, $ d \in \N $, let $u\in C^2([0,T]\times\R^{d},\R)$ be at most polynomially growing, assume for all $t\in (0,T)$, $x\in \R^d$ that
\begin{equation}\llabel{pde}
\textstyle
  \frac{ \partial }{ \partial t } u( t, x )
  +
  \frac{ 1}{ 2 }
  \Delta_x u(t,x)
  = 0,
  \end{equation}
  let $\polycons\in [1,\infty)$ satisfy $\sup_{x\in \R^d} ((1+\|x\|)^{-\polycons}|u(T,x)|)<\infty$,
 let $ ( \Omega, \mathcal{F}, \P ) $ be a probability space,
let $ \mathbb{T}_n \colon \Omega \to [0,T] $, $n\in \N$, be independent uniformly distributed random variables, let $D\subseteq\R^d$ be open and bounded,
let $ \mathbb{X}_n \colon \Omega \to D $, $n\in \N$, be independent uniformly distributed random variables,
let $ W^n \colon [0,T] \times \Omega \to \R^d $, $n\in \N$, be \iid\ standard Brownian motions,  
assume that $ (\mathbb{T}_n)_{n\in \N} $, $ (\mathbb{X}_n)_{n\in \N} $, and $ (W^n)_{n\in \N} $ are independent, for every $k\in \N$ let $M^k = ( M^k_1, M^k_2 ) \in \N^2$, $L_k\in \N\backslash\{1\}$, $\ell_k=(\ell_k^0,\ell_k^1,\dots,\ell_k^{L_k}) \in ( \{ d + 1 \} \times \N^{ L_k - 1 } \times \{ 1 \} )$, $R_k \in (1,\infty)$ and let $ \mathbb{F}_{k} \colon \R^{ \ffd( \ell_k ) } \to \R $
satisfy
for all
$ \theta \in \R^{ \ffd( \ell_k ) } $
that
\begin{equation}\label{def: bbF: main theorem 3}
\textstyle
  \mathbb{F}_{k}( \theta )=
  \frac {1}{M_1^k} \sum_{m=1}^{M_1^k}
    \big|
      \reli{ \ell_k }{ L }{ \theta }( \mathbb{T}_m, \mathbb{X}_m )
      - \bigl[ \frac{ 1 }{ M_2 ^k} \sum_{ n = 1 }^{ M_2^k } u(T, \bbX_{ m } + W^{ m M_2^k + n }_{ T-\bbT_m } ) \bigr]
    \big|^2,
\end{equation}
assume $\liminf_{ k \to \infty } \min\{ \ell_k^1, \ell_k^2,\allowbreak \dots,\allowbreak \ell_k^{ L_k-1 } \} = \infty$, assume $\limsup_{k\to\infty}((R_k)^{4L_k+4}[\textstyle\max\{\ell_k^1,\ell_k^2,\dots,
\allowbreak\ell_k^{L_k-1}\}]^{8L_k+16\polycons+4}(M_1^k)^{-1})=0$, and let $\vartheta_k\colon\Omega\allowbreak\to [-R_k,R_k]^{\fd(\ell_k)}$, $k\in \N$, be random variables which satisfy $\limsup_{ k \to \infty } \E\bigl[ \bbF_k( \vartheta_k ) - \inf_{ \theta \in \R^{ \ffd(\ell_k) } }\allowbreak \bbF_k( \theta ) \bigr] = 0$ \cfload.
Then 
\begin{equation}\label{conclude: main theorem 3}
\begin{split}
\textstyle
  \limsup_{k\to\infty}\E[
    \int_{[0,T]\times D}
      | u(y) - \reli{ \ell_k }{ L}{ \vartheta_k}( y ) |^2
    \, \d y
  ]=0  \ifnocf.
  \end{split}
\end{equation}
\cfout[.]
\end{athm}
\begin{aproof}
    \argument{\cref{cor: linear II relu pinn};}{that there exists $c\in (0,\infty)$ which satisfies for all $k\in \N$ that
    \begin{align}
&\textstyle
  \E[
    \int_{[0,T]\times D}
      | u(y) - \reli{ \ell_k }{ L}{ \vartheta_k }( y ) |^2
    \, \d y
  ]\textstyle\leq c^{L_k}(\fd(\ell_k))^{4\polycons} [R_k\textstyle\max\{\ell_k^1,\ell_k^2,\dots,\ell_k^{L_k-1}\}]^{2L_k+2}(M_1^k)^{-\nicefrac{1}{2}}\notag\\
  &+c [\min\{\ell_k^1, \ell_k^2, \dots, \ell_k^{ L_k - 1 } \}]^{\frac{-2}{d+5}}+T\lambda( D)\,\E\bigl[\textstyle \bbF_{k}( \vartheta_k ) - \inf_{\theta\in \R^{\fd(\ell_k)}} \bbF_{k}( \theta ) \bigr]  \llabel{eq1}\ifnocf.
\end{align}}
\startnewargseq
\argument{the fact that for all $L\in \N\backslash\{1\}$, $\ell=(\ell_0,\ell_1,\dots,\ell_L)\in \N^{L+1} $ it holds that $\fd(\ell)=\sum_{i=1}^L\ell_i(\ell_{i-1}+1)$}{that for all $k\in \N$ it holds that
\begin{equation}\llabel{eq2}
\begin{split}
  \textstyle  \fd(\ell_k)&=\textstyle(d+1)\ell_k^1+\biggl[\sum\limits_{i=2}^{L_k-1}\ell_k^i(\ell_k^{i-1}+1)\biggr]+\ell_k^{L_k-1}+1\\
  &\leq (d+1)[\max\{\ell_k^1,\ell_k^2,\dots,\ell_k^{L_k-1}\}]^2+2(L_k-2)[\max\{\ell_k^1,\ell_k^2,\dots,\ell_k^{L_k-1}\}]^2\\
  &+2[\max\{\ell_k^1,\ell_k^2,\dots,\ell_k^{L_k-1}\}]^2\\
  &\leq (2L_k+d)[\max\{\ell_k^1,\ell_k^2,\dots,\ell_k^{L_k-1}\}]^2.
  \end{split}
\end{equation}}
\argument{\lref{eq2};\lref{eq1}}{for all $k\in \N$ that
\begin{equation}\llabel{eq3}
\begin{split}
&\textstyle
  \E[
    \int_{[0,T]\times D}
      | u(y) - \reli{ \ell _k}{ L}{ \vartheta_k }( y ) |^2
    \, \d y
  ]\\
  &\textstyle\leq c^{L_k}(2L_k+d)^{4\polycons} (R_k)^{2L_k+2}[\textstyle\max\{\ell_k^1,\ell_k^2,\dots,\ell_k^{L_k-1}\}]^{2L_k+8q+2}(M_1^k)^{-\nicefrac{1}{2}}\\
  &+c [\min\{\ell_k^1, \ell_k^2, \dots, \ell_k^{ L_k - 1 } \}]^{\frac{-2}{d+5}}+T\lambda( D)\,\E\bigl[\textstyle \bbF_{k}( \vartheta_k ) - \inf_{\theta\in \R^{\fd(\ell_k)}} \bbF_{k}( \theta ) \bigr]\\
  &\leq (2d)^{8q}c^{L_k}(L_k)^{4\polycons} (R_k)^{2L_k+2}[\textstyle\max\{\ell_k^1,\ell_k^2,\dots,\ell_k^{L_k-1}\}]^{2L_k+8\polycons+2}(M_1^k)^{-\nicefrac{1}{2}}\\
  &+c [\min\{\ell_k^1, \ell_k^2, \dots, \ell_k^{ L_k - 1 } \}]^{\frac{-2}{d+5}}+T\lambda( D)\,\E\bigl[\textstyle \bbF_{k}( \vartheta_k ) - \inf_{\theta\in \R^{\fd(\ell_k)}} \bbF_{k}( \theta ) \bigr]\ifnocf.
  \end{split}
\end{equation}}
\argument{\lref{eq3};the fact that $\sup_{n\in \N} 2^{-n}n^{4q}<\infty$}{that there exists $\scrc\in (0,\infty)$ such that for all $k\in \N$ with $\max\{\ell_k^1,\ell_k^2,\dots,\ell_k^{L_k-1}\}\geq \max\{c,2\}$ it holds that
\begin{samepage}
\begin{align}
&\textstyle
  \E[
    \int_{[0,T]\times D}
      | u(y) - \reli{ \ell_k }{ L}{ \vartheta_k }( y ) |^2
    \, \d y
  ] 
  \leq \scrc(R_k)^{2L_k+2}[\textstyle\max\{\ell_k^1,\ell_k^2,\dots,\ell_k^{L_k-1}\}]^{4L_k+8\polycons+2}(M_1^k)^{-\nicefrac{1}{2}}\notag\\
  &+c [\min\{\ell_k^1, \ell_k^2, \dots, \ell_k^{ L_k - 1 } \}]^{\frac{-2}{d+5}}+T\lambda(D)\,\E\bigl[\textstyle \bbF_{k}( \vartheta_k ) - \inf_{\theta\in \R^{\fd(\ell_k)}} \bbF_{k}( \theta ) \bigr]\llabel{eq4}.
\end{align}
\end{samepage}}
\argument{\lref{eq4};the assumption that 
\begin{equation}
    \liminf_{ k \to \infty } \min\{ \ell_k^1, \ell_k^2,\allowbreak \dots,\allowbreak \ell_k^{ L_k-1 } \} = \infty;
\end{equation}
the assumption that $\limsup_{k\to\infty}((R_k)^{4L_k+4}[\textstyle\max\{\ell_k^1,\ell_k^2,\dots,
\allowbreak\ell_k^{L_k-1}\}]^{8L_k+16\polycons+4}(M_1^k)^{-1})=0$ ; the assumption that $\limsup_{ k \to \infty }\allowbreak \E\bigl[ \bbF_k( \vartheta_k ) \allowbreak- \inf_{ \theta \in \R^{ \ffd(\ell_k) } }\allowbreak \bbF_k( \theta ) \bigr] = 0$}{\cref{conclude: main theorem 3}\dott}
\end{aproof}
\begin{athm}{cor}{main theorem 2}
Let $ T \in (0,\infty) $, $ d \in \N $,
let $ u \in C^{2}( [0,T] \times \R^d, \R ) $ be at most polynomially growing, assume for all $ t \in (0,T) $, $ x \in \R^d $ that
\begin{equation}\llabel{pde}
\textstyle
  \frac{ \partial }{ \partial t } u( t, x )
  +
  \frac{ 1 }{ 2 }
  \Delta_x u(t,x)
  = 0,
\end{equation}
 let $ ( \Omega, \mathcal{F}, \P ) $ be a probability space,
let $ \mathbb{T}_n \colon \Omega \to [0,T] $, $n\in \N$, be independent uniformly distributed random variables, let $D\subseteq\R^d$ be open and bounded,
let $ \mathbb{X}_n \colon \Omega \to D $, $n\in \N$, be independent uniformly distributed random variables,
let $ W^n \colon [0,T] \times \Omega \to \R^d $, $n\in \N$, be \iid\ standard Brownian motions, 
assume that $ (\mathbb{T}_n)_{n\in \N} $, $ (\mathbb{X}_n)_{n\in \N} $, and $ (W^n)_{n\in \N} $ are independent, for every $k\in \N$ let $M^k = ( M^k_1, M^k_2 ) \in \N^2$, $L_k\in \N\backslash\{1\}$, $\ell_k=(\ell_k^0,\ell_k^1,\dots,\ell_k^{L_k}) \in ( \{ d + 1 \} \times \N^{ L_k - 1 } \times \{ 1 \} )$, $R_k \in (1,\infty)$ and let $ \mathbb{F}_{k} \colon \R^{ \ffd( \ell_k ) } \to \R $
satisfy
for all
$ \theta \in \R^{ \ffd( \ell_k ) } $
that
\begin{equation}\label{def: bbF: main theorem 2}
\textstyle
  \mathbb{F}_{k}( \theta )=
  \frac {1}{M_1^k} \sum_{m=1}^{M_1^k}
    \big|
      \reli{ \ell_k }{ L }{ \theta }( \mathbb{T}_m, \mathbb{X}_m )
      - \bigl[ \frac{ 1 }{ M_2 ^k} \sum_{ n = 1 }^{ M_2^k } u(T, \bbX_{ m } + W^{ m M_2^k + n }_{ T-\bbT_m } ) \bigr]
    \big|^2,
\end{equation}
assume $\liminf_{ k \to \infty } \min\{ \ell_k^1, \ell_k^2,\allowbreak \dots,\allowbreak \ell_k^{ L_k-1 } \} = \infty$, assume $\sup_{ c \in [1,\infty) }\limsup_{k\to\infty}([R_k+\textstyle\max\{\ell_k^1,\ell_k^2,\allowbreak\dots,\ell_k^{L_k-1}\allowbreak\}]^{cL_k}(M_1^k)^{-1})=0$, and let $\vartheta_k\colon\Omega\to [-R_k,R_k]^{\fd(\ell_k)}$, $k\in \N$, be random variables which satisfy $\limsup_{ k \to \infty } \E\bigl[ \bbF_k( \vartheta_k ) - \inf_{ \theta \in \R^{ \ffd(\ell_k) } }\allowbreak \bbF_k( \theta ) \bigr] = 0$ \cfload.
Then 
\begin{equation}\label{conclude: main theorem 2}
\begin{split}
\textstyle
  \limsup_{k\to\infty}\E\bigl[
    \int_{[0,T]\times D}
      | u(y) - \reli{ \ell_k }{ L}{ \vartheta_k}( y ) |^2
    \, \d y
  \bigr]=0  \ifnocf.
  \end{split}
\end{equation}
\cfout[.]
\end{athm}
\begin{aproof}
\argument{the assumption that $u$ is at most polynomially growing;}{that there exists $\polycons\in [1,\infty)$ which satisfies
\begin{equation}\llabel{def: q}
    \textstyle\sup_{x\in \R^d}((1+\|x\|)^{-q}|u(T,x)|)<\infty.
\end{equation}}
\startnewargseq
\argument{the Cauchy-Schwarz inequality; the fact that for all $k\in \N$ it holds that $R_k>0$}{that for all $k\in \N$ it holds that
\begin{equation}\llabel{eq0}
    R_k+\max\{\ell_k^1,\ell_k^2,\allowbreak\dots,\ell_k^{L_k-1}\allowbreak\}\geq [R_k\max\{\ell_k^1,\ell_k^2,\allowbreak\dots,\ell_k^{L_k-1}\allowbreak\}]^{1/2}.
\end{equation}}
\argument{\lref{eq0};the assumption that $\sup_{ c \in [1,\infty) }\limsup_{k\to\infty}([R_k+\textstyle\max\{\ell_k^1,\ell_k^2,\allowbreak\dots,\ell_k^{L_k-1}\allowbreak\}]^{cL_k}(M_1^k)^{-1})=0$;}{that 
\begin{equation}\llabel{eq1}
    \limsup_{k\to\infty}([R_k\textstyle\max\{\ell_k^1,\ell_k^2,\dots,\ell_k^{L_k-1}\}]^{(16\polycons+12)L_k}(M_1^k)^{-1})=0.
\end{equation}}
\argument{\lref{eq1};the assumption that for all $k\in \N$ it holds that $R_k\geq 1$; the assumption that for all $k\in \N$ it holds that $L_k\geq 1$}{that
\begin{equation}\llabel{eq2}
    \limsup_{k\to\infty}((R_k)^{4L_k+4}[\textstyle\max\{\ell_k^1,\ell_k^2,\dots,\ell_k^{L_k-1}\}]^{8L_k+16\polycons+4}(M_1^k)^{-1})=0.
\end{equation}}
    \argument{\lref{eq2};\cref{main theorem 3};}{\cref{conclude: main theorem 2}\dott}
\end{aproof}
\subsubsection*{Acknowledgements}
This work has been partially funded by the European Union (ERC, MONTECARLO, 101045811). The views and the opinions expressed in this work are however those of the authors only and do not necessarily reflect those of the European Union or the European Research Council (ERC). Neither the European
Union nor the granting authority can be held responsible for them. We also gratefully acknowledge the Cluster of Excellence EXC 2044-390685587, Mathematics Münster: Dynamics-Geometry-Structure funded by the Deutsche Forschungsgemeinschaft (DFG, German Research Foundation). Most of the specific formulations in the proofs of this work have been created using \cite{Bennoargumentcommand}. Johannes Schmidt-Hieber is gratefully acknowledged for several fruitful discussions.
\bibliographystyle{acm}
\bibliography{bibfile}

\end{document}